\RequirePackage{snapshot} 
\documentclass[11pt]{article}
\pdfoutput=1 
\usepackage{geometry}
\usepackage{snapshot}
\usepackage{amssymb} 
\usepackage{amsmath} 
\usepackage{amsthm}
\usepackage{stmaryrd} 
\usepackage{comment} 
\usepackage{color}
\usepackage{graphicx}
\usepackage{subfig}
\usepackage{tikz}
\usepackage{pgfplots}
\usepackage{pgfplotstable}
\usepackage{algorithm}
\usepackage{algpseudocode}
\usepackage{booktabs}

\usepackage[colorinlistoftodos,prependcaption]{todonotes}

\newcommand{\icol}[1]{
  \left(\begin{smallmatrix}#1\end{smallmatrix}\right)%
}

\newtheorem*{remark}{Remark}

\usepackage{environ} 

\NewEnviron{CommentsFEM}
{ 
}

\newcommand{\jump}[1]{\left\llbracket #1\right\rrbracket}
\newcommand{\norm}[1]{\left\| #1\right\|}

\newcommand{\V}[1]{\underline{#1}}
\newcommand{\M}[1]{\underline{\underline{#1}}}

\newcommand{\Cola}{  \color{black}  }
\newcommand{\Colb}{  \color{black}  }

\newcommand{\Grad}{{ \Cola \M{\text{grad}} }}

\newcommand{\Domain}{{ \Cola \Omega  }}
\newcommand{\DomainHat}{{ \Cola \widehat \Omega  }}
\newcommand{\Interface}{{ \Cola \Gamma  }}

\newcommand{\DDomain}{{ \Cola \partial \Omega }}

\newcommand{\Normal}{{ \Cola \V{n} }}
\newcommand{\DDomainDirichlet}{{ \Cola \partial \Omega^\textrm{u} }}

\newcommand{\DDomainNeumann}{{ \Cola \partial \Omega^\textrm{t} }}

\newcommand{\DDomainNeumannFEM}{{ \Cola  \partial   \Omega ^{\textrm{t}}_\textrm{h} }}

\newcommand{\NeumannBC}{{ \Cola \V{T} }}
\newcommand{\DirichletBC}{{ \Cola \V{U} }}
\newcommand{\BodyForce}{{ \Cola \V{f} }}

\newcommand{\Disp}{{ \Cola  \V{u} }}

\newcommand{\DispHatTrial}{{ \Cola  \V{u}^{*,i}_h }}

\newcommand{\Stress}{{ \Cola \M{\sigma} }}

\newcommand{\DispTest}{{ \Cola \delta \Disp }}

\newcommand{\StrainBasic}{{ \Cola \M{\epsilon} }}

\newcommand{\InterfaceForce}{{ \Colb \V{F}  }}
\newcommand{\InterfaceDisp}{{ \Colb \V{W}  }}

\newcommand{\InterfaceForceHat}{{ \Colb \widehat{ \V{F} }  }}
\newcommand{\InterfaceDispHat}{{ \Colb  \widehat{ \V{W} } }}

\newcommand{\Hooke}{{ \Cola \underset{\widetilde{}}{D} }}

\newcommand{\NumberSubDom}{n_\text{d}}

\newcommand{\IdentityTensor}{{ \Cola \M{I}_d }}

\newcommand{\bmesh}{\mathcal{T}_h}

\title{A stable and optimally convergent LaTIn-Cut Finite Element Method for multiple unilateral contact problems
}

\author{S. Claus$^{1}$\footnote{clauss2@cardiff.ac.uk} , P. Kerfriden$^{1}$
\\ \\
$^{1}$ Cardiff University, School of Engineering, \\ The Parade, CF243AA Cardiff, United Kingdom 
}
\date{}
  
\begin{document}

\maketitle

\begin{abstract}
In this paper, we propose a novel unfitted finite element method for the simulation of multiple body contact. The computational mesh is generated independently of the geometry of the interacting solids, which can be arbitrarily complex. The key novelty of the approach is the combination of elements of the CutFEM technology, namely the enrichment of the solution field via the definition of overlapping fictitious domains with a dedicated penalty-type regularisation of discrete operators, and the LaTIn hybrid-mixed formulation of complex interface conditions. Furthermore, the novel P1-P1 discretisation scheme that we propose for the unfitted LaTIn solver is shown to be stable, robust and optimally convergent with mesh refinement. Finally, the paper introduces a high-performance 3D level-set/CutFEM framework for the versatile and robust solution of contact problems involving multiple bodies of complex geometries, with more than two bodies interacting at a single point.\\

\noindent \textbf{keywords:} unilateral contact, LaTIn, nonconforming finite element, CutFEM, ghost penalty, multiple level sets, composite materials
\end{abstract}



\section{Introduction}


Unfitted or non-conforming finite element methods uncouple the description of the geometry from the representation of the solution field itself. Typically, the geometry of the computational domain is projected over a regular background grid. In this setting, boundaries or interfaces between objects cut through elements of the corresponding mesh. Non-conforming methods are attractive for applications where coupling analysis codes and third party meshing libraries is either impractical, numerically expensive and/or prone to errors. Nonetheless, a number of specific challenges needs to be addressed for unfitted numerical solutions to be computable. Firstly, integrals need to be calculated over cut elements, which requires specialised numerical quadratures. Secondly, unfitted approaches require stabilisation. This is because combinations of degrees of freedom may be poorly controlled in regions where contributions from cut elements to integral forms are small. Finally, in the context of multiple interacting materials, enrichment of the finite element solution space is required in order to allow for numerical jumps or kinks to develop over embedded interfaces. Failure to do so may severely impair the convergence rate of the  unfitted finite element solver.

Over the last two decades, several encompassing frameworks have been developed to provide guidance for the development of unfitted finite element solvers. The eXtended Finite Element Method (XFEM) \cite{moesdolbow1999,FriesBelytschko2006,bordasmoran2006,guidaultallix2008} relies on the Partition of Unity Method \cite{melenkbabuska1996} to enrich the approximation space. Within this popular framework, non-conforming discontinuities can be introduced in an elegant and robust manner. Field discontinuities can also be enabled by making use of Embedded Discontinuity Methods \cite{olivercervera1999,moslermeschke2003,linderarmero2007}, whereby discontinuous shape functions are directly defined at the element level. The CutFEM framework \cite{hansbohansbo2002,burmanclaus2015}, which is fundamental to the developments reported in this paper, uses an alternative concept to introduce non-conforming field discontinuities, and to perform computations involving non-conforming geometrical domains. In the context of multiple interacting materials, CutFEM solvers enrich the solution space by overlaying (\textit{i.e.} doubling or tripling) cut elements. Subsequently, penalty techniques such as the Nitsche method are employed to disallow unwanted discontinuities, for instance to represent jumps in fluxes whilst preserving the continuity of the primal field \cite{hansbohansbo2002}. CutFEM solvers are usually equipped with a ghost-penalty stabilisation technique \cite{Burman_2010_a}, which results in optimal convergence rates whilst ensuring that the condition number of the system matrix degrades at the same rate as that of the conforming finite element method when the mesh is refined.

Extending unfitted finite element methods to nonlinear problems, and especially to problems involving complex interface behaviour, remains a challenge that needs to be addressed on a case-by-case basis. The present paper focuses on problems involving unilateral contact between elastic solids. This class of problems is important to engineers, and they are typically encountered when modelling joints between components in large-scale engineering assemblies. Numerically, unilateral contact is notoriously difficulty to treat, even in the context of fitted finite element approaches. References on the matter include \cite{wriggers2002,laursen2003,simolaursen1992,alartcurnier1991,mcdevittlaursen2000}. 
Of particular interest to us are LaTIn-based solvers \cite{champaneycognard1999,ladevezenouy2002a,boucardchampaney2003,liuborja2008,ladevezepassieux2009,allixkerfriden2010}, which treat unilateral contact through a two-field formulation of the interface equations, and resolve these fields iteratively using a dedicated two-stage solver. A number of researchers have demonstrated the versatility of LaTIn solvers, owing to the fact that the predictor stage (so-called ``linear stage"), is weakly dependent on the type of interface conditions, and completely independent of the current interface state. As a result, switching between nonlinear interface conditions is relatively painless \cite{kerfridenallix2009}. Moreover, dedicated space-time Model Order Reduction techniques may take advantage of the predictor invariance to reduce the cost of time-dependent or multi-query computations \cite{ladevezepassieux2009,kerfridenallix2009,giacomadureisseix2014}. However, a noteworthy issue related to the hybrid-mixed LaTIn formulation is the difficulty of discretising displacement and traction interface fields in a consistent manner. Whilst a piecewise constant discretisation of interface fields seems to be physically sound when associated with P1 finite element spaces for the bulk equations, such a discretisation scheme allows for the onset of numerically uncontrolled deformation modes whose associated amplitude creep up with increasing LaTIn iteration count, and eventually alter the convergence rate of the finite element solver. A claimed remedy to this numerical plague is to introduce additional degrees of freedom for the bulk fields, by means of either h-refinement or p-refinement. An alternative to this numerically expensive procedure has recently been proposed in \cite{gravouilpierres2011}, as discussed later on. Although both approaches seem to regularise the contact problem efficiently, the effect on the convergence of the mixed finite element solver with mesh refinement has not been thoroughly studied.

Over the years, several attempts to develop stable unfitted finite element solvers for unilateral contact have been met with success. XFEM-based solutions have been proposed in \cite{dolbowmoes2001,khoeinikbakht2006,ribeaucourtbaiettodubourg2007,elguedjgravouil2007,liuborja2008,gravouilpierres2011,muellerhoeppe2012}. Of particular importance to us are LaTIn-based unfitted solvers, as the LaTIn is also a fundamental technological element for the present contribution. The first of these solvers was proposed in \cite{dolbowmoes2001} to solve crack propagation problems (\textit{i.e.} the historically favoured playground for XFEM). In \cite{ribeaucourtbaiettodubourg2007}, the XFEM-LaTIn idea was extended in order to simulate fatigue crack propagation. In \cite{liuborja2008}, the authors compared a LaTIn formulation to a penalty/Newton approach, and claimed the superiority of the latter in terms of convergence rate and stability. Finally, the authors of \cite{gravouilpierres2011} proposed a regularised, LaTIn-type formulation of contact for XFEM. The main idea of the approach is to relax the kinematic continuity condition between the bulk primal field and its interface trace and to enforce it weakly via a regularised mortar method.

In contrast to XFEM-based developments, the CutFEM framework is relatively underdeveloped in the context of unfitted unilateral contact. We mention the recent work described in \cite{choulyhild2013,choulyhild2015}, which extends the Nitsche-type treatment of interface conditions that is usually employed in CutFEM solvers to the case of unilateral contact. The resulting fitted finite element solver is proven to be optimally convergent, stable and numerically efficient. A very first nonconforming CutFEM strategy for unilateral contact is discussed in \cite{burmanhansbo2017}, whereby the contact conditions are treated by an augmented Lagrange multiplier formulation, and the discrete equations are stabilised by a ghost-penalty regularisation.


The present paper proposes a novel LaTIn-CutFEM solver for unilateral contact between multiple interacting solid bodies. The solid bodies will be described by linearised elasticity, which implies that the location of contact interfaces is known \textit{a priori}. Our algorithm relies on the traditional CutFEM enrichment strategy \cite{hansbohansbo2002}. However, contact conditions over nonconforming interfaces will be enforced via the LaTIn method. The resulting algorithm is a two-stage solver, whereby the linear stage (\textit{i.e.} linear predictor) consists in solving a series of independent linear problems for each of the solid bodies, whilst the local stage (\textit{i.e.} interface corrector) is nonlinear. The method naturally inherits the coarse-grain parallel characteristics of LaTIn-based domain decomposition methods \cite{ladevezedureisseix1998,ladevezenouy2002a,kerfridenallix2009}, together with the previously mentioned versatility of LaTIn solvers with regards to the type of interface conditions being dealt with (\textit{i.e.} unilateral contact with Coulomb friction \cite{boucardchampaney2003,champaneyboucard2008}, frictional sliding \cite{annavarapuhautefeuille2013,annavarapuhautefeuille2014}, cohesive fracture models with transition to unilateral frictionless contact \cite{kerfridenallix2009}).
 
Two algorithmic elements of the proposed solver require particular care. Firstly, the LaTIn-CutFEM scheme needs to be stabilised. We regularise each of the subproblems of the linear stage using a ghost penalty technique \cite{Burman_2010_a}. We will show that, when using the appropriate scaling for the ghost penalty terms, the condition number of the system matrices ``degrades optimally" with mesh refinement (\textit{i.e.} it degrades at the same rate as that of the conforming finite element method). 
Secondly, we need to be careful when discretising the contact interface fields. As mentioned before, a naive discretisation of the LaTIn hybrid-mixed formulation leads to the onset of spurious interface solutions. We will show that the P1-P1 discretisation scheme that we propose (\textit{i.e.} P1 for the bulk primal field, P1 for both its trace and the trace of its normal flux) yields optimal convergence rates with mesh refinement, and eliminates the onset of instabilities. Our scheme ``nonlocalises" the local stage of the LaTIn, as it relies on a two-scale treatment of the interface fields and nonlinear equations. We will also show that the nonlinear solver of the LaTIn converges fast, in the sense that it quickly leads to iterates whereby the discretisation error is larger than the algorithmic error, typically within twenty iterations. To our best knowledge, these convergence results have never been reported in the context of LaTIn solvers, even in the context of conforming geometries.

Our formal developments are accompanied by a high-performance computer implementation. The core of our implementation is the finite element C++/Python library FEniCS \cite{fenics2015}, which, in particular, proposes a range of high-level tools to rapidly developed finite element solvers. The CutFEM C++ library, LibCutFEM, developed by S. Claus and A. Massing, and partially described in \cite{burmanclaus2015}, defines additional tools that are specific to unfitted finite elements. This library forms the basis for the LaTIn-CutFEM code. We will present results in two and three space dimensions, which is made seamless by the FEniCS framework. Some of our examples exhibit a high level of geometric complexity, which is handled through a robust implementation of multiple level sets. In particular, we will study the convergence of contact problems exhibiting multiple interacting bodies (\textit{i.e.} more than two) at a single point, (see other contributions on the topic, in a different context, in \cite{annavarapuhautefeuille2013,annavarapuhautefeuille2014}).

This paper is organised as follows. The LaTIn hybrid-mixed formulation of contact is presented in a continuous setting in Section 2, together with the associated two-stage nonlinear algorithm. The non-conforming regularised LaTIn-CutFEM formulation is presented in Section 3. Section 4 is dedicated to numerical investigations, where we  demonstrate the optimal convergence and stability properties of the proposed framework. We also showcase the versatility of our high-performance computer implementation.


\section{LaTIn method for unilateral contact: the continuous setting}
\label{sec:LaTin}



\subsection{Domain decomposition}

We consider a problem of linear elasticity over a 2 or 3-dimensional spatial domain $\Domain $ composed of $ \NumberSubDom $ non-overlapping subdomains $\Domain^{i}$, where $ i \in \mathcal{I}_\Domain :=  \{ 1, \dots, \NumberSubDom \} $. The domain decomposition satisfies $\bigcup_{ i \in \mathcal{I}_\Domain } \Domain^{i} = \Domain$. Let us denote the set of all subdomains by $\Xi_\Domain$. Subdomains may be multiply-connected. The boundary of subdomain $\Domain^i$ is denoted by $\partial\Domain^i$, for any $i \in \mathcal{I}_\Domain$, whilst $\DDomain$ denotes the boundary of computational domain $\Domain$. Furthermore, we define interface domains $\Interface^{i,j} := \DDomain^{i} \cap \DDomain^{j}$ for any $(i,j) \in \mathcal{I}_\Domain \times \mathcal{I}_\Domain$ such that the two indexes satisfy $i<j$ and the result of intersection $\DDomain^{i} \cap \DDomain^{j}$ is of non-zero measure. The union of these interfaces is denoted by $\Interface$. The set of pairs $(i,j)$ of indexes such that $\Interface^{i,j}$ is defined as a valid interface is denoted by $\mathcal{I}_{\Interface}$, and the corresponding set of interfaces is denotes by $\Xi_\Interface$.
We also denote the union of all the interfaces of subdomain $\Domain^{i}$ by $\Interface^{i} = \bigcup_{(i,j) \in \mathcal{I}_\Interface, (k=i \text{ or } l=i) }\Interface^{k,l}$, for any subdomain index $ i \in \mathcal{I}_\Domain$. We denote the outward normal pointing from $\Domain^i$ to $\Domain^j$ on $\Interface^{i,j}$ by $\Normal^{i,j}$. We further note that we have $\Normal^{i,j} = - \Normal^{j,i}$, and we denote the outward normal to $\Domain^i$ over $\Interface^{i}$ by $\Normal^i$. 
\begin{figure}[htb]
\centering
 \includegraphics[width=.25 \textwidth]{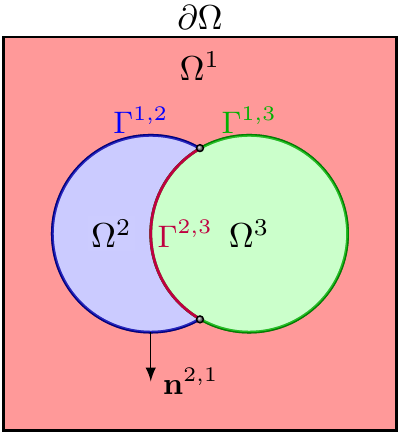}
 \caption{Schematic representation of a composite material made of three distinct phases.}
 \label{fig: schematic notation}
\end{figure}
As an illustration for this set of definitions, consider the case of the multiphase domain depicted in Figure~\ref{fig: schematic notation}. The cubic domain $\Domain$ comprises three subdomains $\Domain^{1}$, $\Domain^{2}$ and $\Domain^3$. There are three interfaces between the domains: $\Interface^{1,2}$, $\Interface^{1,3}$ and $\Interface^{2,3}$. Therefore, the interface of individual subdomains reads as  $\Interface^{1} = \Interface^{1,2} \cup \Interface^{1,3}$, $\Interface^{2} = \Interface^{1,2} \cup \Interface^{2,3}$ and $\Interface^{3} = \Interface^{1,3} \cup \Interface^{2,3}$. 

As further detailed in the next subsection, subdomain field variables will be governed by the equations of linear elasticity, while contact equations to be solved over interfaces will couple subdomain problems to one another. The LaTIn framework \cite{ladeveze1999} is particularly well suited to the treatment of such problems where complex nonlinear equations need to be solved over embedded interfaces \cite{boucardchampaney2003,ladevezepassieux2009,allixkerfriden2010b}. In this particular context, the LaTIn method relies on two main ingredients:
\begin{itemize}
\item a hybrid description of bulk and interface fields, and a mixed representation of the interface unknowns,
\item a two-step iterative solution algorithm for the solution of the global nonlinear problem formulated in its hybrid-mixed form.
\end{itemize}
Both these steps are described below.

\subsection{Hybrid-Mixed LaTIn formulation of unilateral contact problems}
\label{sec:Formulation}

\paragraph{Hybrid fields}

For every subdomain $\Domain^{i}$, where $i \in \mathcal{I}_\Domain$, we introduce a displacement field $\Disp^{i}:  \Domain^{i} \rightarrow \mathbb{R}^d$, where $d=2$ (2D problem) or $d=3$ (3D problem). \\
For every domain interface $\Interface^{i}$ with $i \in \mathcal{I}_\Domain$, we define two field variables:
\begin{itemize}
\item a force field $\InterfaceForce^{i}:  \Interface^{i} \rightarrow \mathbb{R}^d$ that represents a density of tractions (\textit{i.e.} flux) applied to subdomain $i$ by adjacent subdomains through interface $\Interface^{i}$.
\item a displacement field $\InterfaceDisp^{i}:  \Interface^{i} \rightarrow \mathbb{R}^d$ that represents the trace of $\Disp^{i}$ over $\Interface^{i}$.
\end{itemize}
As detailed below, $\Disp^{i}$ will be governed by linear elasticity, while over an interface $\Interface^{i,j}$ between two domains $\Domain^i$ and $\Domain^j$, where $(i,j) \in \mathcal{I}_\Interface$, interface models will provide a link between the restriction of fields $\InterfaceForce^{i}$, $\InterfaceDisp^{i}$ to $\Interface^{i,j}$ and the restriction of fields $\InterfaceForce^{j}$, $\InterfaceDisp^{j}$ to $\Interface^{i,j}$.

\paragraph{Bulk equations and boundary data} For every subdomain $\Domain^{i}$, $i \in \mathcal{I}_\Domain$, the displacement field $\Disp^{i} \in  [H^1(\Domain^i)]^d$ is required to be sufficiently regular (continuous in particular), and to satisfy the following equations:
\begin{itemize}
\item Static equilibrium equation and constitutive relation: \\
\textit{For all $ \DispTest^i \in [H^1(\Domain^i)]^d$, the displacement field $\Disp^i  \in  [H^1(\Domain^i)]^d$ satisfies} 
\begin{equation}
\int_ {\Domain^{i}} \Stress(\Disp^{i};\Hooke^i) : \StrainBasic(\DispTest^i)   \, d \Omega = \int_ {\Domain^{i}}   \BodyForce  \cdot   \DispTest^i  \, d \Omega +   \int_ {\DDomain^i}  \left( \Stress(\Disp^{i};\Hooke^i)   \cdot \Normal^{i} \right)  \cdot    \DispTest^i  \, d \Gamma \, .
\label{eq: linear elasticity}
\end{equation}
In equation \eqref{eq: linear elasticity}, $ \BodyForce$ is a body force and $\Normal^{i}$ is the outer normal of domain $\Domain^{i}$. The stress function is defined as $\Stress(\Disp^{i};\Hooke^i) := \Hooke^i : \StrainBasic(\Disp^{i})$, where the strain is the symmetric part of the displacement gradient, \textit{i.e.} $\StrainBasic(\Disp^{i}) :=  \frac{1}{2} \left ( \Grad \, \Disp^{i} +  (\Grad \, \Disp^{i} )^T \right)$ and  $\Hooke^i$ is a fourth-order Hooke tensor. In the case of isotropic elasticity, we have that  $\Hooke^i : \StrainBasic(\Disp^{i}) = \lambda^i \, \text{Tr}(\StrainBasic(\Disp^{i}))\, \IdentityTensor + 2 \, \mu^i \, \StrainBasic(\Disp^{i})$, where $\lambda^i$ and $\mu^i$ are the two Lam\'e coefficients. These two coefficients are related to the Poisson and Young's modulii, respectively denoted by $\nu$ and $E^i$, through 
$\lambda^i = \frac{E^i \nu}{(1 + \nu)(1 - 2 \nu)}$ and $\mu^i = \frac{E^i}{2(1 + \nu)} $. 
All these coefficients are assumed to be constant in $\Domain^i$. The dependency of the stress function to the Hooke tensor field will be omitted when possible. However, remember that the Hooke tensor may depend on the spatial coordinates. In particular, it may be discontinuous across an interface between connected subdomains. Such cases will be treated in the numerical example section.
  \item Boundary Conditions: \\
We shall assume that the boundary $\DDomain$ of domain $\Domain$ is decomposed into a Neumann part and a Dirichlet part, \it{i.e.} \normalfont $\DDomain = \DDomainNeumann  \cup \DDomainDirichlet$. 
For Neumann boundary conditions, we specify that
\begin{equation}
\Stress(\Disp^{i}) \cdot \Normal^{i} = \NeumannBC \qquad \text{over } \DDomainNeumann \cap \DDomain^{i} \, ,
\label{eq: neumann bc}
\end{equation}
whilst for Dirichlet boundary conditions, we set 
\begin{equation}
\Disp^{i} = \DirichletBC  \qquad \text{over }  \DDomainDirichlet   \cap \DDomain^{i} \, .
\label{eq: dirichlet bc}
\end{equation}
Here, $\NeumannBC$ is a density of prescribed tractions per unit of surface, $\DirichletBC$ is a field of prescribed displacement.
\end{itemize}

\paragraph{Interface governing equations} Any interface $\Interface^{i,j}$, where $(i,j) \in \mathcal{I}_\Interface$, will exhibit behaviour relating the surface tractions to the displacement fields. Although extensions are relatively straightforward, we consider only the case of frictionless unilateral contact, which is  characterised by the following equations:
\begin{itemize}
\item Newton's third law,
\begin{equation}
\label{eq:ThirdLaw}
 \InterfaceForce^{i} + \InterfaceForce^{j} = \V{0} \,,
\end{equation}
\item Signorini's law of unilateral contact, which reads as
\begin{equation}
\begin{aligned}
\left(  \InterfaceDisp^{j} - \InterfaceDisp^{i} \right) \cdot \Normal^{i} &\geq 0 \, , \\
 \InterfaceForce^{i} \cdot  \Normal^{i}   &\leq 0 \, , \\
\left(  (  \InterfaceDisp^{j} -  \InterfaceDisp^{i} ) \cdot \Normal^{i} \right) \cdot \left(   \InterfaceForce^{i}  \cdot\Normal^{i}  \right)  &= 0 \, ,
\end{aligned}
\label{eq:Signorini}
\end{equation}
\item the condition that the tangential component of tractions must vanish
\begin{equation}
(\IdentityTensor - \Normal^{i} \otimes \Normal^{i} ) \cdot \InterfaceForce^{i} = \V{0} \, .
\label{eq: traction contact}
\end{equation}
\end{itemize}

\paragraph{Interface compatibility conditions for the hybrid-mixed formulation} In order to close the system of equations, the interface fields must be linked to the bulk variables. This is done by enforcing the continuity condition
\begin{equation}
\InterfaceForce^{i} = \sigma(\Disp^{i}) \cdot \Normal^{i} \, ,
\label{eq: interface force compatibility}
\end{equation}
together with the kinematic condition
\begin{equation}
 \InterfaceDisp^{i} = \Disp^{i} \, 
 \label{eq: interface displacement compatibility}
\end{equation}
over all interfaces $\Interface^{i}$, with $i \in \mathcal{I}_\Domain$.

\subsection{Nonlinear iterative solver}

\begin{figure}[htb]
\centering
\includegraphics[width=0.45 \textwidth]{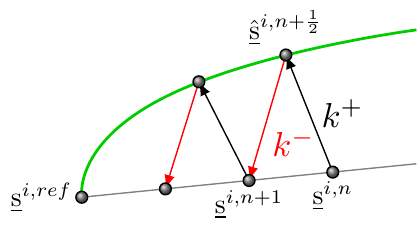}
\caption{Schematic of LaTIn algorithm.}
\label{fig:microstructure}
\end{figure}

We employ the standard LaTIn algorithm to solve the contact problem written in its hybrid-mixed form. The LaTIn algorithm consists in searching for mixed interface fields $(\InterfaceForce^{i},\InterfaceDisp^i)$ in two stages (see Figure \ref{fig:microstructure}). In the so-called linear stage, we search for interface fields that satisfy the bulk equations and the interface compatibility equations. The resulting problem is global, at least over each individual subdomain, but remains linear (similar to a predictor in a Newton algorithm). In the so-called local stage, we seek interface fields that satisfy the interface contact laws. The resulting problems are potentially non-linear (if contact is unilateral) but they are local  \textit{a priori} (this is similar to an nonlinear update in a Newton solver). As both primal and dual interface fields are unknown, search directions must be defined to close each of these two problems. Details about the algorithm are provided below.

\paragraph{Local Stage.} For all $\Interface^i$ with $i \in \mathcal{I}_\Domain$, we search for interface fields $\hat{\underbar{s}}^{i,n+\frac{1}{2}}:=(\InterfaceForceHat^{i,n+\frac{1}{2}}, \InterfaceDispHat^{i,n+\frac{1}{2}})$ satisfying unilateral contact equations \eqref{eq:ThirdLaw}-\eqref{eq: traction contact}, closed by the ascent search direction 
\begin{equation}
\left(\InterfaceForceHat^{i,n+\frac{1}{2}}  -  \InterfaceForce^{i,n}  \right) - k^{+} \left( \InterfaceDispHat^{i,n+\frac{1}{2}} - \InterfaceDisp^{i,n} \right) = \V{0} \, .
\end{equation}
Here, $n \in \mathbb{R}^{+}$ denotes the iteration index, $.^{n+\frac{1}{2}}$ denotes a quantity predicted by the $(n+1)-\text{th}$ local stage of the algorithm (\textit{i.e.} half iterate). The search direction  ``stiffness"  $k^+$ is a strictly positive scalar parameter.

\paragraph{Linear Stage.} For all $\Interface^i$ with $i \in \mathcal{I}_\Domain$, we search for interface fields $\underbar{s}^{i,n+1}:=(\InterfaceForce^{i,n+1},\InterfaceDisp^{i,n+1})$ satisfying the interface compatibility conditions \eqref{eq: interface force compatibility} - \eqref{eq: interface displacement compatibility} and the subdomain equations \eqref{eq: linear elasticity} - \eqref{eq: dirichlet bc}, closed by the descent search direction
\begin{equation}
\left( \InterfaceForce^{i,n+1} - \InterfaceForceHat^{i,n+\frac{1}{2}} \right) + k^{-} \left( \InterfaceDisp^{i,n+1} - \InterfaceDispHat^{i,n+\frac{1}{2}} \right) = \V{0} \, .
\end{equation}
We choose $k^- = k^+ $, which is the classical setting of the LaTIn algorithm (conjugate search directions \cite{ladeveze1999}). We further set $k^- = \xi \frac{E}{L}$, where $\xi$ is a non-dimensional strictly positive parameter (set to $1$ unless stated otherwise), $E$ is a reference Young's modulus (set to be the Young modulus of the matrix phase, \textit{i.e.} of the composite materials that we will analyse in all the examples of the numerical section), and $L$ is a reference length of the problem.

In a variational setting, and dropping iteration index $n$, the linear stage can be written as a set of linear problems that read: \textit{For any $i \in \mathcal{I}_\Domain$, 
find $\Disp^{i} \in [H^1(\Domain^i)]^d$ and compatible with the Dirichlet conditions such that for all $ \DispTest^i \in [H^1_0(\Domain^i)]^d$}
\begin{equation}
a_\text{D}^i(\Disp^{i},\DispTest^i) + a_\text{k}^i(\Disp^{i},\DispTest^i) = l_\text{f}^i(\DispTest^i) +  l_\text{k}^i(\DispTest^i)  \, ,
\label{eq:LinearStageContinuous}
\end{equation}
where the bilinear forms $a^i_{\text{D}}$ and $a_\text{k}^i$ are given by
\begin{equation}
a_\text{D}^{i}(\Disp^i,\delta \Disp^i) = \int_{\Domain^{i}} \StrainBasic(\Disp^i) : \Hooke^i : \StrainBasic(\DispTest^i) \, d\Domain
\end{equation}
\begin{equation}
a_\text{k}^{i}(\Disp^i,\delta \Disp^i) = \int_{\Interface^{i}} k^{-} \, \Disp^{i} \cdot \DispTest^i \, d\Interface \, ,
\end{equation}
and the linear forms $l_\text{f}^i$ and $l_\text{k}^i$  are given by
\begin{equation}
l_\text{f}^{i}(\delta \Disp^i) = \int_{\Domain^{i}} \BodyForce \cdot \DispTest^i \, d\Domain  + \int_{\partial \Domain^{i} \cap \DDomainNeumann}  \NeumannBC^i \cdot \DispTest^i \, d \Interface \, ,
\end{equation}
\begin{equation}
l_\text{k}^{i}(\delta \Disp^i) = \int_{ \Interface^{i} } \left( \InterfaceForceHat^{i} + k^{-} \, \InterfaceDispHat^{i} \right) \cdot  \DispTest^i   \, d \Interface \, .
\end{equation}
\begin{remark}
Problem \eqref{eq:LinearStageContinuous} is always well-posed (at least at the continuous level), owing to the effect of the LaTIn augmentation terms $a_\text{k}^{i}$, which ensure that the bilinear forms $a^{i}$ remain strictly coercive. 
\end{remark}

\paragraph{Convergence.} A relaxation step is performed at the end of the linear stage to control the convergence properties of the algorithm (the interested reader is referred to \cite{ladeveze1999} for proofs of convergence of the algorithm). The relaxation step reads as
\begin{equation}
\left\{ \begin{array}{l}
\displaystyle \InterfaceForce^{i,n+1} \leftarrow \eta \,  \InterfaceForce^{i,n+1}  + (1-\eta) \, \InterfaceForce^{i,n}, \\
\displaystyle \InterfaceDisp^{i,n+1} \leftarrow \eta \,  \InterfaceDisp^{i,n+1}  + (1-\eta) \, \InterfaceDisp^{i,n}
\end{array} \right.
\label{equ: relaxation linear stage}
\end{equation}
for all $i \in \mathcal{I}_\Domain$.
The robustness of the algorithm is not very sensitive to the choice of relaxation parameter $\eta \in [0 \, 1]$, but the convergence rate is. Based on our experience of the method, we set $\eta = 0.85$. 



\section{LaTin-CutFEM solver for embedded unilateral contact problems}



We present the proposed unfitted finite element strategy by firstly detailing the method that we use to describe the geometry of the subdomains. And secondly, we explain how the LaTIn hybrid-mixed formulation of unilateral contact problems are discretised and stabilised in the context of unfitted unilateral contact problems.

\subsection{Description and approximation of the geometry}
\label{subsection: fixed bg grids}

The geometry of the contact problems will be described independently of the background mesh using multiple level set functions whose zero level sets describe the interface locations (see similar applications of multiple level set approaches in \cite{merriman1994,smith2002,saye2012,starinshak2014}, amongst other contributions). Our general strategy is to describe each subdomain of our problem as the union of auxiliary subdomains that can be easily defined using level set functions. These functions can be given in an analytical form or in a numerical form. This approach is well suited to the description of contact problems with nonconforming embedded interfaces and/or nonconforming external boundaries. We will simply interpolate our level set functions in the finite element space, which will automatically yield an analysis suitable, piecewise polynomial approximation of the geometry.

\subsubsection{Formal domain decomposition using multiple level sets}

\begin{figure}
\centering
\subfloat[Level set description.]{\label{subfig: levelsetdescription}
\includegraphics[width=.6\textwidth]{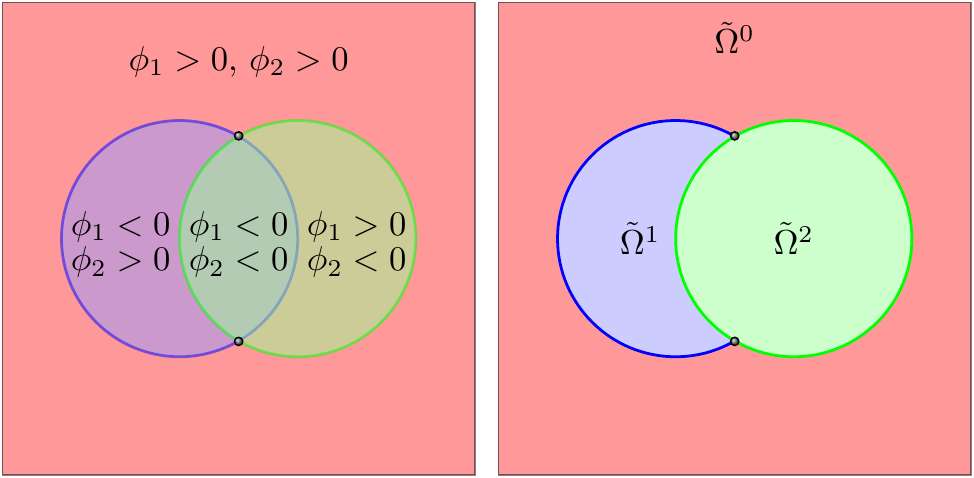}}
\subfloat[Piecewise linear approximation of the physical domains.]{\label{subfig: geometryapprox}\includegraphics[width=.3\textwidth]{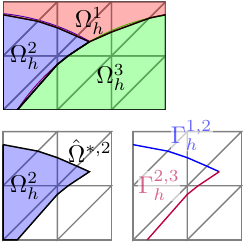}}
\caption{Schematic of multiple level set description of domain and its piecewise linear approximation.}
\label{fig: level set domains}
\end{figure}

Let $\Omega^{\square}$ be a domain that covers the physical domain $\Omega$, \textit{i.e.} $\Omega \subset \Omega^{\square}$, and let us decompose $\Omega^{\square}$ into $n_d^{\square}$ non-overlapping subdomains $\tilde{\Omega}^i$ with $n_d^{\square} \geq n_d$. $\Omega^{\square}$ will be meshed later on, and we therefore require that its geometry is of simple nature. 
Furthermore, let us assume that any physical subdomain $\Omega^i$ is exactly represented by one of the subdomains $\tilde{\Omega}^i$ or by the union of several of these auxiliary subdomains. 

The decomposition of $\Omega^{\square}$ into subdomains $\tilde{\Omega}^i$ is defined through a multiple level set approach. We employ  $n_d^{\square}-1$ level set functions $\phi^i: \mathbb{R}^d \rightarrow \mathbb{R}$, $i\in \{ 1, \dots n_d^{\square}-1 \}$ to describe the $n_d^{\square}$ subdomains  $\tilde{\Omega}^i$. This is done as follows. Let $I_{\phi^-}(\pmb{x})$ denote the set of indices of level set functions, which are smaller than zero in $\pmb{x}$, i.e. $I_{\phi^-}(\pmb{x})=\{ i \in \{1,\dots,n_d^{\square}-1 \}: \phi_i(\pmb{x}) <0 \}$. Then for each point $\pmb{x} \in \Omega^{\square}$, the domain that $\pmb{x}$ belongs to is defined as
\begin{equation}
\pmb{x} \in 
\begin{cases}
\tilde{\Omega}^0 &\mbox{if } I_{\phi^-} = \emptyset, \\
\tilde{\Omega}^{i} \mbox{, } i = \max(I_{\phi^-}) &\mbox{else}.
\end{cases}
\end{equation}
This definition implies that there exists a hierarchical relationship between level sets. Indeed, the domain defined by the region of negative level set values of level set $\phi^i$ overlaps all the domains defined by the negative regions of level sets $\phi^j$ with $j<i$. Also, domain $\tilde{\Domain}^0$ is defined as the complement to all other domains $\tilde{\Domain}^i$, i.e. $\tilde{\Domain}^0 := \Omega^{\square} \setminus \bigcup_{i \in I_{\tilde{\Omega}}} \tilde{\Domain}^i$. In our numerical examples, domain $\tilde{\Domain}^0$ represents the matrix phase of a composite material. Figure~\ref{subfig: levelsetdescription} shows a graphical representation of how the domains are described using multiple level sets. 

We now define the interface between $\tilde{\Omega}^i$ and $\tilde{\Omega}^j$, $\tilde{\Interface}^{i,j}$ ($i<j$), as
\begin{equation} 
\begin{aligned}
\tilde{\Interface}^{i,j}&:= \{\pmb{x} \in \Domain^{\square}: \phi^{j}(\pmb{x}) = 0 \mbox{ and } \phi^{i}(\pmb{x}) < 0 \} &\mbox{ for }  i \in \{1,\dots,n_d^{\square}-1 \}, \\
\tilde{\Interface}^{0,j}&:= \{\pmb{x} \in \Domain^{\square}: \phi^{j}(\pmb{x}) = 0  \mbox{ and } \phi^k(\pmb{x}) >0 \,\, \forall k \neq j \}. &
\end{aligned}
\end{equation}
Finally, we can formally map the subdomains of $\Omega^{\square}$ onto that of $\Omega$ with the mapping
\begin{equation}
\Omega^{i} = \mathcal{M}_d^{i} \left(\{ \tilde{\Omega}^{j}  \}_{  j \in \llbracket 0 , n_d^{\square}-1 \rrbracket} \right) =  
 \bigcup_{  j \in \tilde{\mathcal{I}}^{i} \subset \llbracket 0 , n_d^{\square}-1 \rrbracket, 
  } \tilde{\Omega}^{j}.
\end{equation}
This mapping means that each of the physical domains $\Omega^{i}$ is defined as the union of particular subsets of subdomains of $\{ \tilde{\Omega}^{i} \} $. Importantly, the union of subdomains $\tilde{\Omega}^{i}$ that are included in $\bar{\Omega} := \Omega^{\square} \backslash \Omega$ are void, \textit{i.e.} there is no partial differential equation written over $\bar{\Omega}$. In fact, $\bar{\Omega} \neq \emptyset$ means that the boundary of the domain will be treated in an implicit way, while $\bar{\Omega} = \emptyset$ means that the boundary $\partial \Omega$ is meshed explicitly.

The interface sets $\tilde{\Interface}^{i}$ and $\tilde{\Interface}^{i,j}$ are completely described by the definition of their counterpart  $\Interface^{i}$ and $\Interface^{i,j}$ and the definition of the mapping of $\tilde{\Domain}^{i}$. However, it should be pointed out that interfaces between subdomains of $\bar{\Omega}$ and subdomains of $\Omega^\square$ that also belong to $\Omega$ are implicitly defined boundaries. For these interfaces, the mapping can be written as
\begin{equation}
{\Interface}^{i,j} =  \bigcup_{k,l  \ | \ \tilde{\Omega}^{k} \in \left( \mathcal{M}_d^{k} \right)^{-1}(\Omega^{i} ) , \, \tilde{\Omega}^{l} \in  \left( \mathcal{M}_d^{l} \right)^{-1}(\Omega^{j} )}
\tilde{\Interface}^{k,l},
\end{equation}
\begin{align}
\partial \Omega &= \bigcup_{k, l \ | 
\begin{subarray}{c}
\tilde{\Omega}^k  \in \bar{\Omega},   \tilde{\Omega}^{l} \in  \left( \mathcal{M}_d^{l} \right)^{-1}(\Omega^{j} ) \\\text{ or } \tilde{\Omega}^l \in \bar{\Omega},  \tilde{\Omega}^{k} \in \left( \mathcal{M}_d^{k} \right)^{-1}(\Omega^{i} )
\end{subarray}
}\tilde{\Interface}^{k,l},
\end{align}
\textit{i.e.} any interface between two physical domains is a physical interface and anything between a physical domain and a "void" domain is a physical boundary.

This apparently complex formalism allows us to easily and generically define our approximate geometry by simply discretising all our level sets, as explained in the following.


%
%

\subsubsection{Discretisation of the level set functions and approximate geometry}
 
In this section, we describe the approximation of the physical subdomains $\Omega^i$ and their boundaries and interfaces.  Let $\bmesh$ denote a tessellation of our background domain $\Omega^{\square}$. Now, lets consider a piecewise linear approximation of our level set functions $\phi_h^i$ resulting in  a piecewise linear (in 2D) or piecewise planar (in 3D) interface approximation $\tilde{\Interface}^{i,j}_h$. This yields piecewise linear approximations of $\Omega^i$, $\Omega_h^i$,  and  $\Interface^{i,j}$, $\Interface^{i,j}_h$, as illustrated in Figure~\ref{subfig: geometryapprox}.

To ensure that the background mesh is fine enough to resolve the interfaces sufficiently, we impose the condition on the background mesh that each $\Interface^{i,j}_h$ intersects any face of each element $K \in G^i$ at most once and that the mesh is quasi-uniform. 

Note that the integrals of the discretised weak formulation introduced in the next section are evaluated over $\Domain_h^i$ and $\Gamma^{i,j}_h$, which includes the evaluation of integrals over arbitrary shaped areas in  interface regions and integrals over linear/planar interface parts. These integrals over the areas $K_i = K \cap \Domain_h^i$ are performed using quadrature rules which are generated from a standard sub-triangulation as detailed in \textit{e.g.} \cite{burmanclaus2015}. 

\subsection{Mixed Finite element formulation}

Using the piecewise polynomial description of our physical subdomains, we can now formulate our LaTIn-CutFEM method. The key features, which make our method stable and optimally convergent are: 
\begin{itemize}
\item a representation of both primal and dual interface quantities in the trace of the finite element space defined over the subdomains. That means, we develop a stable hybrid P1-P1 scheme, \textit{i.e.}  we employ a piecewise linear approximation of the displacement field in the bulk of the subdomains, and piecewise linear approximations of the displacement and force fields over interfaces between subdomains;
\item the regularisation of the equations of elasticity over each of the subdomains using the ghost-penalty methodology;
\item a weakened enforcement of the contact equation through a stabilised two-scale representation of the interface fields.
\end{itemize}

\subsubsection{Bulk and interface finite element spaces}

\paragraph{Fictitious bulk finite element spaces.}

\begin{figure}[htb]
\centering
\includegraphics[width=.9\textwidth]{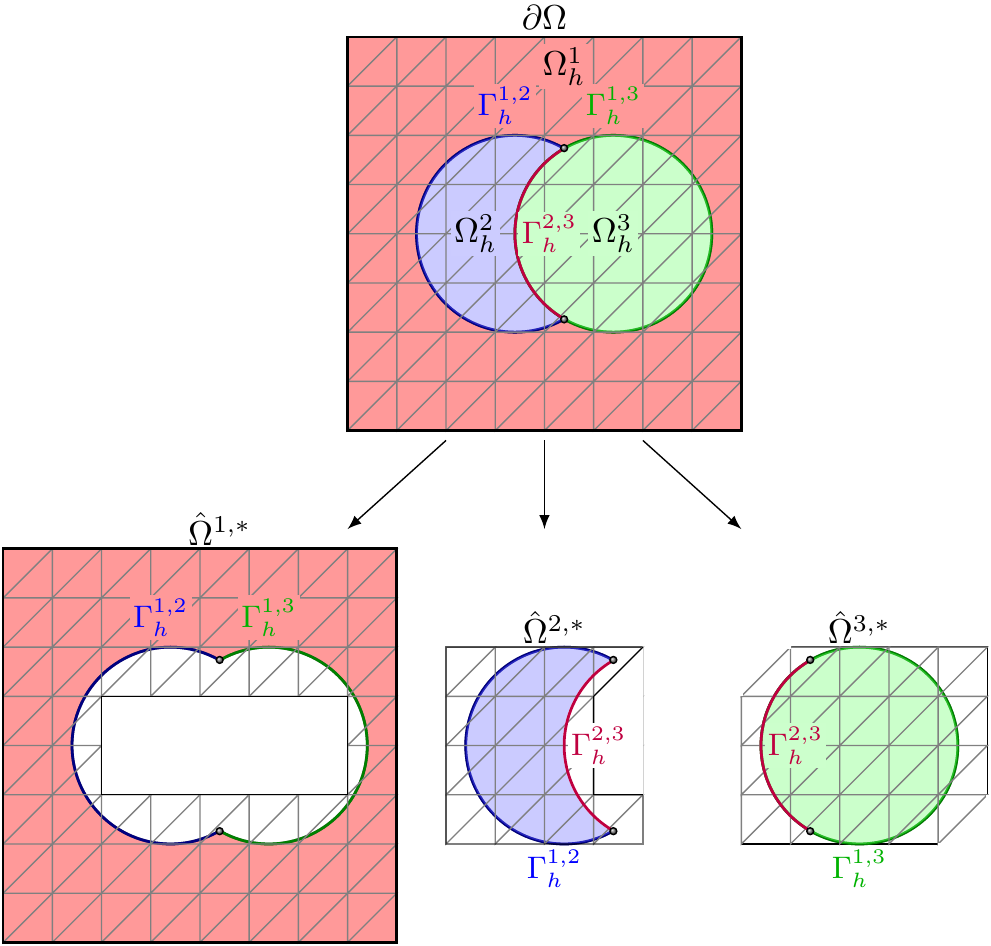}
\caption{Schematic representation of fictitious domain meshes $\hat{\Domain}^{i,*}$.}
\label{fig: fictitious domain meshes}
\end{figure}

\begin{figure}[htb]
\centering
\includegraphics[width=.4\textwidth]{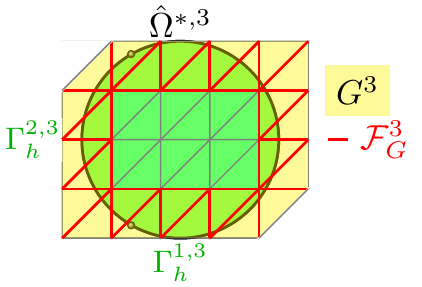}
\caption{Schematic representation of ghost penalty faces.}
\label{fig: ghost face notation}
\end{figure}

We integrate the equations of linear elasticity over each subdomain $\Omega^i_h$ (see Figure~\ref{subfig: geometryapprox}) in our weak finite element formulation. In the case of small overlaps of the physical domain $\Omega^i_h$ with an element $K$, i.e. $K_i = \Omega^i_h \cap K$ small, the resulting finite element system matrix can become ill-conditioned. To circumvent this problem, we will solve the problem of linear elasticity on an extension of the physical domain (the so-called fictitious domains) instead of only on the domain $\Omega^i_h$. More precisely, these fictitious domains are local background meshes for each of the subdomains $i \in \mathcal{I}_\Domain$ defined as the set of elements in the background mesh $\bmesh$ that are fully or partially covered by the physical domain $\Omega^i_h$, \textit{i.e.}
\begin{equation}
\DomainHat^{*,i} := \{ K \in \bmesh: K \cap \Domain^i_h \neq \emptyset\}.
\end{equation}
Figure~\ref{fig: fictitious domain meshes} illustrates the three fictitious domains for the three phase particulate composite example mentioned in the previous section. \\
To make this extension from the physical domain $\Omega^i_h$ onto the fictitious domain $\DomainHat^{*,i}$ without loosing weak consistency, we will employ a regularisation technique on the set of elements intersected by interface $\Interface^i_h$. We denote this set of elements by
\begin{equation}
G^{i} := \{ K \in \bmesh: K \cap \Interface^i_h \neq \emptyset\} 
\end{equation}
and define a set of element faces associated with $G^{i}$ 
\begin{equation}
\mathcal{F}^i_G := \{ F \in \mathcal{F}(K,K'):  K \in  G^{i} \text{ or } K' \in G^{i}\},
\label{equ: ghost penalty facets}
\end{equation}
which we call \textit{ghost penalty faces}. Here,  $\mathcal{F}(K,K') = K \cap K'$ denotes the face between  element $K \in \bmesh$ and its neighbouring element $K' \in \bmesh$. The ghost penalty faces, $\mathcal{F}^i_G$, and the set of intersected elements, $G^{i}$, are illustrated in Figure~\ref{fig: ghost face notation}.

Now, we define our finite element space over the fictitious domains as: Let $\mathcal{U}_h$ denote the vector valued space of continuous piecewise linear polynomials defined on the background mesh $\bmesh$. Then, for each $i \in I_{\Domain}$, we define the space 
\begin{equation}
\mathcal{U}^i_h =  \mathcal{U}_h |_{\DomainHat^{*,i}} \, ,
\end{equation} 
\textit{i.e.} $\mathcal{U}^i_h$ is defined as the restriction of $\mathcal{U}_h$ to the fictitious domain mesh $\DomainHat^{*,i} \subset \bmesh$. 

Recall that two different fictitious domains overlap in cells that are intersected by interfaces. This overlapping, subdomain-based definition of approximation spaces, which is consistent with the classical CutFEM paradigm \cite{hansbohansbo2004,burmanclaus2015}, naturally leads to a displacement field that may be discontinuous across interfaces.

\paragraph{Piecewise linear spaces for the LaTIn interface quantities.}

We will seek continuous piecewise linear approximations of the interface displacements $\InterfaceDisp^{i,j}_h$ and interface forces $\InterfaceForce^{i,j}_h$, for all linearly interpolated subdomain interfaces $\Interface_h^{i,j}$, with $(i,j) \in  \mathcal{I}_\Interface$. Similarly to the bulk equations, we will extend the interface quantities from the interface $\Interface_h^{i,j}$ onto the entire band of elements which are intersected by the interface $\Interface_h^{i,j}$ defined as 
\begin{equation}
G^{i,j} := \{ K \in \bmesh: K \cap \Interface^{i,j}_h \neq \emptyset\} 
\end{equation}
 to prevent ill-conditioning. These extensions will be denoted by $\InterfaceDisp^{*,i,j}_h$ and  $\InterfaceForce^{*,i,j}_h$, respectively. In particular, these extensions fulfil 
\begin{equation}
\InterfaceDisp^{*, i,j}_h |_{\Interface^{i,j}_h} = \InterfaceDisp^{i,j}_h |_{\Interface^{i,j}_h}, \,\,\, \InterfaceForce^{*, i,j}_h |_{\Interface^{i,j}_h} = \InterfaceForce^{i,j}_h |_{\Interface^{i,j}_h}.
\end{equation} 
Algebraic operations on interface quantities $\InterfaceDisp^{i,j}_h$ and $\InterfaceForce^{i,j}_h$ are then performed implicitly through manipulating the nodal values of their volume extensions. 
Now let us formally define the continuous piecewise linear finite element space in which we will seek the interface quantity extensions $\InterfaceDisp^{*, i,j}_h$ and $\InterfaceForce^{*, i,j}_h$ as 
\begin{equation}
\mathcal{V}_h^{i,j} = \mathcal{U}_h|_{G^{i,j}},
\label{equ: FEM space interface}
\end{equation}
\textit{i.e.} the vector valued space of continuous piecewise linear polynomials on $G^{i,j}$. Let $N_{d,G^{i,j}}$ denote the number of degrees of freedom in $\mathcal{V}_h^{i,j}$.







\subsubsection{Regularised formulation of the bulk equations}

\begin{figure}
\centering
\includegraphics[width=.5\textwidth]{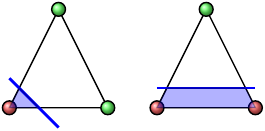}
\caption{Cases of interface-element intersections that lead to ill conditioning of the system matrices in the abscence of appropriate regularisation.}
\label{fig: bad cut cases}
\end{figure}
 
As briefly mentioned above, integration over element parts $K_i=K \cap \Omega_h^i$ in the interface region can lead to conditioning issues. For instance, cutting an element in an arbitrary location may result in elements that have an almost-zero intersection with the physical domain (see Figure~\ref{fig: bad cut cases}). As a consequence, degrees of freedom of such elements can move in an uncontrolled manner, which results in severe ill-conditioning and convergence issues for linear and nonlinear solution algorithms.


In the context of our LaTIn-CutFEM method, we will address this issue through a ghost penalty extension of the coercivity of the linear elasticity problems to the entire fictitious domains.  This is done by penalising the jump of fictitious tractions (\textit{i.e.} defined over the entire fictitious domain) across elements that are cut by an unfitted LaTIn interface. When the penalty term is scaled correctly, this consistent approach ensures  that the convergence of the condition number with mesh refinement is that of a conforming FEM. The interested reader is referred to  \cite{Burman_2010_a,burmanclaus2015} for extensive derivations and applications of the ghost penalty approach in the context of Nitsche's method.

The ghost-penalty regularised variational statements associated with the static equilibrium of the subdomains read as
\begin{center} 
\textit{For all $i  \in \mathcal{I}_\Omega$,  field  $\DispHatTrial \in \mathcal{U}^i_\text{h}$ defined over the fictitious subdomain $\hat{\Domain}^{*,i}$ must satisfy, for any variation $\DispTest \in \mathcal{U}^i_\text{h}$ ,}
\end{center}
\begin{equation}
\begin{array}{rcl}
a_\text{D}^{i}(\DispHatTrial,\DispTest) & + & a_\text{t}^{i}(\DispHatTrial,\DispTest) + \displaystyle a_\text{N}^{i}(\DispHatTrial,\DispTest) + j_\text{u}^{i}(\DispHatTrial,\DispTest)   \\
& = & \displaystyle
l_\text{f}^{i}(\DispTest) +  l_\text{N}^{i}(\DispTest)   \, .
\end{array}
\label{equ: bulk equations discrete}
\end{equation}
The terms introduced in the previous formulation are listed and commented below. 
\begin{itemize}
\item \textbf{``Mechanical" virtual work.} The virtual work of the internal forces in domain $\Domain^{i}_h$, $i \in \mathcal{I}_\Omega$, reads as
\begin{equation}
a_\text{D}^{i}(\DispHatTrial,\DispTest) =  \int_{\Domain^{i}_h } \StrainBasic(\DispHatTrial) : \Hooke : \StrainBasic(\DispTest) \, d\Domain \, ,
\end{equation}
the work of the tractions from neighbouring subdomains reads as 
\begin{equation}
a_\text{t}^{i}(\DispHatTrial,\DispTest) = - \int_{\Gamma^{i}_h} \left( \left( \Hooke : \StrainBasic(\DispHatTrial) \right) \cdot \V{n}^i \right)   \cdot \DispTest \, d \Gamma \, ,
\label{equ: surface traction stress}
\end{equation}
and the work of the external forces, excluding that of the interface, reads as
\begin{equation}
l_\text{f}^{i}(\DispTest) =  
 \int_{\Domain^{i}_h} \BodyForce \cdot \DispTest \, d\Domain  + 
\int_{ \DDomainNeumannFEM \cap \partial \Domain^{i}_h } \NeumannBC \cdot  \DispTest  \, d\Interface
  \, .
\end{equation}
\item \textbf{Nitsche terms for non-conforming Dirichlet boundaries.} On approximate boundary $\partial \Omega^{\text{u}}_h$, Nitsche's method allows us to apply the Dirichlet conditions weakly. The corresponding terms read as
\begin{equation}
\begin{aligned}
  \displaystyle a_\text{N}^{i}(\DispHatTrial,\DispTest) = & \displaystyle
- \int_{ \partial \Omega^{\text{u}}_h }  
 \DispTest  \cdot 
\left( \Hooke : \StrainBasic(\DispHatTrial) \right) \cdot \Normal 
  \, d\Interface -  \displaystyle
\int_{ \partial \Omega^{\text{u}}_h }  
 \DispHatTrial \cdot 
\left( \Hooke : \StrainBasic(  \DispTest ) \right) \cdot \Normal
  \, d\Interface 
\\
&+  \displaystyle
\int_{ \partial \Omega^{\text{u}}_h}  
\frac{\alpha E^i}{h} \,    \DispHatTrial   \cdot \DispTest 
  \, d\Interface
  \, ,
\end{aligned}
\end{equation}
and
\begin{equation}
l_\text{N}^{i}(\DispTest)
=
\int_{ \partial \Omega^{\text{u},\text{h}} }  
\frac{\alpha E^i}{h}  \,  \DirichletBC^\text{h} \cdot \DispTest 
  \, d\Interface
  \, .
\end{equation}
The $1/h$ scaling ensures stability and optimal convergence of the symmetric Nitsche's method (see \textit{e.g.}  \cite{burmanclaus2015}). Note that, when using penalty-type formulations to enforce Dirichlet conditions, the test field $\DispTest$ is not required to vanish over Dirichlet boundaries. 

\item \textbf{The ghost penalty term} reads as

\begin{equation}
j_\text{u}^{i}(\DispHatTrial,\DispTest)
= \sum_{F \in \mathcal{F}^i_G }
  \int_{F}
 \frac{\gamma_g \, h}{E^i} \, \jump{ \left(\Hooke : \StrainBasic(\DispHatTrial) \right)  \cdot \Normal_F} \cdot  \jump{ \left( \Hooke : \StrainBasic(\DispTest) \right)  \cdot \Normal_F  }
\, ds
 \, .
\label{equ: ghost penalty terms}
\end{equation}
Here, $\jump{x \cdot \Normal_F }$ denotes the normal jump of the quantity $x$
over the face, $F$, defined as $\jump{x \cdot \Normal_F } = \left. x
\right|_{K}{n}_F  - \left. x
\right|_{K'}{n}_F$, where ${n}_F$ denotes a unit
normal to the facet $F$ with fixed but arbitrary orientation. The ghost penalty parameter, $\gamma_g >0$, is chosen to be sufficiently large to regularise the solution in the interface regions. In the previous definition, we have penalised the jump of stress across all ghost penalty faces as defined in \eqref{equ: ghost penalty facets}. These include intersected faces and faces of intersected elements that couple the element to an element in the interior of $\Omega^i_h$. This way, we extend the traction vector into the fictitious domain to prevent ill-conditioning in a consistent way (i.e. for smooth tractions the jump terms vanish).
 
\end{itemize}

\subsubsection{Multiscale split of interface quantities and weak compatibility conditions}
\label{subsec: multiscale heart local}

\begin{figure}[htb]
\centering
\includegraphics[width=0.8 \textwidth]{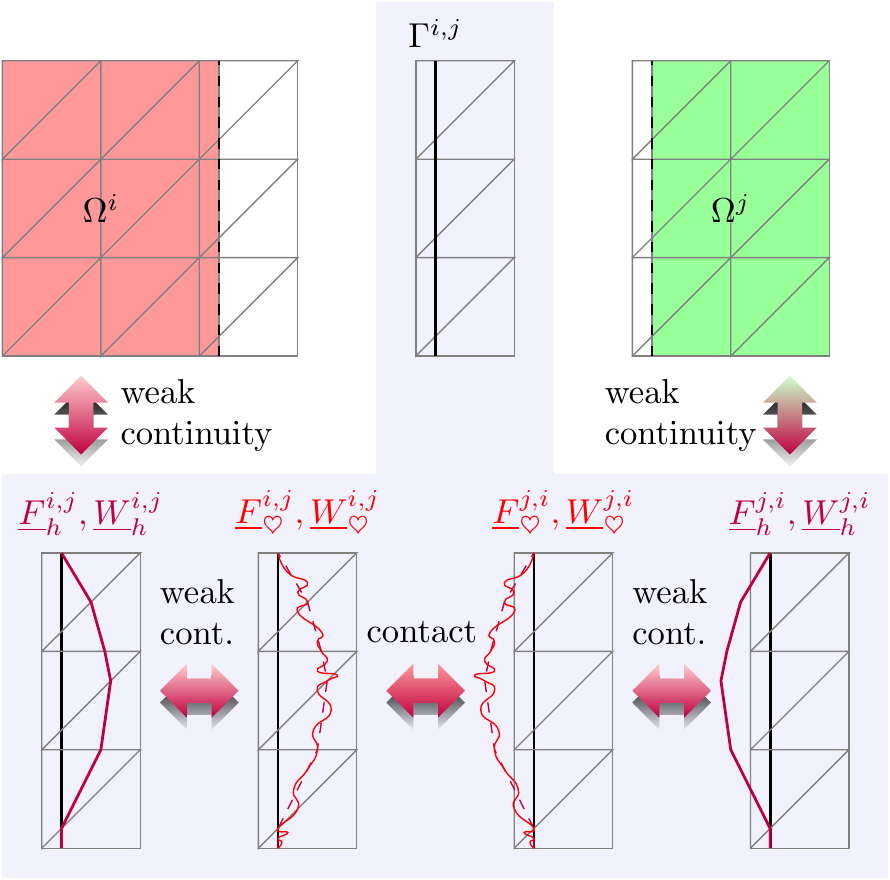}
\caption{Schematic of LaTIn-CutFEM interface/bulk and interface/interface coupling.}
\label{fig:LaTInCutFEM}
\end{figure}

\paragraph{Heart quantities}
Forcing piecewise linear quantities $\{ (\InterfaceForce^{i,j}_h,\InterfaceDisp^{i,j}_h) \}$ to satisfy the nonlinear equations of unilateral contact is doomed to failure due to the intrinsically nonlinear nature of unilateral contact. Recall that the mainstream approach is to choose these quantities as piecewise constant, which is compatible with complex interface laws. However, such approaches lead to numerical instabilities. 

We propose to introduce a new pair of dual interface fields $\{  \InterfaceDisp^{i,j}_{\heartsuit} \}$, $\{ \InterfaceForce^{i,j}_{\heartsuit} \}$, chosen in $L^2(\Gamma)$. These quantities are the ones that are required to satisfy the contact conditions, as represented in Figure \ref{fig:LaTInCutFEM}. They are then weakly connected to  $\{ (\InterfaceForce^{i,j}_h,\InterfaceDisp^{i,j}_h) \}$. More precisely, a regularised $L^2$ projection of $\{ \InterfaceDisp^{i,j}_{\heartsuit} \}$, $\{ \InterfaceForce^{i,j}_{\heartsuit} \}$ into the finite element piecewise linear interface space is required to equate $\{ (\InterfaceForce^{i,j}_h,\InterfaceDisp^{i,j}_h) \}$ , whilst the remainder, which we term ``fluctuation", will satisfy a Robin condition.

In practice, the contact conditions are satisfied by the heart quantities at a set of quadrature points, more precisely the quadrature points of the subtriangulation of the cut interface (3 quadrature points per element in the sub-triangulation of interface $\{ \Gamma^{i,j}_h \}$).


\paragraph{Compatibility between interface and bulk}

For all interfaces $\Gamma^{i,j}_h$, $(i,j) \in \mathcal{I}_{\Interface}$, we couple the interface and bulk quantities through  approximate compatibility conditions 
\begin{equation}
\begin{aligned}
\InterfaceDisp^{i,j}_h  &= \DispHatTrial,  \quad \InterfaceDisp^{j,i}_h  = \Disp_h^{*,j}  \quad &\mbox{ on } \Gamma^{i,j}_h, \\ 
\int_{\Gamma^{i,j}_h} \InterfaceForce^{i,j}_h  \cdot  \delta \Disp &=  \int_{\Gamma^{i,j}_h} \left( \left( \Hooke : \StrainBasic(\DispHatTrial) \right) \cdot \V{n}^i \right)   \cdot \DispTest \, d \Gamma \quad &\forall \DispTest \in  \mathcal{U}^i_\text{h}, \\ 
\int_{\Gamma^{i,j}_h} \InterfaceForce^{j,i}_h  \cdot  \delta \Disp &=  \int_{\Gamma^{i,j}_h} \left( \left( \Hooke : \StrainBasic(\Disp^{*,j}_h) \right) \cdot \V{n}^i \right)   \cdot \DispTest \, d \Gamma \quad &\forall \DispTest \in  \mathcal{U}^j_\text{h}.
\end{aligned}
\end{equation}
These are standard conditions whereby the kinematic continuity is enforced exactly, whilst the continuity of dual quantities is only enforced on average, \textit{i.e.} in a finite element sense. When solving our problem using the LaTIn algorithm, these two conditions will be introduced in \eqref{equ: bulk equations discrete} to yield our working expression of the linear stage.


\paragraph{Compatibility between interface and heart}

We project the heart quantities onto the continuous piecewise linear approximation space $\mathcal{V}^{i,j}_h$ \eqref{equ: FEM space interface} for each $\Gamma^{i,j}_h$ using a stabilised $L^2$-projection as follows. Find the extensions $\InterfaceForce^{*,i,j}_h \in \mathcal{V}^{i,j}_h$ and $\InterfaceDisp^{*,i,j}_h \in \mathcal{V}^{i,j}_h$ such that
\begin{equation}
\begin{aligned}
\int_{\Interface^{i,j}_h} (\InterfaceForce^{*, i,j}_h -\InterfaceForce^{i,j}_{\heartsuit} ) \,\, \delta \InterfaceForce^{i,j}_h \,d\Gamma+  j_F(\InterfaceForce^{*,i,j}_h, \delta \InterfaceForce^{i,j}_h)  = \V{0}\, \quad \forall  \delta \InterfaceForce^{i,j}_h \in \mathcal{V}^{i,j}_h,\\
\int_{\Interface^{i,j}_h} (\InterfaceDisp^{\star,i,j}_h -\InterfaceDisp^{i,j}_{\heartsuit} ) \,\, \delta \InterfaceDisp^{i,j}_h \,d\Gamma +  j_F(\InterfaceDisp^{*,i,j}_h, \delta \InterfaceDisp^{i,j}_h)  = \V{0}\, \quad \forall  \delta \InterfaceDisp^{i,j}_h \in \mathcal{V}^{i,j}_h,
\end{aligned}
\label{eq:ContinuityLocal}
\end{equation}
where
\begin{equation}
 j_F( \InterfaceForce^{*,i,j}_h, \delta  \InterfaceForce^{i,j}_h ) =  \sum_{F \in \mathcal{F}_I^{i,j}} \gamma_{\Pi} h^2 \int_{F} \jump{\nabla   \InterfaceForce^{*,i,j}_h \cdot \Normal_F} \cdot  \jump{\nabla \delta \InterfaceForce^{i,j}_h, \cdot \Normal_F} ds
\end{equation}
regularises the interface fields with a penalty parameter $\gamma_{\Pi}>0$ in the band of elements intersected by $\Gamma^{i,j}_h$. Here,  
\begin{equation}
\mathcal{F}^{i,j}_I :=\{ F \in \mathcal{F}(K,K'):  K \in  G^{i,j} \text{ and } K' \in  G^{i,j}\}.
\end{equation} 
denotes the set of faces intersected by $\Gamma^{i,j}_h$.  

\paragraph{Heart closure}

We close our hybrid-mixed formulation of contact by requiring that the fluctuation of the heart quantities around $\InterfaceForce^{i,j}_h$ and $\InterfaceDisp^{i,j}_h$ satisfy a Robin condition, which reads as 
\begin{equation}
\left( \InterfaceForce^{i,j}_h -\InterfaceForce^{i,j}_{\heartsuit} \right) - \beta \left( \InterfaceDisp^{i,j}_h -\InterfaceDisp^{i,j}_{\heartsuit} \right) = 0 \, ,
\end{equation}
where $\beta$ is an algorithmic parameter that is homogeneous to a stiffness. If we set it equal to the LaTIn search direction parameter, the resulting local stage of the LaTIn solver is a classical local stage followed by a ``coarse-scale" filtering step, as described in the next section. 

This closure introduces a lack of consistency in the formulation. However, we expect the associated consistency error to vanish with mesh refinement with optimal order, as it directly affects the remainder of the heart fields after projection in the finite element space. This will be shown numerically through examining the convergence properties of the method. \\

It is noticeable that the non-local continuity equations \eqref{eq:ContinuityLocal} make the local stage non-local. In our opinion, this is not a significant limitation of the method. This is because the non-linearities can still be treated locally and semi-explicitly as per usual, in a first sub-stage, whilst the second sub-stage requires solving a unique set of independent (and small) linear problems, each of these problems corresponding to an interface between two adjacent subdomains.



\subsection{Iterative algorithm for the fully discrete mixed cut finite element problem}

We are now in the position to use the regularisation and modifications introduced previously to apply the LaTIn algorithm to the mixed interface quantities $\{ (\InterfaceForce^{i,j}_h,\InterfaceDisp^{i,j}_h) \}$.

\subsubsection{Regularised linear stage over the fictitious domains}

\paragraph{Linear system of equations}

\begin{center} 
\textit{For all $i  \in \mathcal{I}_\Omega$, find the displacement field  $\DispHatTrial \in \mathcal{U}^i_\text{h}$ defined over the fictitious subdomain $\hat{\Domain}^{*,i}$ such that, for any variation $\DispTest \in \mathcal{U}^i_\text{h}$ ,}
\end{center}
\begin{equation}
\begin{array}{rcl}
a_\text{D}^{i}(\DispHatTrial,\DispTest) & + & \displaystyle a_\text{k}^{i}(\DispHatTrial,\DispTest) +  a_\text{N}^{i}(\DispHatTrial,\DispTest) + j_\text{u}^{i}(\DispHatTrial,\DispTest) \\
& = & \displaystyle
l_\text{f}^{i}(\DispTest) + \hat{l}_\text{k}^{i}(\DispTest) +  l_\text{N}^{i}(\DispTest)  \, .
\end{array}
\label{equ: bulk equation linear stage}
\end{equation}
Two terms corresponding to the LaTIn approximation of the virtual work of the forces applied by neighbouring subdomains have been introduced. These forces are found iteratively as described in Section \ref{sec:LaTin}. At a given iteration of the LaTIn solver, the approximate equilibrium of the subdomain is described via an augmentation term (\textit{i.e.} ``penalty" or ``regularisation" boundary term induced by the descent search direction), which reads
\begin{equation}
a_\text{k}^{i}(\DispHatTrial,\DispTest)  = \int_{\Gamma^{i}_h }  k^- \,  \DispHatTrial \cdot \DispTest   \, d\Interface
 \, .
\end{equation}
Mixed quantities from the previous half-iterations appear on the right-hand side as follows
\begin{equation}
l_\text{k}^{i}(\DispTest)  = \int_{\Gamma^{i}_h }    
\left( \InterfaceForceHat^{i} + k^{-} \, \InterfaceDispHat^{i} \right) \cdot  \DispTest  
 \, d\Interface
 \, .
\end{equation}
Note that the interface quantities $\InterfaceForceHat^{i}$, $\InterfaceDispHat^{i}$ over the interface $\Interface^i_h$ are given by the interface quantities  $\InterfaceForceHat^{i,j}$, $\InterfaceDispHat^{i,j}$ by the following sum 
\begin{equation}
\begin{aligned}
\int_{\Gamma^{i}_h }    
\left( \InterfaceForceHat^{i}_h + k^{-} \, \InterfaceDispHat^{i}_h \right) \cdot  \DispTest_h  
 \, d\Interface  
 = \sum_{\substack{j \in \mathcal{I}_{\Interface^i}}} \left( \int_{\substack{\Interface^{k,l}_h}} \InterfaceForceHat^{i,j}_h + k^{-} \, \InterfaceDispHat^{i,j}_h \right) \cdot  \DispTest_h  
 \, d\Interface, 
 \end{aligned}
\end{equation}
where  $(k,l)=(i,j)$  if  $i<j$ and  $(k,l)=(j,i)$ if  $i>j$.
Here, 
\begin{equation}
\mathcal{I}_{\Gamma^i} = \{ j \in \mathcal{I}_{\Domain}: \Gamma^{i,j}_h \neq \emptyset \mbox{ if } i<j \mbox{; or } \Gamma^{j,i}_h \neq \emptyset \mbox{ if } i>j \}.
\end{equation}
For example, in Figure~\ref{fig: schematic notation}, $\Gamma^3 = \Gamma^{1,3} \cup \Gamma^{2,3}$, $\mathcal{I}_{\Gamma^3}=\{1,2\}$. 

\paragraph{Post-processing of interface quantities}

\begin{figure}
\centering 
\includegraphics[width=.5\textwidth]{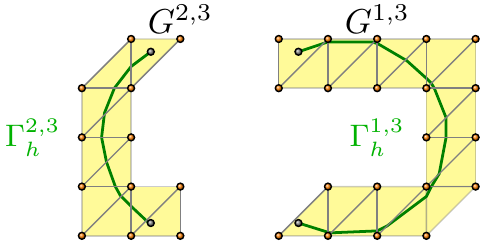}
\caption{Decomposition of intersected cells for extended interface quantities into meshes $G^{i,j}$.}
\label{fig: meshes local stage}
\end{figure}

For all $(i,j)\in\mathcal{I}_{\Interface}$, we seek the extended interface displacement fields, $\InterfaceDisp^{*,i,j}_h \in \mathcal{V}^{i,j}_h$,  $\InterfaceDisp^{*,j,i}_h \in \mathcal{V}^{j,i}_h$, which fulfil
\begin{equation}
\begin{aligned}
\InterfaceDisp^{*,i,j}_h|_{\Interface^{i,j}_h}  &= \InterfaceDisp^{i,j}_h|_{\Interface^{i,j}_h}  = \DispHatTrial |_{\Interface^{i,j}_h}, \\
 \InterfaceDisp^{*,j,i}_h|_{\Interface^{i,j}_h} & =  \InterfaceDisp^{j,i}_h|_{\Interface^{i,j}_h} = \Disp_h^{*,j}|_{\Interface^{i,j}_h}.   \\ 
\end{aligned}
\label{equ: strong interface displacement}
\end{equation}
We obtain these extended interface displacement fields by setting $\InterfaceDisp^{*,i,j}_h$ equal to $ \DispHatTrial$ and  $\InterfaceDisp^{*,j,i}_h$ equal to $ \Disp^{*,j}_h$ in each degree of freedom in the band of intersected elements $G^{i,j}$, \textit{i.e.}
\begin{equation}
\begin{aligned}
\InterfaceDisp^{*,i,j}_h(\pmb{x}_i)  = \DispHatTrial(\pmb{x}_i) \quad \forall  i=1, \dots, N_{d,G^{i,j}}, \\
\InterfaceDisp^{*,j,i}_h(\pmb{x}_i)  = \Disp^{*,j}_h(\pmb{x}_i) \quad \forall  i=1, \dots, N_{d,G^{i,j}}.
\label{equ: extended interface disp local}
\end{aligned}
\end{equation}
This means, as $\InterfaceDisp^{*,i,j}_h$ and $\InterfaceDisp^{*,j,i}_h$ is equal to $\DispHatTrial $ and $\Disp^{*,j}_h$  in all degrees of freedom around the interface, it is also equal to the $\DispHatTrial$ and $\Disp^{*,j}_h$ on the restriction to the interface $\Gamma^{i,j}_h$, hence we satisfy \eqref{equ: strong interface displacement}. The location of the degrees of freedom, $\pmb{x}_i$,  in the surface band mesh $G^{i,j}$ are illustrated in Figure~\ref{fig: meshes local stage}. 

%
Next, we can obtain the extended interface forces $\InterfaceForce^{*,i,j}_h \in \mathcal{V}^{i,j}_h$ and $\InterfaceForce^{*,j,i}_h \in \mathcal{V}^{i,j}_h$ for each $(i,j) \in I_{\Interface}$ through the search direction 
\begin{equation}
\begin{aligned}
\InterfaceForce^{*,i,j}_h(\pmb{x}_i) = \InterfaceForceHat^{*,i,j}_h(\pmb{x}_i) + k^- (\InterfaceDispHat^{*,i,j}_h(\pmb{x}_i) - \InterfaceDisp^{*,i,j}_h(\pmb{x}_i)) \quad \forall i=1, \dots, N_{d,G^{i,j}}, \\
\InterfaceForce^{*,j,i}_h(\pmb{x}_i) = \InterfaceForceHat^{*,j,i}_h(\pmb{x}_i) + k^- (\InterfaceDispHat^{*,j,i}_h(\pmb{x}_i) - \InterfaceDisp^{*,j,i}_h(\pmb{x}_i)) \quad \forall i=1, \dots, N_{d,G^{i,j}}.
\end{aligned}
\label{equ: interface force linear stage}
\end{equation}

Note that this gives us a piecewise linear continuous representation of the interface quantities. 

\paragraph{Multiscale local stage}

Having obtained the approximation for the interface force, $\InterfaceForce^{*,i,j}_h$, and the interface displacement, $\InterfaceDisp^{*,i,j}_h$, we proceed to calculating local quantities following the multi-scale strategy introduced in section~\ref{subsec: multiscale heart local}. 

Firstly, to enforce frictionless unilateral contact,  we determine the heart quantities for each interface $\Gamma^{i,j}_h$ in each quadrature point $\pmb{x}_q$ along the interface $\Gamma^{i,j}_h$ by


\begin{equation}
\InterfaceForce_{\heartsuit}^{i,j}(\pmb{x}_q) =  \frac{1}{2}(\InterfaceForce^{*,i,j}_h(\pmb{x}_q)-\InterfaceForce^{*,j,i}_h(\pmb{x}_q) + k^+ (\InterfaceDisp^{*,j,i}_h(\pmb{x}_q) -  \InterfaceDisp^{*,i,j}_h(\pmb{x}_q))) \cdot \Normal_K^{i,j}(\pmb{x}_q).
\end{equation}
Here, $\Normal_K^{i,j}$ denotes a piecewise constant normal obtained from the level-set function $\phi^j_h$
\begin{equation}
\Normal_K^{i,j}(\pmb{x}_q) = - \left. \frac{\nabla \phi^{j}_h(\pmb{x})}{|\nabla \phi^{j}_h(\pmb{x})|}\right|_{\pmb{x}=\pmb{x_q}}.
\label{equ: compute normals} 
\end{equation}
To fulfil frictionless contact, we test in each quadrature point, if  $\InterfaceForce_{\heartsuit}^{i,j}(\pmb{x}_q) >0$ (i.e. no contact), in which case we set $\InterfaceForce_{\heartsuit}^{i,j}(\pmb{x}_q) =0$. 
The heart quantity now satisfies frictionless contact in each point and we can determine the local interface force extension $\InterfaceForceHat^{*,i,j}_h \in \mathcal{V}^{i,j}_h$ from the stabilised $L^2$-projection \eqref{eq:ContinuityLocal}.
Next, we use $\InterfaceForceHat^{*,i,j}_h$ to obtain
\begin{equation}
\InterfaceForceHat^{*,j,i}_h(\pmb{x}_i) = - \InterfaceForceHat^{*,i,j}_h(\pmb{x}_i)
\label{local stage force ji}
\end{equation}
\begin{equation}
\begin{aligned}
\InterfaceDispHat^{*,i,j}_h(\pmb{x}_i) = \InterfaceDisp^{*,i,j}_h(\pmb{x}_i) + \frac{1}{k^+} ( \InterfaceForceHat^{*,i,j}_h(\pmb{x}_i) -  \InterfaceForce^{*,i,j}_h(\pmb{x}_i)) \\
\InterfaceDispHat^{*,j,i}_h(\pmb{x}_i) = \InterfaceDisp^{*,j,i}_h(\pmb{x}_i) + \frac{1}{k^+} ( \InterfaceForceHat^{*,j,i}_h(\pmb{x}_i) - \InterfaceForce^{*,j,i}_h(\pmb{x}_i)). \\
\end{aligned}
\label{equ: local stage contact displacements}
\end{equation}
for all $i=1, \dots, N_{d,G^{i,j}}$. Note that these algebraic relations together with the stabilised $L^2$-projection \eqref{eq:ContinuityLocal} gives us the continuous and piecewise linear extended interface quantities that we are seeking, \textit{i.e.} $\InterfaceForceHat^{*,j,i}(\pmb{x}_i) \in \mathcal{V}^{j,i}_h$,  $\InterfaceDispHat^{*,i,j}_h(\pmb{x}_i) \in \mathcal{V}^{i,j}_h$,  $\InterfaceDispHat^{*,j,i}_h(\pmb{x}_i) \in \mathcal{V}^{j,i}_h$. 

\subsubsection{Pseudo-code of the LaTIn-CutFEM method} 

Algorithm~\ref{alg: cutfem latin algorithm} summaries the LaTin-CutFEM algorithm presented above. 

\begin{algorithm}
    \caption{LaTin-CutFEM algorithm for unilateral contact}
    \label{alg: cutfem latin algorithm}
    \begin{algorithmic}[1] 
    \State Approximate level set functions through linear interpolation $\to \phi^i_h$
    \State Compute intersection of level set functions with background mesh and obtain linear approximations of $\Interface^{i,j}_h$ and $\Domain^i_h$
    \State Set $\InterfaceDispHat^{*,i,j}_h=\mathbf{0}$, $\InterfaceDispHat^{*,j,i}_h=\mathbf{0}$ and $\InterfaceForceHat^{*,i,j}_h = \mathbf{0}$, $\InterfaceForceHat^{*,j,i}_h = \mathbf{0}$, $\forall (i,j) \in \mathcal{I}_{\Interface}$.
    \While{$it<it_{max}$} \Comment{$it_{max}$ is maximum iteration count}
    \Procedure{Linear Stage}{$\InterfaceDispHat^{*,i,j}_h,\InterfaceForceHat^{*,i,j}_h$} 
    \State Find $\DispHatTrial$ by solving bulk equation \eqref{equ: bulk equation linear stage} $\forall i \in \mathcal{I}_{\Domain}$
    \State Determine $\InterfaceDisp^{*,i,j}_h$ from \eqref{equ: extended interface disp local} and $\InterfaceForce^{*,i,j}_h$ from \eqref{equ: interface force linear stage} using $\DispHatTrial$.
    \State Perform relaxation \eqref{equ: relaxation linear stage} for $\InterfaceDisp^{*,i,j}_h$ and $\InterfaceForce^{*,i,j}_h$.  
    \State \textbf{return} $\DispHatTrial$, $\InterfaceDisp^{*,i,j}_h$,$\InterfaceForce^{*,i,j}_h$ 
         \EndProcedure
       \Procedure{Local Stage}{$\DispHatTrial$, $\InterfaceDisp^{*i,j}_h$,$\InterfaceForce^{*,i,j}_h$} 
       \State for \textbf{contact interfaces} $\Interface^{i,j}_h$ :
       \State \quad Compute normals $\Normal_K^{i,j}$ \eqref{equ: compute normals} 		  
       \State \quad Obtain $\InterfaceForceHat^{*,i,j}_h$ from \eqref{eq:ContinuityLocal}, \eqref{local stage force ji} and $\InterfaceDispHat^{*,i,j}_h$ from \eqref{equ: local stage contact displacements}
        \State \textbf{return} $\InterfaceForceHat^{*,i,j}_h$, $\InterfaceDispHat^{*,i,j}_h$ 
       \EndProcedure
       \State it +=1

    \EndWhile\label{while}
    \end{algorithmic}
\end{algorithm}


\section{Numerical investigations}


\subsection{Controlled Condition Number: Square Domain with a Crack} 

\begin{figure}
[htbp]
\centering
\subfloat[$\epsilon_x=0.5$, $\epsilon_y=0.25$.]{
\begin{tikzpicture}
\node[inner sep=0pt] (cont) at (0,0)
    {\includegraphics[width=.4\textwidth]{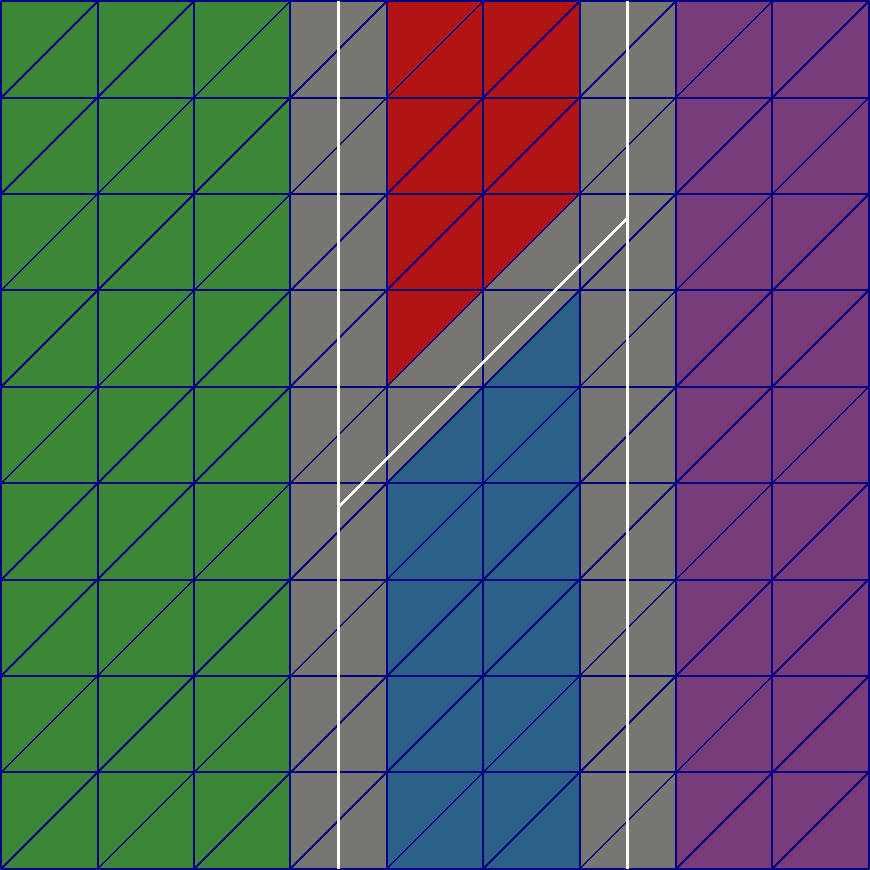}};
\node[white] at (0.4,2) {{\Large $\Omega_1$}};   
\node[white] at (0.4,-1.5) {{\Large $\Omega_2$}};  
\node[white] at (-2,0) {{\Large $\Omega_3$}}; 
\node[white] at (2.4,0) {{\Large $\Omega_4$}};  
\draw[arrows=-latex,thick,yellow] (0.35,0.5) -- (0.35,-0.2);
\end{tikzpicture}
\label{subfig: simply bad cut case 1}}\hspace{.2cm}
\subfloat[$\epsilon_x=0.5$, $\epsilon_y=10^{-10}$.]{\includegraphics[width=.4\textwidth]{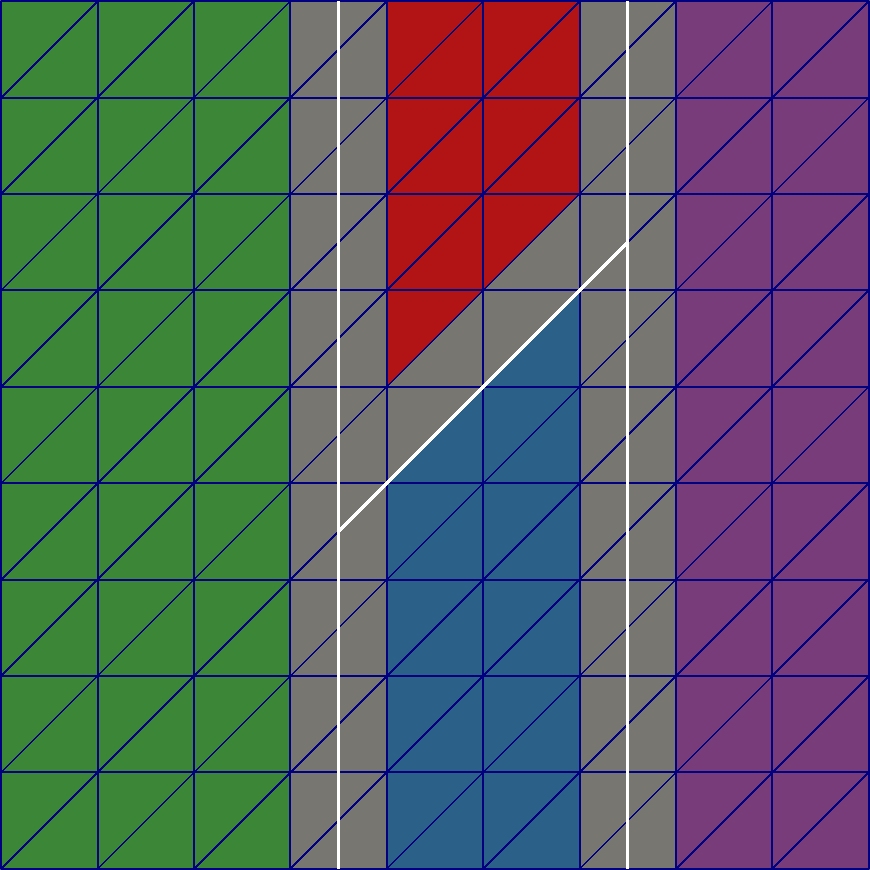}
\label{subfig: simply bad cut case 2}}\\

\subfloat[$\epsilon_x=\epsilon_y$, $\epsilon_y=0.25$.]{
\begin{tikzpicture}
\node[inner sep=0pt] (cont) at (0,0)
    {\includegraphics[width=.4\textwidth]{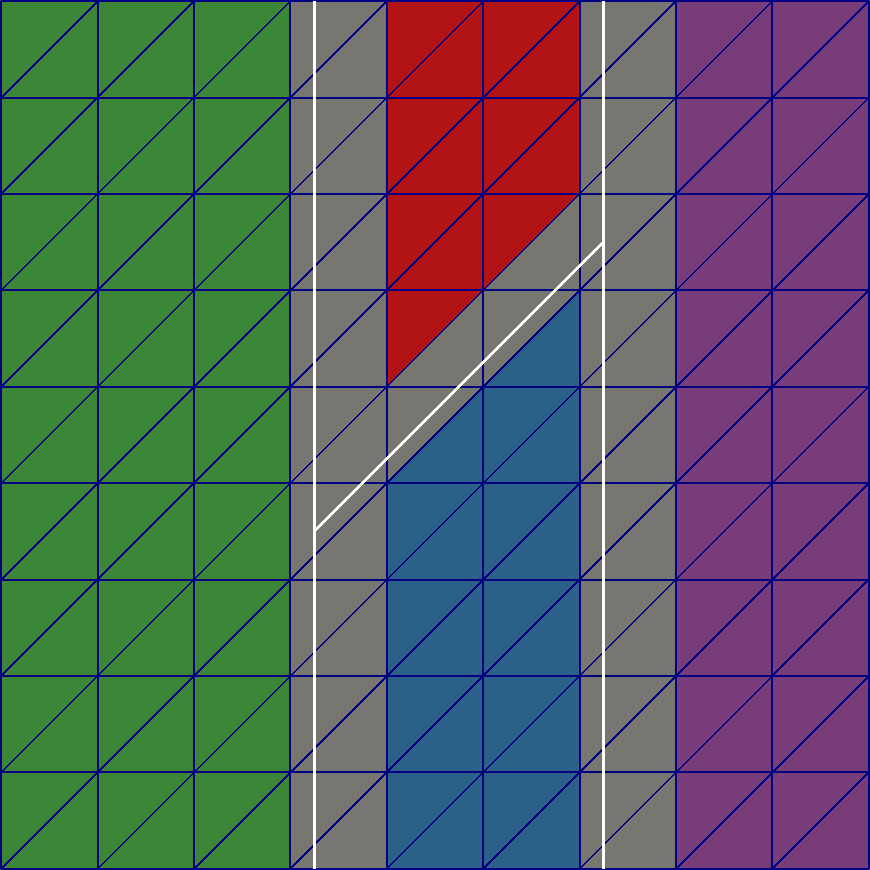}};
\node[white] at (0.4,2) {{\Large $\Omega_1$}};   
\node[white] at (0.4,-1.5) {{\Large $\Omega_2$}};  
\node[white] at (-2,0) {{\Large $\Omega_3$}}; 
\node[white] at (2.4,0) {{\Large $\Omega_4$}};  
\draw[arrows=-latex,thick,yellow] (1.2,-0.3) -- (0.5,-0.3);
\draw[arrows=latex-,thick,yellow] (-1.55,-0.3) -- (-0.85,-0.3);
\draw[arrows=-latex,thick,yellow] (0.35,0.5) -- (0.35,-0.2);
\end{tikzpicture}
\label{subfig: double bad cut case 1}}\hspace{.2cm}
\subfloat[$\epsilon_x=\epsilon_y$, $\epsilon_y=10^{-10}$.]{\includegraphics[width=.4\textwidth]{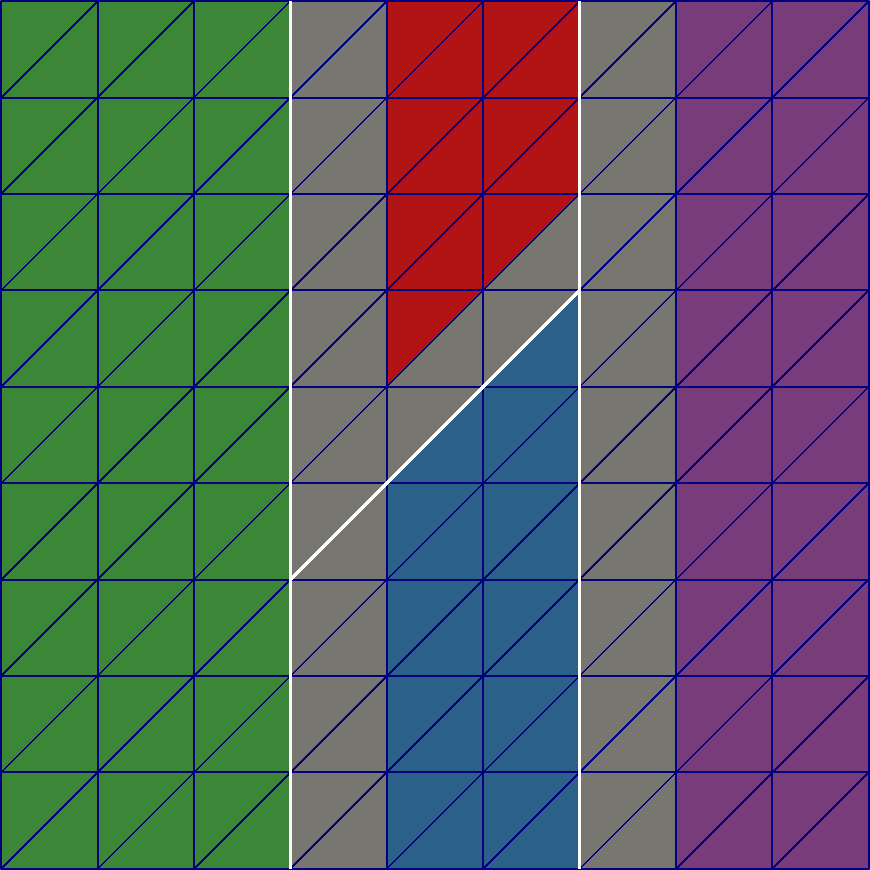}\label{subfig: double bad cut case 2}}
\caption{Domains and interface cut locations for the two condition number test cases (i): \protect\subref{subfig: simply bad cut case 1}, \protect\subref{subfig: simply bad cut case 2}  and (ii): \protect\subref{subfig: double bad cut case 1}, \protect\subref{subfig: double bad cut case 2}.}
 \label{fig: cuts condition number}
\end{figure}

In our first numerical example, we demonstrate that the ghost penalty terms introduced in \eqref{equ: ghost penalty terms} prevent ill conditioning of the system matrices associated with the linear stage of our LaTIn solver (recall equation \eqref{equ: bulk equation linear stage}). To fully control ``bad cuts" (see Figure \ref{fig: bad cut cases} for two typical bad cut cases), we design a crack test inspired by previous studies published in \cite{liuborja2008} and \cite{annavarapuhautefeuille2014}.
The computational domain of our test case is decomposed into 4 subdomains with two vertical and one diagonal ``crack" interface (see Figure~\ref{fig: cuts condition number}). ``Good" (large $K_i$) and ``bad" (small $K_i$) cut cases can be constructed through shifting these interfaces. 

\begin{figure}
[htbp]

\subfloat[$\epsilon_x=0.5$, $\epsilon=\epsilon_y$.]{\includegraphics[width=.5\textwidth]{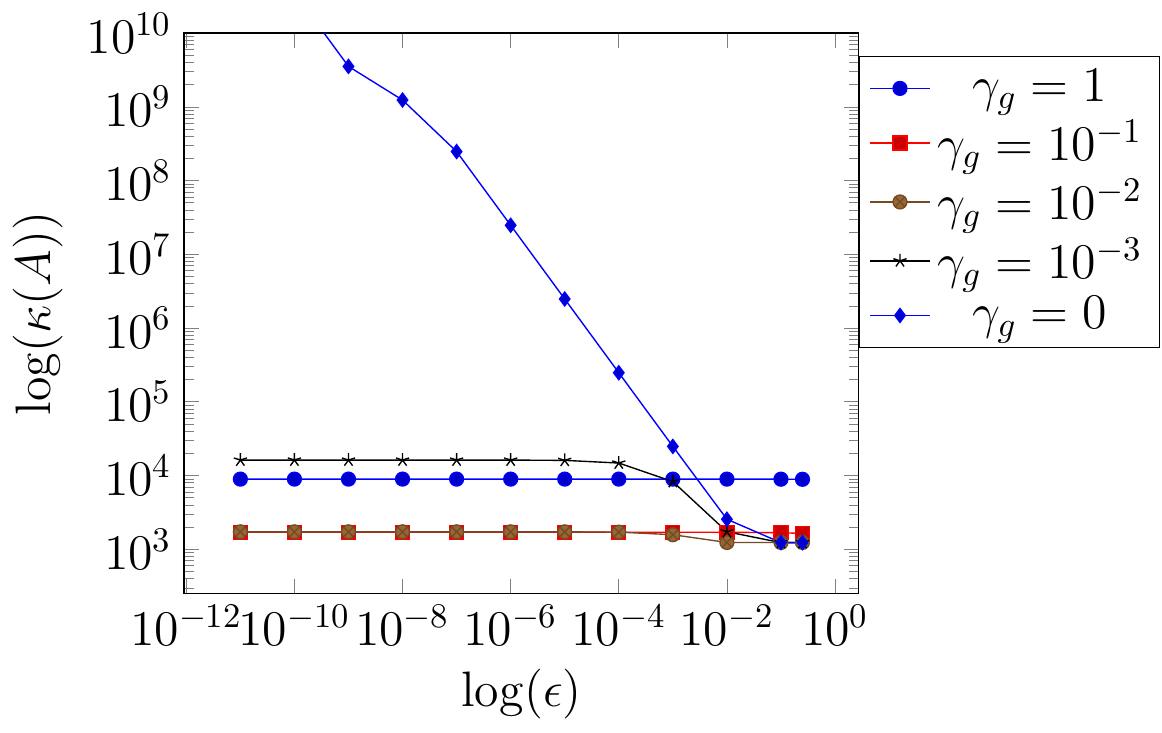}\label{subfig: condition number a}}
\subfloat[$\epsilon = \epsilon_x=\epsilon_y$.]{\includegraphics[width=.5\textwidth]{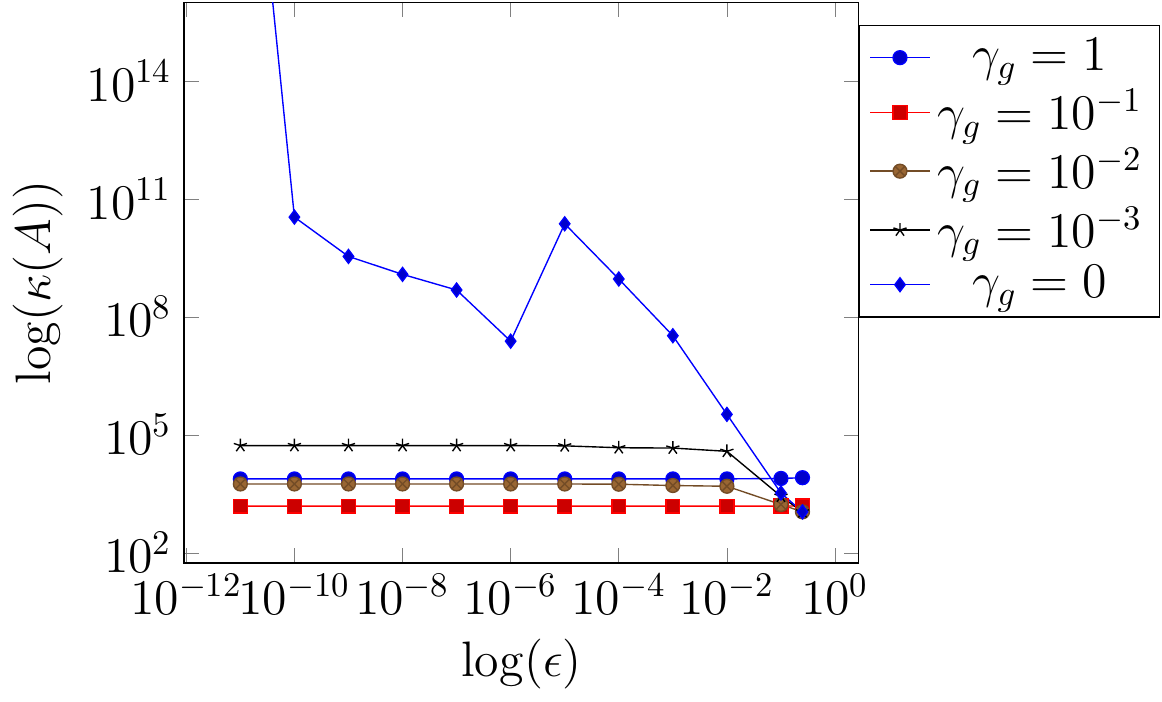}\label{subfig: condition number c}}
\caption{Condition number dependence on cut location for test case (i): \protect\subref{subfig: condition number a} and test case (ii): \protect\subref{subfig: condition number c}.}
 \label{fig: condition number}
\end{figure}

\begin{figure}
\centering
\includegraphics[width=.5\textwidth]{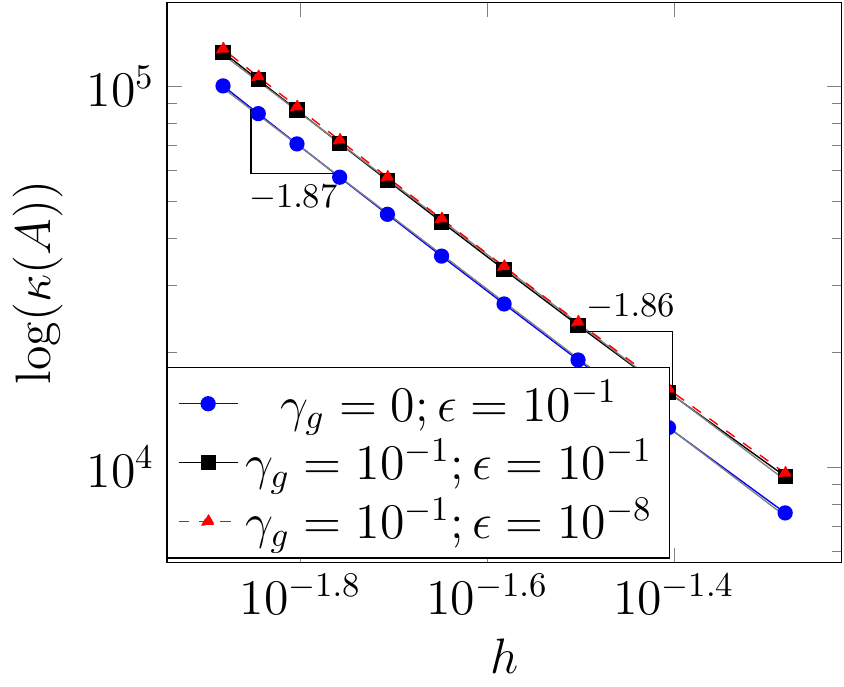}
\caption{Condition number dependence on mesh refinement.}
 \label{fig: condition number refinement}
\end{figure}

The geometry is defined by the following three level set functions
\begin{align}
\phi_1(x,y) &= y -x - \epsilon_y h \, , \\
\phi_2(x,y) &= x - \frac{1}{3} - \epsilon_x h \, , \\
\phi_3(x,y) &= \frac{2}{3} + \epsilon_x h - x \, ,
\end{align}
where the parameters $\epsilon_x$ and $\epsilon_y$ are chosen to range from $0.5$ (good cut) to $10^{-11}$ (bad cut). We investigate two test cases: (i)  the vertical interfaces are kept fixed at a good cut position with $\epsilon_x=0.5$ and the diagonal crack is moved in the interval $10^{-11}\leq\epsilon_y\leq 0.25$ (\textit{i.e.} ``simply bad" cut, see Figure~\ref{subfig: simply bad cut case 1}, \ref{subfig: simply bad cut case 2}) ; (ii) both $\epsilon_x$ and $\epsilon_y$ are varied in the interval $10^{-11} \leq \epsilon_x=\epsilon_y \leq 0.25$ (\textit{i.e.}  ``doubly bad" cut, see Figure~\ref{subfig: double bad cut case 1}, \ref{subfig: double bad cut case 2}).

Our LaTIn-based domain decomposition strategy leads to uncoupled sets of linear systems of equations. Therefore, we will only report the conditioning of the worst subproblem, which we denote by $\kappa$. Figure~\ref{fig: condition number} shows this condition number for test case (i) (see Figure \ref{subfig: condition number a}) and test case (ii) (see Figure \ref{subfig: condition number c}). In both cases, if the ghost penalty regularisation is switched off, \textit{e.g.} $\gamma_g=0$, then the condition number  increases drastically with decreasing $\epsilon$ (\textit{i.e.} worsening the cut scenario). This problem appears despite the regularisation effect of the LaTIn augmentation terms. This is because the search direction coefficient is chosen so as to scale with macroscopic properties of the problem at hand. In other words, this term is chosen to ensure the fast convergence of the LaTIn solver, and not to repair instabilities created by local bad cuts. However, with ghost penalty regularisation the increasingly unstable behaviour with decreasing $\epsilon$ is cured.  Even a very small value of the ghost penalty parameter, \textit{e.g.} $\gamma_g=10^{-3}$, is sufficient to guarantee a bounded condition number for very small $\epsilon$. 

Most importantly, the ghost penalty regularisation does not impact the scaling of the condition number of the system matrices with mesh refinement. This is demonstrated in Figure~\ref{fig: condition number refinement}. In fact, this scaling is typical for standard P1 finite element matrices. In other words, through the scaling that we have chosen for the ghost penalty regularisation, we ensure that the regularisation properties stated previously are consistent throughout the mesh refinement process.


\subsection{Convergence with mesh refinement}

We will now proceed to studying the convergence properties of the LaTIn-CutFEM algorithm, through three test cases: (i) a problem with a simple elliptical inclusion in contact with a matrix, (ii) a problem with two inclusions, which exhibits multiply connected points, and (iii) a problem that features a large elastic contrast between interacting elastic phases. 

\subsubsection{Elliptical inclusion}

\begin{figure}
[htbp]
\centering
\subfloat[Schematic and boundary conditions.]{
\includegraphics[width=.4\textwidth]{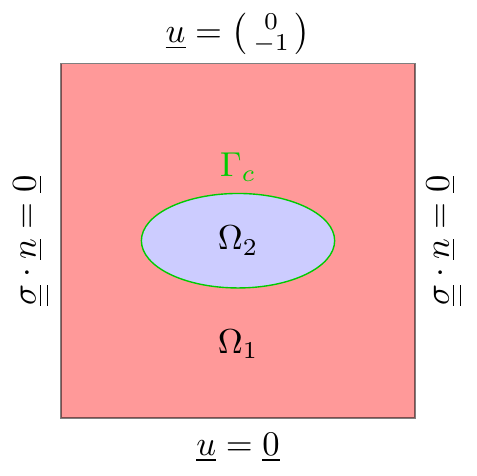} \label{subfig: schematic elliptical inclusion}}\hspace{.1cm}
\subfloat[Coarsest Mesh.]{
\includegraphics[width=.35\textwidth]{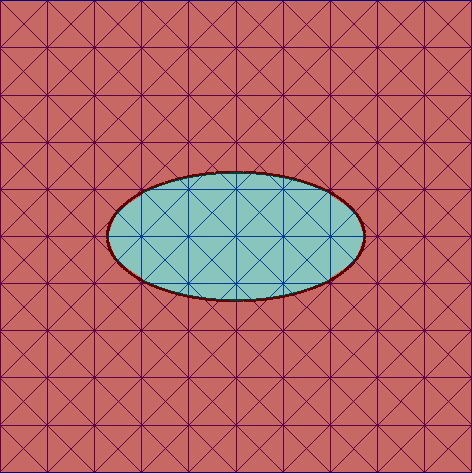}
 \label{subfig: mesh elliptical inclusion}
} \\
\subfloat[First hierarchical refinement.]{
\includegraphics[width=.35\textwidth]{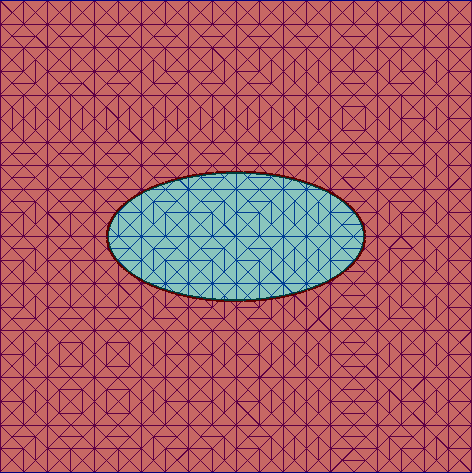}
 \label{subfig: mesh elliptical inclusion 1}
}
\subfloat[Finest mesh.]{
\includegraphics[width=.35\textwidth]{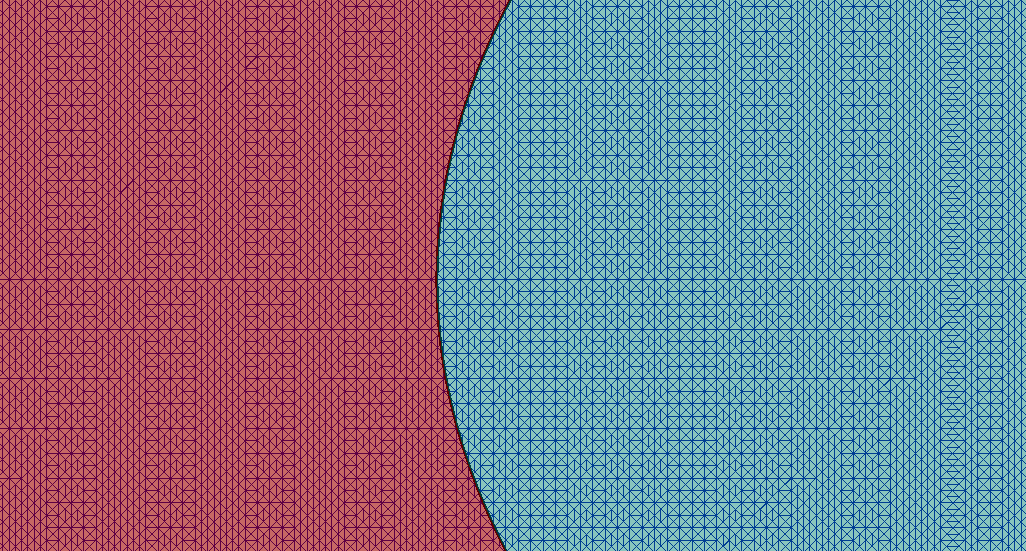}
 \label{subfig: mesh elliptical inclusion 2}
}
\caption{Schematic of elliptical inclusion with boundary conditions and regular background mesh refinement strategy.}
 \label{fig: schematic elliptical inclusion}
\end{figure}

\begin{figure}
[htbp]
\centering
\subfloat[Horizontal displacement $u_x$.]{\includegraphics[width=.4\textwidth]{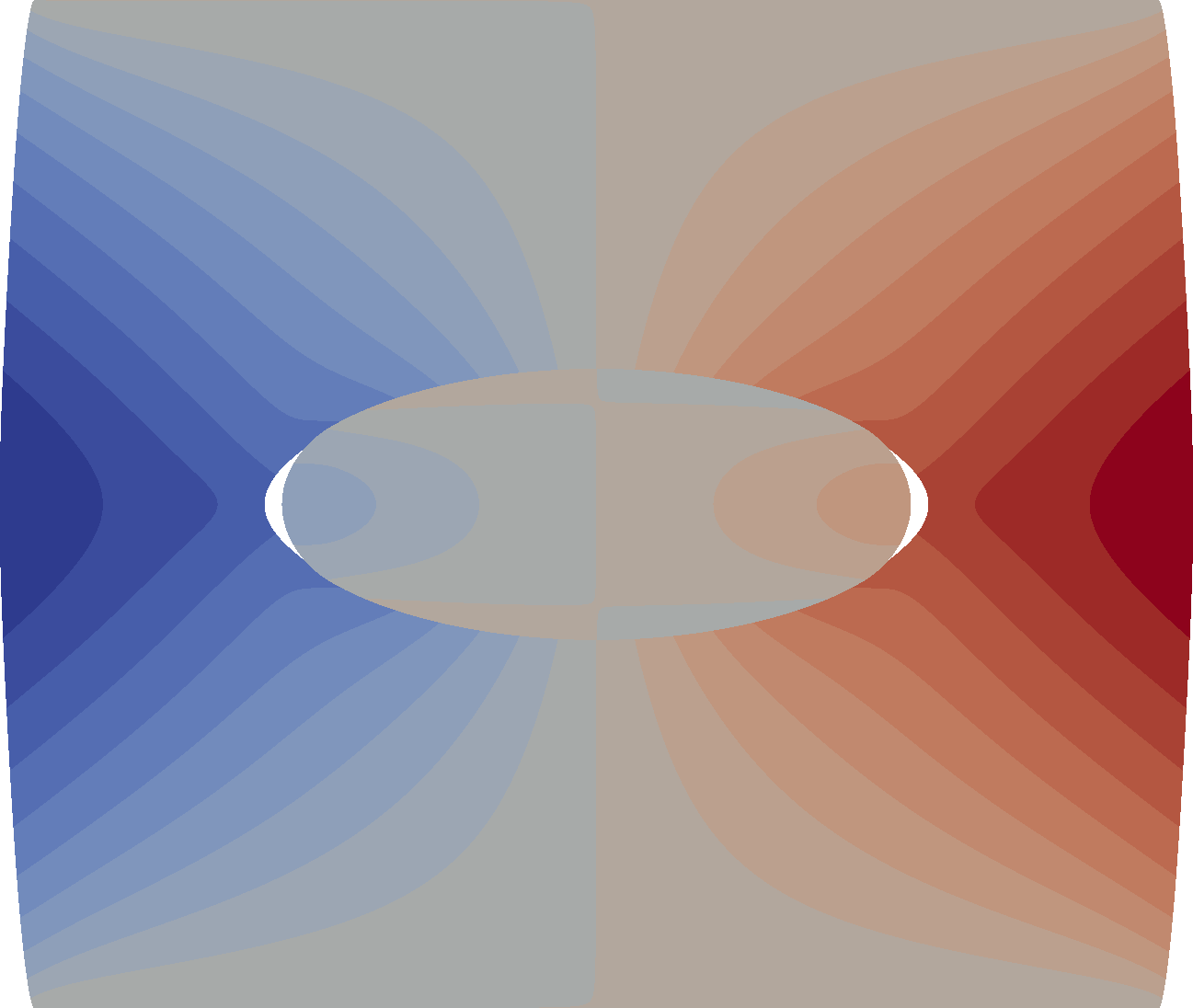}\,\,
\includegraphics[width=.1\textwidth]{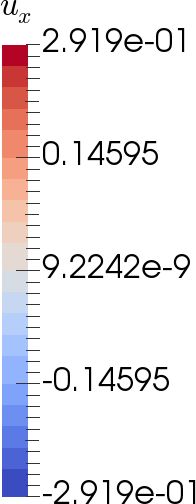}} 
\subfloat[Vertical displacement $u_y$.]{\includegraphics[width=.4\textwidth]{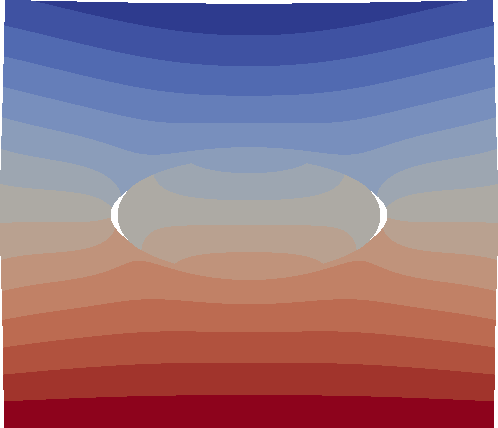}\,\,
\includegraphics[width=.1\textwidth]{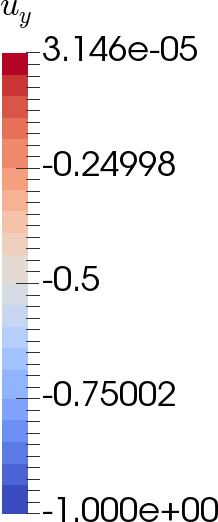}} \\
\subfloat[Stress component $\sigma_{xx}$.]{\includegraphics[width=.4\textwidth]{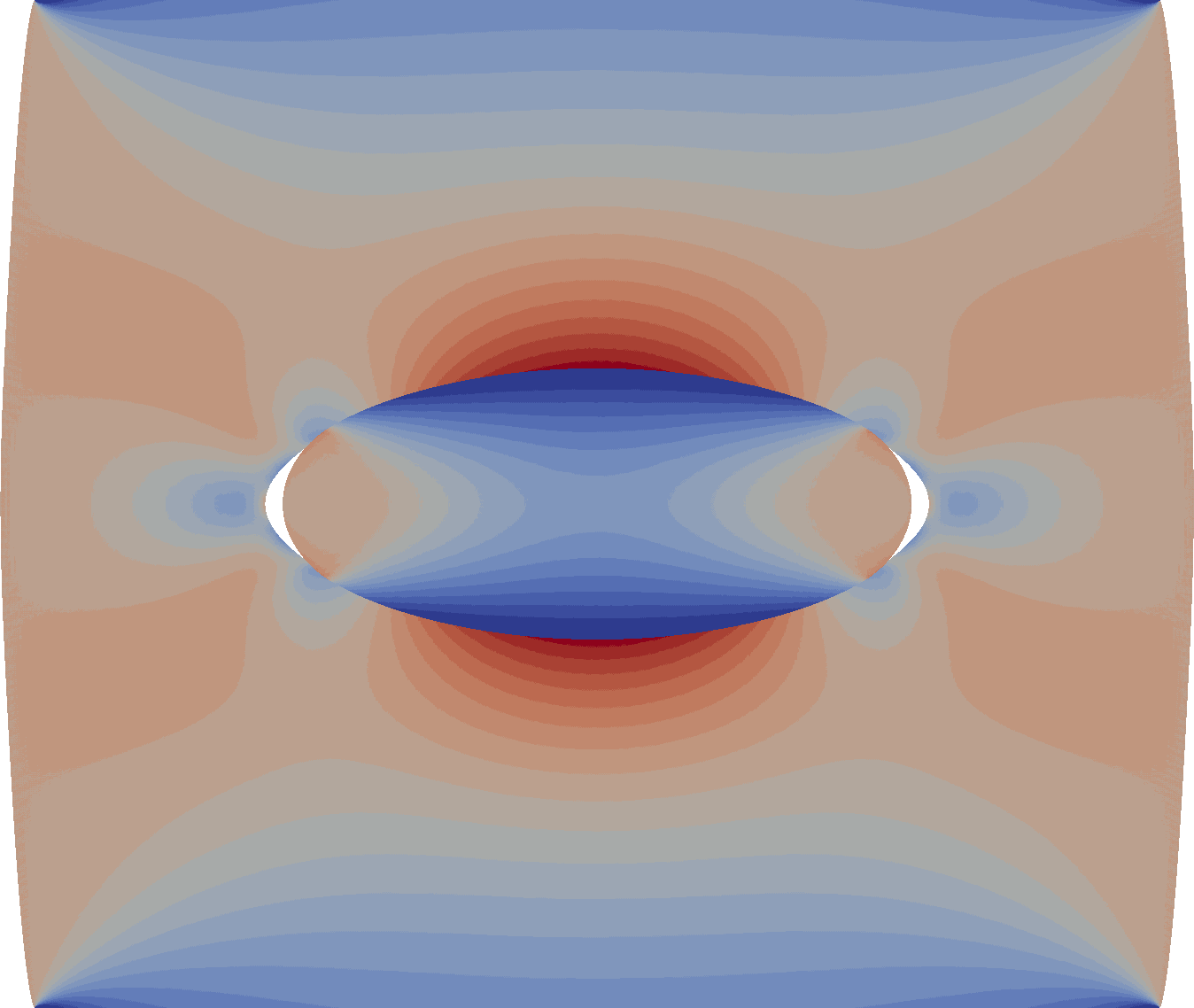}\,\,
\includegraphics[width=.1\textwidth]{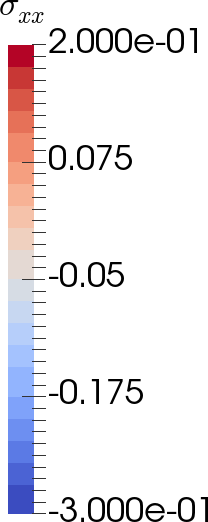}} 
\subfloat[Stress component $\sigma_{xy}$.]{\includegraphics[width=.4\textwidth]{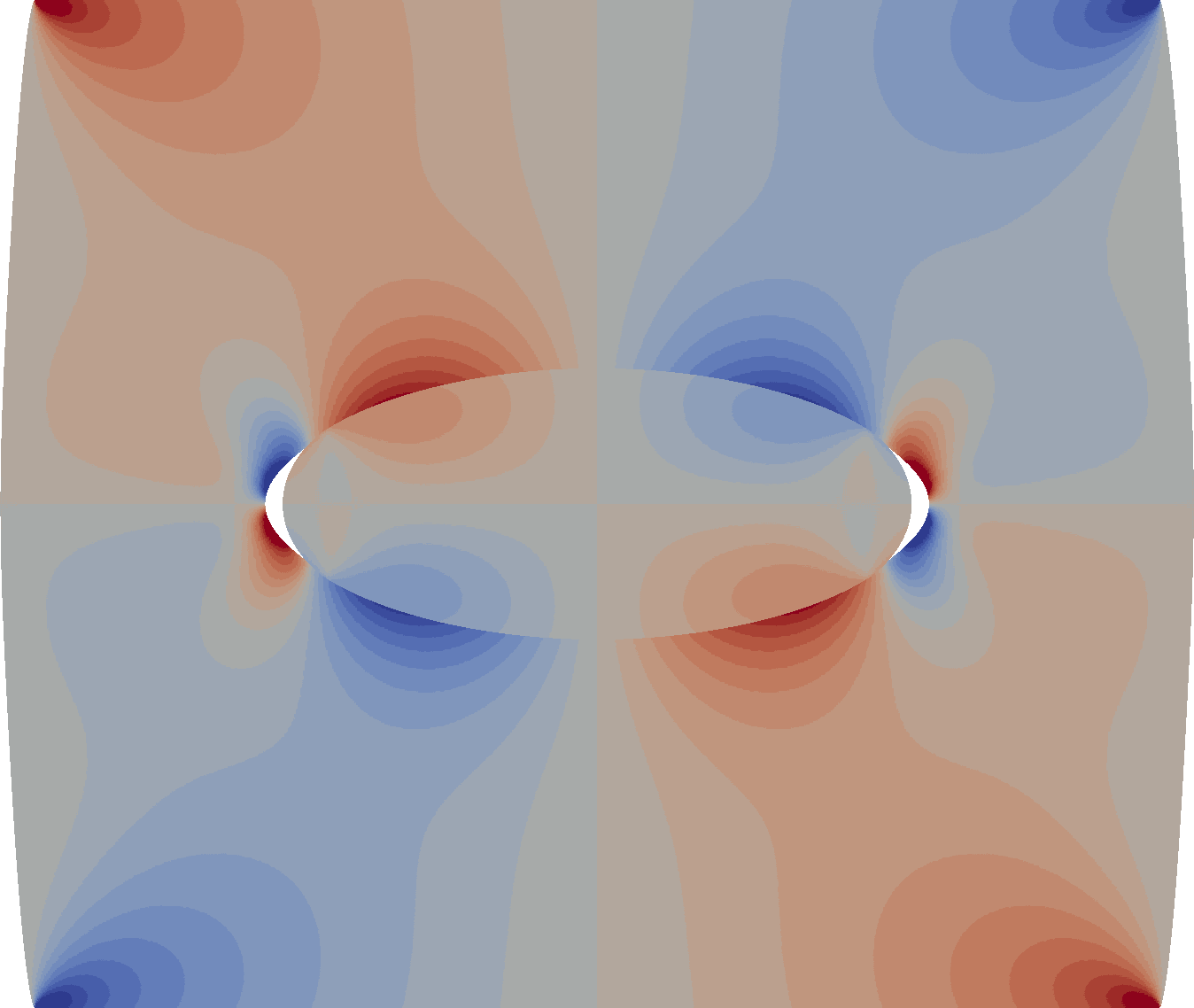}\,\,
\includegraphics[width=.1\textwidth]{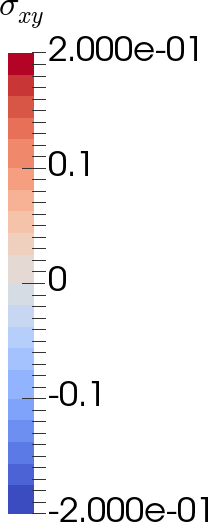}} \\
\subfloat[Stress component $\sigma_{yy}$.]{\includegraphics[width=.4\textwidth]{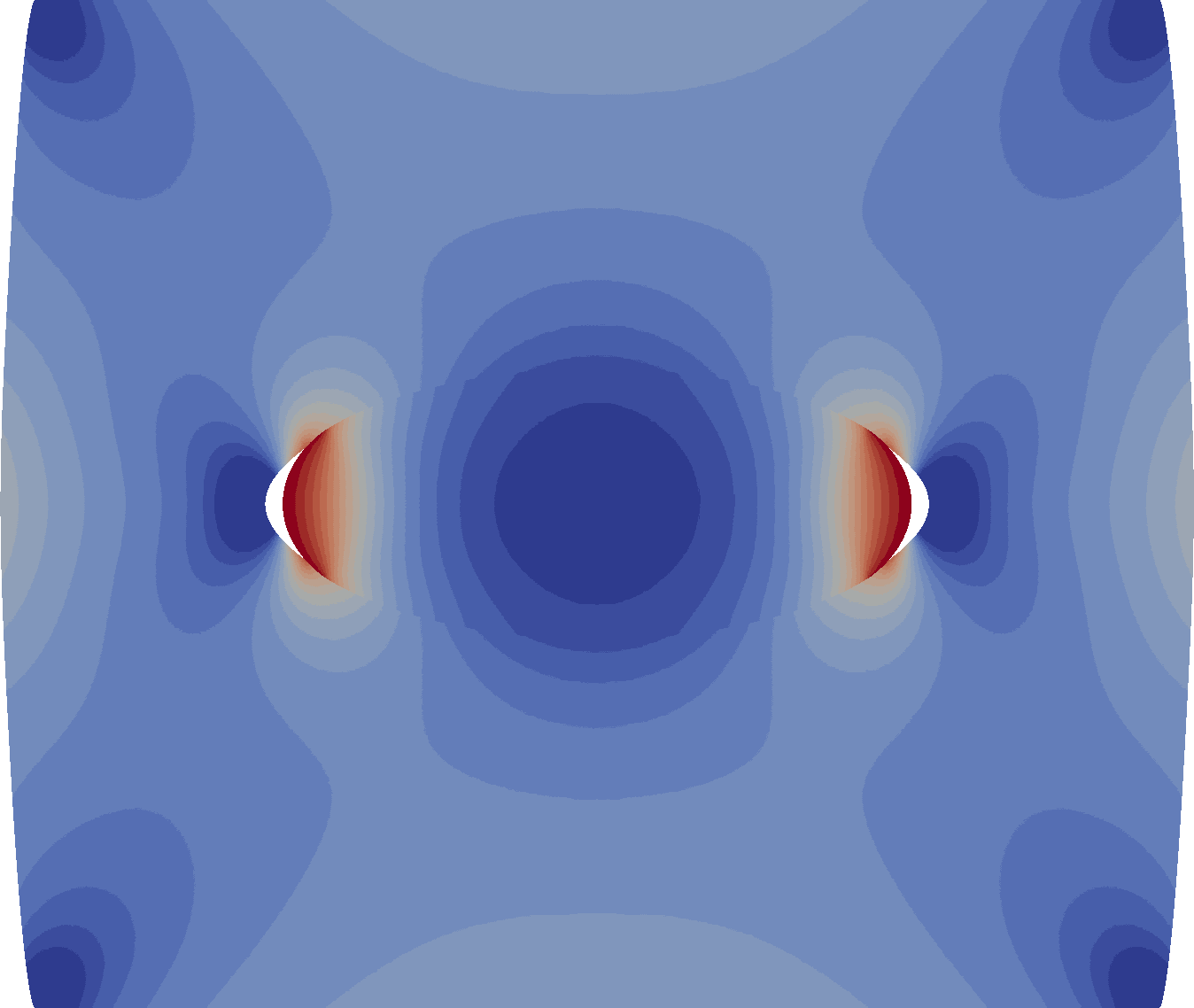}\,\,
\includegraphics[width=.1\textwidth]{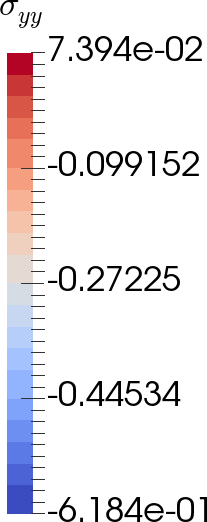}} 
\caption{Displacement and stress components of elliptical inclusion for the finest mesh.}
 \label{fig: elliptical inclusion unknowns}
\end{figure}

We start with our first convergence study. Let us consider an elliptical inclusion $\Omega_1$ in a rectangular domain $\Omega_2=\Omega \setminus \Omega_1$ with $\Omega = [-1.2,1.2] \times [-1.2,1.2]$. The elliptical inclusion is described by the level set function 
\begin{equation}
\phi_1(x,y)=\sqrt{\left(\frac{x}{a}\right)^2+\left(\frac{y}{b}\right)^2}- r,
\end{equation}
where $r=0.654545$, $a=1$, $b=0.5$. We apply a displacement of $\V{u} = \icol{0\\-1}$ on the top of the outer domain, a zero displacement at the bottom and zero Neumann conditions at the side of the domain
(see Figure~\ref{subfig: schematic elliptical inclusion}). We choose $E^1=E^2=1.0$, $k^+=k^-=1.0$, $\gamma_g=0.1$,  $\gamma_{\Pi}=0.1$ and $\alpha=10$.
The elliptical inclusion and the rectangular domain interact through unilateral contact interface $\Gamma_C$. Figure~\ref{fig: elliptical inclusion unknowns} shows the displacement and stress components for a fine mesh with $h=0.00375$ after $200$ LaTIn iterations (the LaTIn algorithm has reached a converged solution state, as shown later on). We observe that the inclusion and the background block material are not in contact on the left and right of the ellipse and are in contact on the top. There are stress concentrations where the contact boundary changes from ``in contact" to ``not in contact". We utilise this fine mesh solution as a numerical reference to investigate the convergence of our LaTin-CutFEM algorithm with mesh refinement and with the number of LaTin iterations. 

\begin{figure}
[htbp]
\centering
\subfloat[Energy norm error vs iterations and mesh size.]{\includegraphics[width=.6\textwidth]{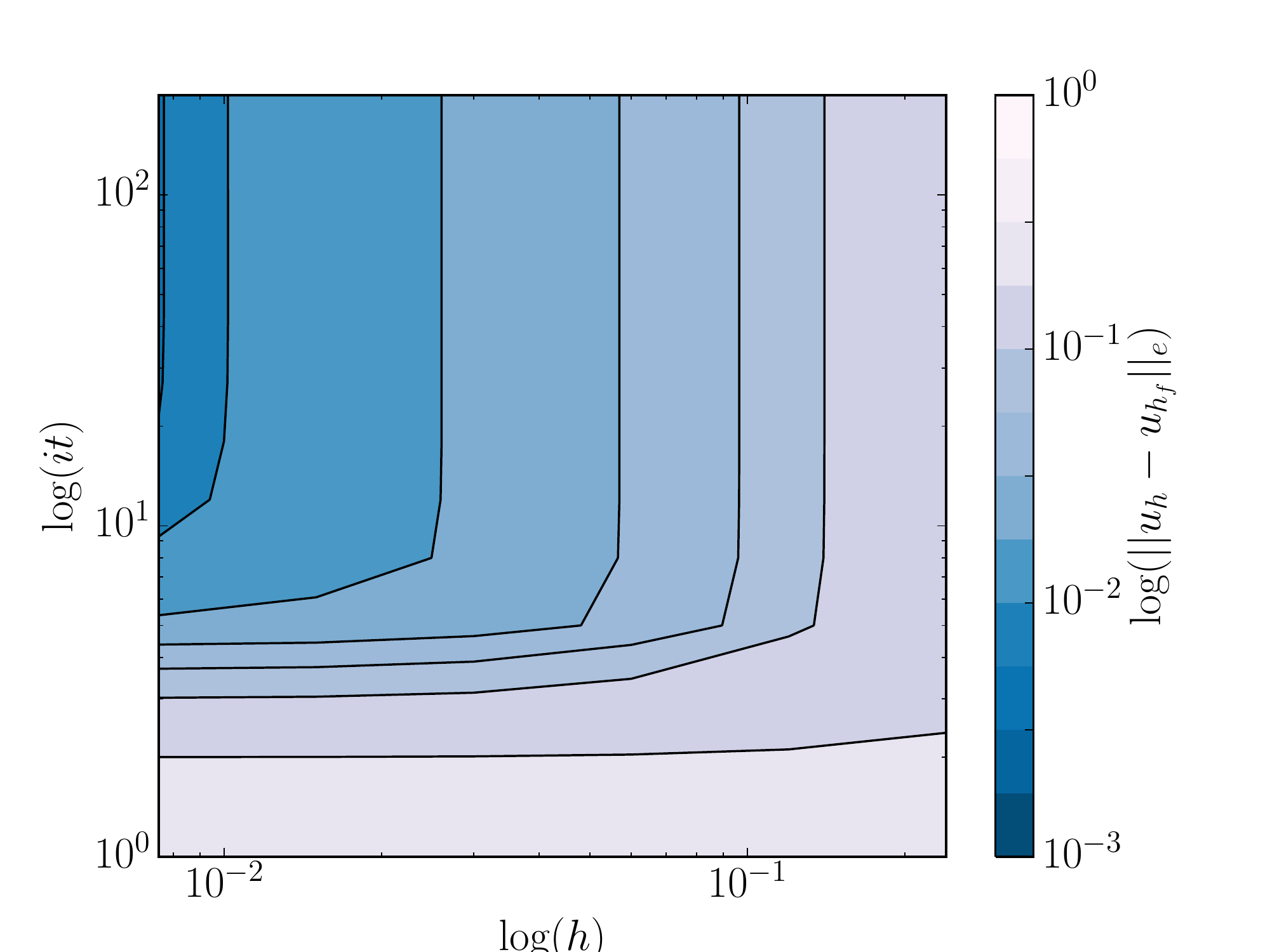}
\label{subfig: ellipse energy norm 3D}
} \\
\subfloat[Convergence rates for 2 quadrature points per interface segment.]{\includegraphics[width=.4\textwidth]{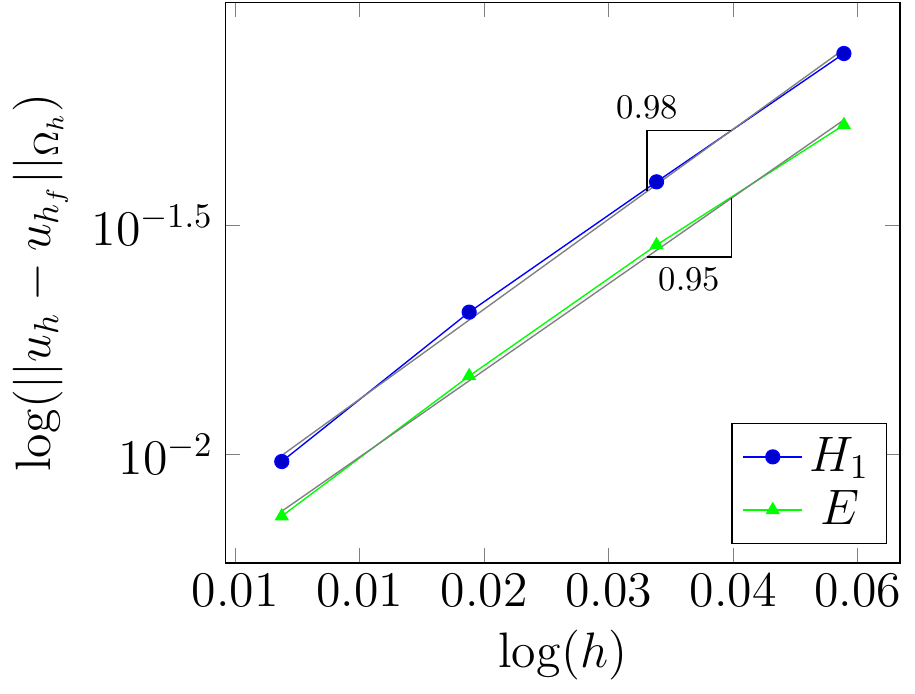}
\label{subfig: ellipse convergence rate 2}
} \,\,
\subfloat[Convergence rates for 4 quadrature points per interface segment.]{\includegraphics[width=.4\textwidth]{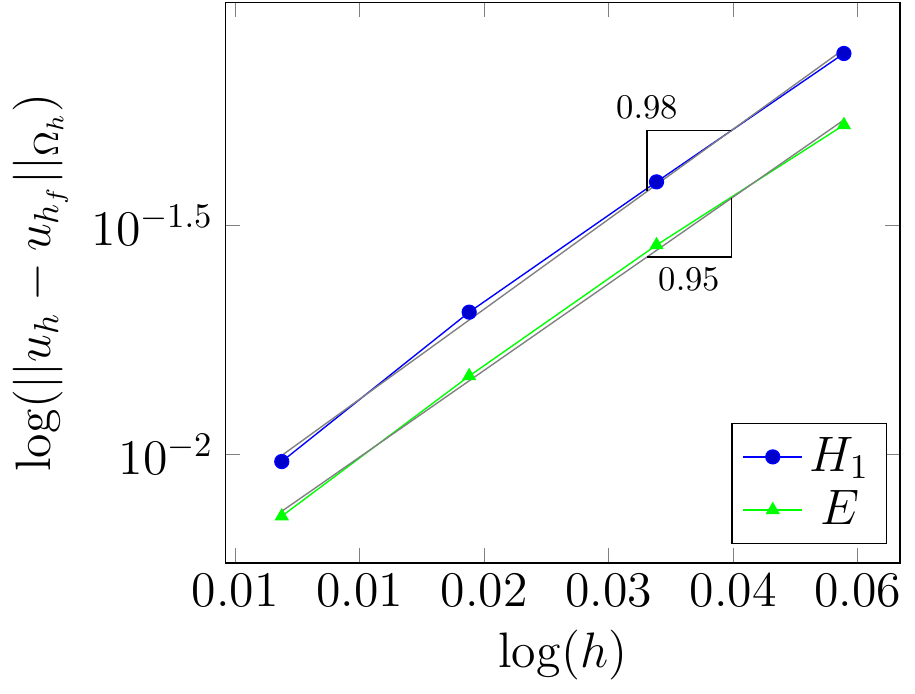}
\label{subfig: ellipse convergence rate 4}
}\\
\subfloat[Energy norm error with iteration count for coarser mesh sizes.]{\includegraphics[width=.5\textwidth]{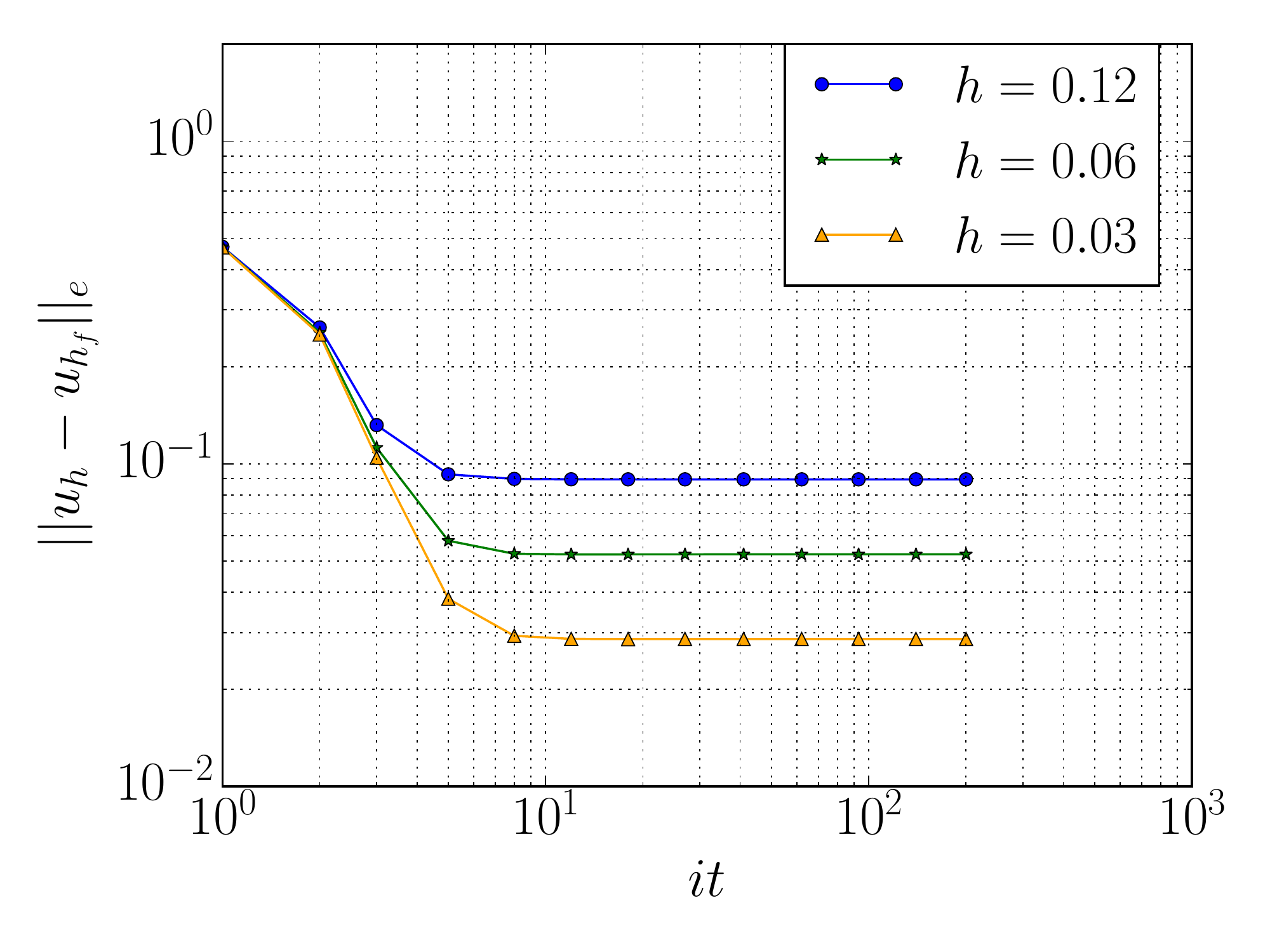}
\label{subfig: ellipse energy norm it coarse}}
\subfloat[Energy norm error with iteration count for finest solution.]{\includegraphics[width=.5\textwidth]{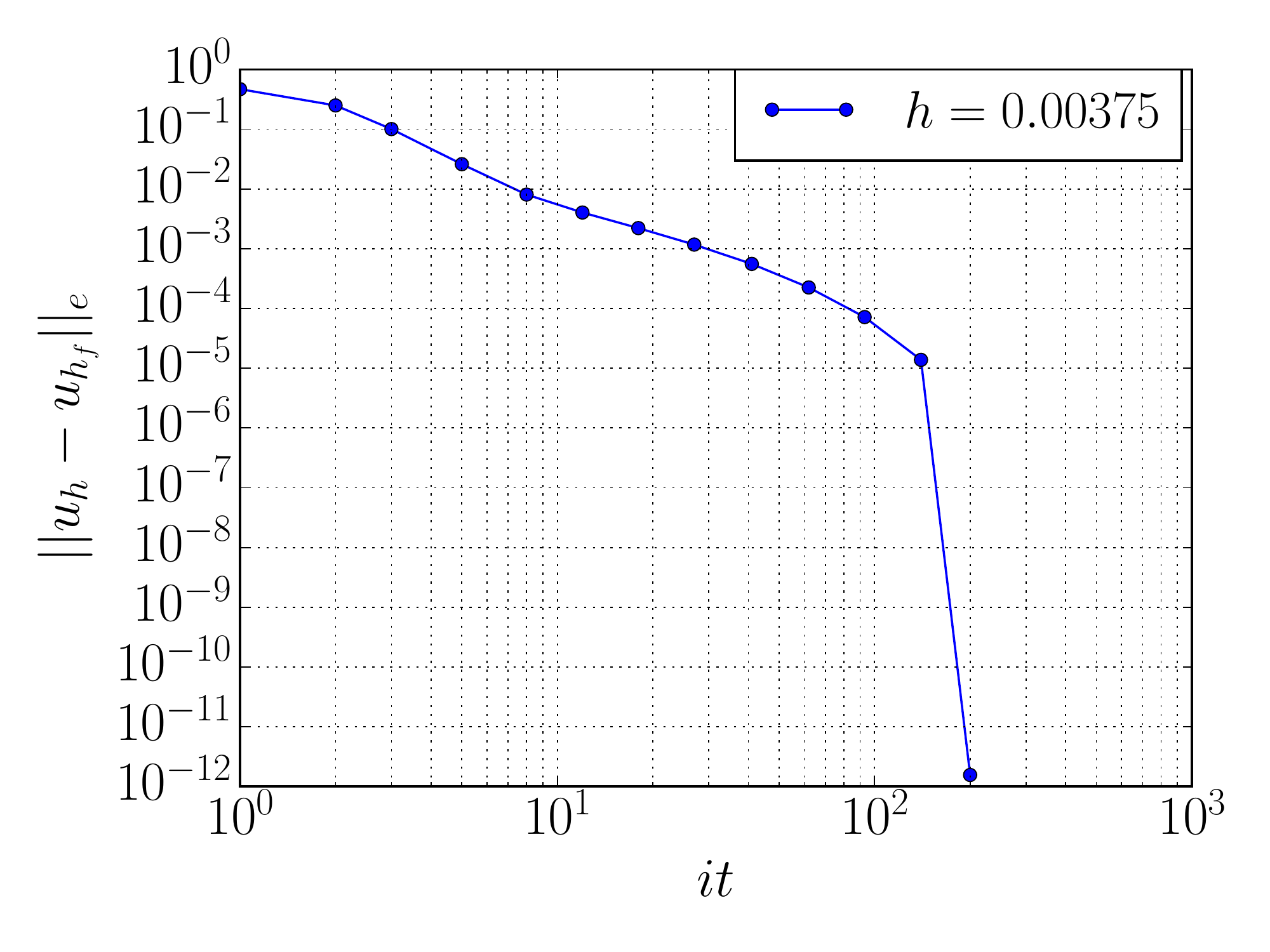}
\label{subfig: ellipse energy norm it fine}}
\caption{Convergence rates and energy norm error with iteration count for ellipse-shaped inclusion.}
 \label{fig: elliptical inclusion convergence}
\end{figure}


\begin{table}
\centering
\subfloat[Two quadrature points.]
{\begin{filecontents*}{./pgfplots/convergence/circle-convergence/nquad2/u_R7_it200_error_table.txt}
N		hmax		L2		H1		E
3		0.06		0.00304983		0.0750873		0.0524061
4		0.03		0.000943873		0.0393787		0.0286505
5		0.015		0.000211899		0.0204461		0.0148274
6		0.0075		5.7897e-05		0.00964451		0.00732947
\end{filecontents*}
\pgfplotstabletypeset[
every head row/.style={before row=\hline,after row=\hline\hline},
every last row/.style={after row=\hline},
every first column/.style={
column type/.add={|}{}
},
every last column/.style={
column type/.add={}{|}
},
columns={hmax,H1,E},
columns/hmax/.style={
column name=$h$,
sci,sci zerofill,sci sep align,precision=1,sci 10e},
columns/H1/.style={
column name=$H_1$-error,
sci,sci zerofill,sci sep align,precision=4,sci 10e},
columns/E/.style={
column name=$E$-error,
sci,sci zerofill,sci sep align,precision=4,sci 10e}
]
{./pgfplots/convergence/circle-convergence/nquad2/u_R7_it200_error_table.txt}
}

\subfloat[Four quadrature points.]
{
\begin{filecontents*}{./pgfplots/convergence/circle-convergence/nquad4/u_R7_it200_error_table.txt}
N		hmax		L2		H1		E
3		0.06		0.00299211		0.0751062		0.0524135
4		0.03		0.000944123		0.0393794		0.0286507
5		0.015		0.000211864		0.0204461		0.0148274
6		0.0075		5.76218e-05		0.00964443		0.00732944

\end{filecontents*}
\pgfplotstabletypeset[
every head row/.style={before row=\hline,after row=\hline\hline},
every last row/.style={after row=\hline},
every first column/.style={
column type/.add={|}{}
},
every last column/.style={
column type/.add={}{|}
},
columns={hmax,H1,E},
columns/hmax/.style={
column name=$h$,
sci,sci zerofill,sci sep align,precision=1,sci 10e},
columns/H1/.style={
column name=$H_1$-error,
sci,sci zerofill,sci sep align,precision=4,sci 10e},
columns/E/.style={
column name=$E$-error,
sci,sci zerofill,sci sep align,precision=4,sci 10e}
]
{./pgfplots/convergence/circle-convergence/nquad4/u_R7_it200_error_table.txt}
}
\caption{$H_1$ and energy norm error values for two and four quadrature points per interface segment of elliptical inclusion problem.}
\label{tab: quad points elliptical}
\end{table}

\begin{figure}
[htbp]
\centering
\subfloat[P1/P0]{\includegraphics[width=.33\textwidth]{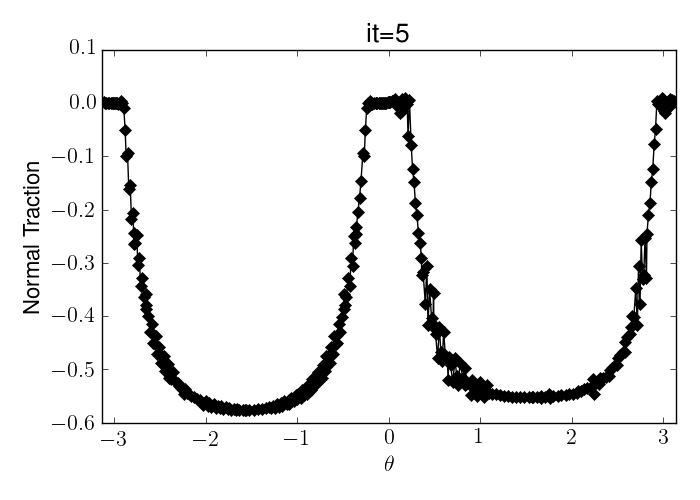}}
\subfloat[P1/P0]{\includegraphics[width=.33\textwidth]{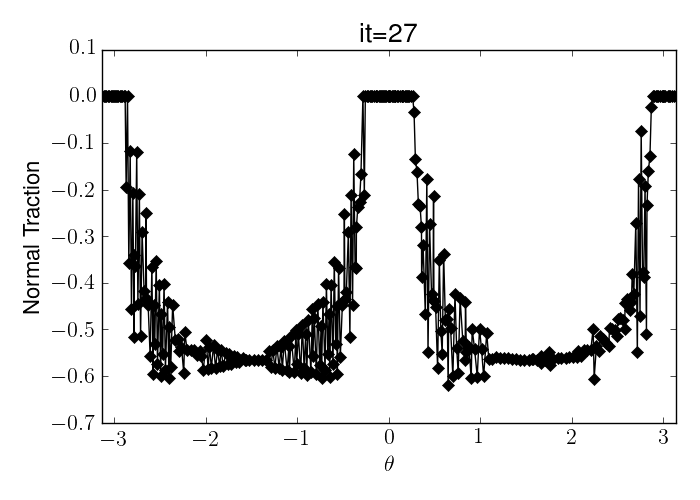}}
\subfloat[P1/P0]{\includegraphics[width=.33\textwidth]{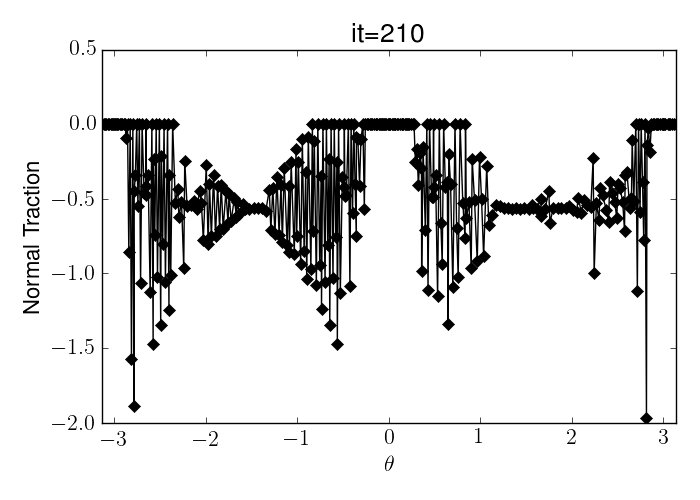}} \\ 
\subfloat[P1/P1]{\includegraphics[width=.33\textwidth]{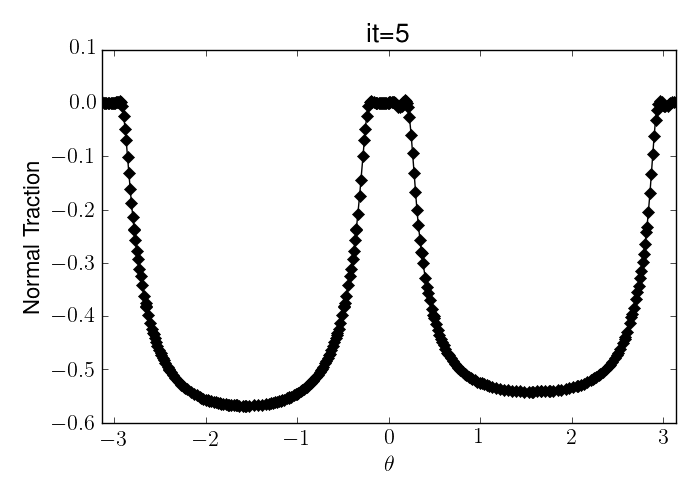}}
\subfloat[P1/P1]{\includegraphics[width=.33\textwidth]{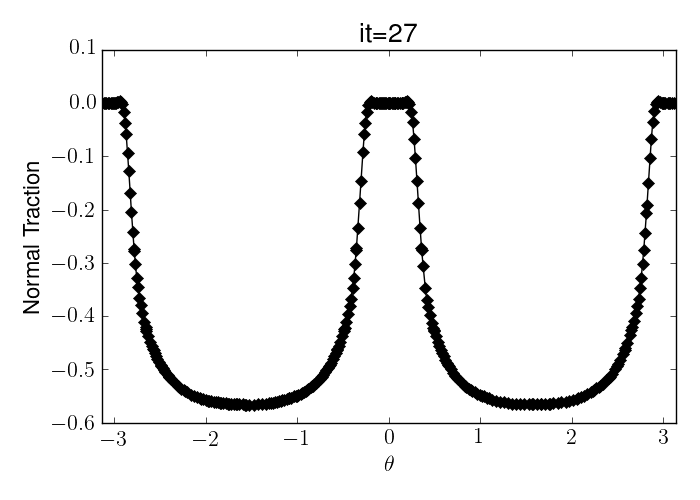}}
\subfloat[P1/P1]{\includegraphics[width=.33\textwidth]{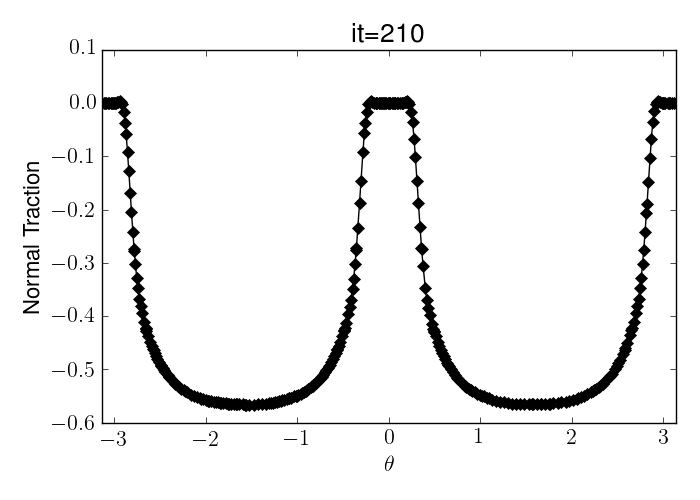}} 
\caption{Comparison of current stable P1/P1 discretisation versus unstable P1/P0 discretisation for it$=\{5,27,210\}$ LaTIn iterations.}
\label{fig: P1P1 vs P1P0}
\end{figure}

\paragraph{Convergence of algorithmically converged solutions with mesh refinement.} The coarsest background mesh is the regular mesh shown in Figure~\ref{subfig: mesh elliptical inclusion}. This mesh is refined in a hierarchical manner (each triangle of the coarser mesh is subdivided into 4 triangles) to produce the finer levels of discretisation. We perform 200 iterations of the LaTIn solver for each of these meshes. We now investigate the rate of convergence of the numerical solution in the $H_1$-norm, which is defined as
\begin{equation}
\norm{\Disp_h-\Disp_{h_f}}_1 = \sqrt{\sum_{i \in \mathcal{I}_{\Domain}} \int_{\Domain^i_{h_f}} (( \Disp^i_{h} - \Disp^i_{h_f})^2 + (\nabla \Disp^i_{h} - \nabla \Disp^i_{h_f})^2 ) d\Domain}
\end{equation} 
and in the energy norm, which reads 
\begin{equation}
\norm{\Disp_h-\Disp_{h_f}}_E  = 
\sqrt{ \sum_{i \in \mathcal{I}_{\Domain}} \int_{\Domain^i_{h_f}} \StrainBasic(\Disp_h-\Disp_{h_f}) : \Hooke : \StrainBasic(\Disp_h-\Disp_{h_f}) \, d\Domain } \, ,
\end{equation}
where $\Disp_{h_f}$ denotes the reference solution (\textit{i.e.} finest mesh). To evaluate these norms, we interpolate the solution on the coarser mesh into the finest mesh and then evaluate the error norm over the linear geometry approximation, $\Domain^i_{h_f}$, of the finest mesh. Figure~\ref{subfig: ellipse convergence rate 2}, \ref{subfig: ellipse convergence rate 4} shows that the convergence rate of our scheme is $0.98$ in the $H_1$ and $0.95$ in the energy norm. This is optimal, in the sense that the proposed mixed numerical solver for unilateral contact uses P1 standard finite elements to solve the bulk equations and the optimal convergence rate in the $H^1$ and energy norm is therefore expected to be one. 

\paragraph{Effect of the quadrature order for the multiscale filter.} The ``heart" interface fields are evaluated at a discrete set of points. This is required to integrate the right-hand side of projection \eqref{eq:ContinuityLocal}. We set those points to be the quadrature points of a standard Gauss quadrature rule, applied to each of the line segments (2D) or triangles (3D) of the subtessalation associated with the embedded contact interfaces.
The purpose of the present study is to choose an appropriately large quadrature order for the remainder of our work. 
Considering Figure~\ref{subfig: ellipse convergence rate 2}, \ref{subfig: ellipse convergence rate 4} and Table~\ref{tab: quad points elliptical}, we see that switching from 2 points (\textit{i.e.} Gauss quadrature of third order) to 4 (\textit{i.e.} order 7) leads to only marginal changes in the error norms and to no measurable change in the overall convergence order of the LaTIn-CutFEM algorithm. In other words, order 3 is sufficiently accurate to yield results that cannot be distinguished from that obtained when using an overkill integrator. Therefore, we choose to use order 3 quadrature rules in all our numerical examples. We emphasise that this study indicates that the stable behaviour of the proposed P1-P1 LaTIn scheme is not due to some under-integration of the contact laws, but to the two-scale nature of the (``non-localised") local stage.

\paragraph{Convergence of the LaTIn algorithm.}
Next, we show that the LaTIn solver converges well with the number of LaTIn iterations. We investigate the evolution of the error in energy norm with LaTIn iterations for a fixed coarser mesh $h>h_f$. More precisely, we measure the error of our numerical solution at LaTIn iteration $it<it_{max}$ with respect to the finest solution at the maximum iteration $it_{max}:=200$. The results plotted in Figure~\ref{subfig: ellipse energy norm it coarse} indicate the LaTin-CutFEM algorithm reaches an algorithmically converged solution at $it\approx 20$, which is indicated by the flat part of the error curves. We can also observe that the error in the solutions remains constant after $it\approx 20$, \textit{i.e.} the remainder is the discretisation error. 

We emphasise that these results are non-classical (as far as we are aware, such convergence results for LaTIn-based solvers have never been published). Most contributions only report the convergence of the algorithmic error of LaTIn solvers. For our finest mesh, the algorithmic error is the one reported in Figure~\ref{subfig: ellipse energy norm it fine}, as the reference solution mesh and that of the approximate one are the same. However, the error curves of Figure~\ref{subfig: ellipse energy norm it coarse} represent the \textit{combined effect of algorithmic and discretisation errors}. When seen from this perspective, the convergence of the LaTIn solver is faster than usually thought, as the algorithmic error quickly drops to the level of the discretisation error, which can only be overcome by mesh refinement. As a corollary, the algorithmic error alone decreases slowly (this has been acknowledged in several contributions), but there is in fact no reason to make it converge to level of errors that are below the level of pure discretisation error.

\paragraph{Convergence of the LaTIn algorithm for P1/P0 approximation versus P1/P1 approximation.} 
Figure~\ref{fig: P1P1 vs P1P0} shows the comparison of our stable P1/P1 LaTin-CutFEM discretisation with a P1/P0 scheme. In the latter case, the bulk displacement is approximated by continuous piecewise linear polynomials but both interface fields are approximated as piecewise constant (the continuity is enforced in a weak sense, by enforcing that the projection of the corresponding residual in the space of piecewise constant interface fields vanishes). Figure~\ref{fig: P1P1 vs P1P0} shows the normal traction, $\InterfaceForce^{2,1}_h \cdot \Normal^{2,1}$ , along the elliptical contact interface $\Interface_C$ in terms of the angle $\theta = \arctan{(y/x)}$. We can see that with increasing LaTIn iterations the solution for the P1/P1 scheme has reached a converged state at $it=27$ and the solution remains virtually unchanged at $it=210$. However, for the P1/P0 scheme, we observe an onset of oscillations in the normal traction at $it=5$, which grow substantially with LaTIn iterations. At $it=27$, we see substantial pollution of the normal tractions through these oscillations and at $it=210$, these oscillation completely destroy the normal tractions.

\subsubsection{Multiply-connected interface points: effect of simultaneous enrichments}

\begin{figure}
[htbp]
\centering
\subfloat[Schematic.]{\includegraphics[width=.4\textwidth]{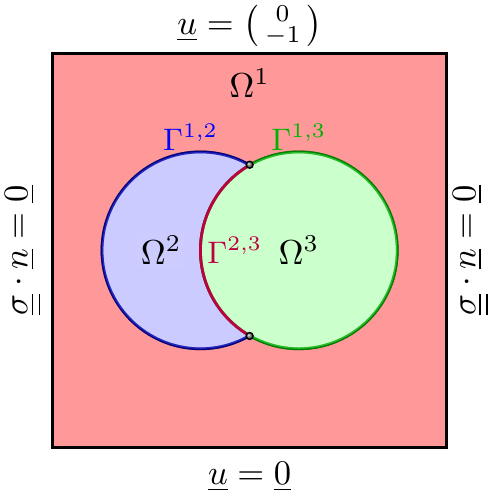}
\label{subfig: schematic two inclusions}}
\subfloat[Coarsest mesh $h=0.26$.]{\includegraphics[width=.35\textwidth]{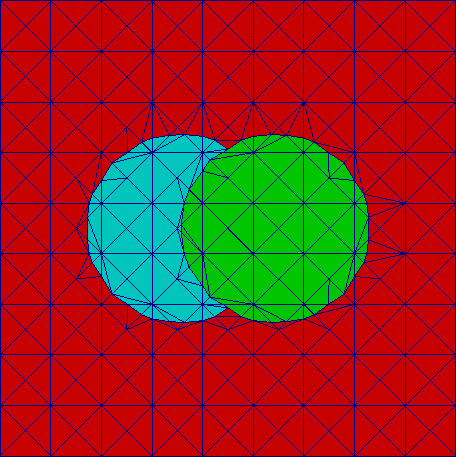}
\label{subfig: schematic two inclusions mesh}} 
\caption{Schematic of two inclusions with boundary conditions and piecewise linear approximation of the geometry in coarsest mesh used.}
 \label{fig: schematic two inclusions}
\end{figure}

\begin{figure}
[htbp]
\subfloat[$\Disp_y$]{\includegraphics[width=.4\textwidth]{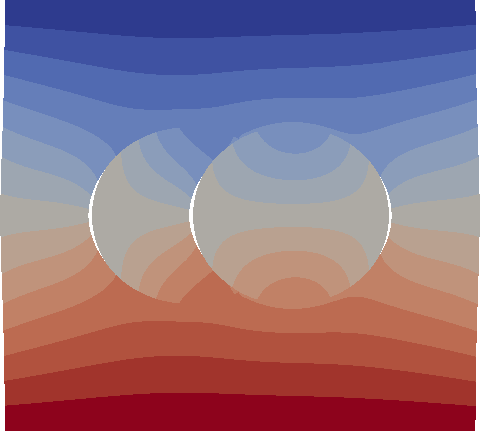}\,\,
\includegraphics[width=.1\textwidth]{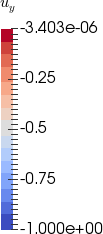}}
\subfloat[$\Disp_y$]{\includegraphics[width=.4\textwidth]{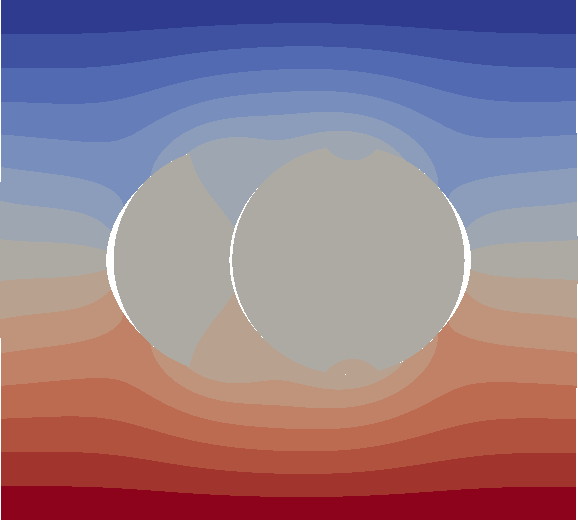}\,\,
\includegraphics[width=.1\textwidth]{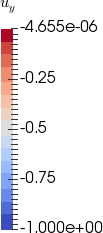}} \\
\subfloat[$\sigma_{xy}$]{\includegraphics[width=.4\textwidth]{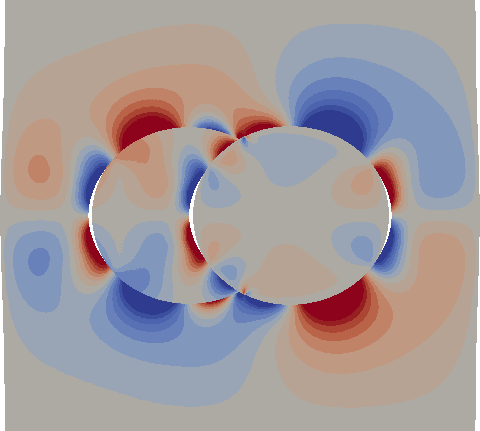}\,\,
\includegraphics[width=.1\textwidth]{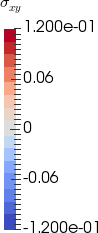}} 
\subfloat[$\sigma_{xy}$]{\includegraphics[width=.4\textwidth]{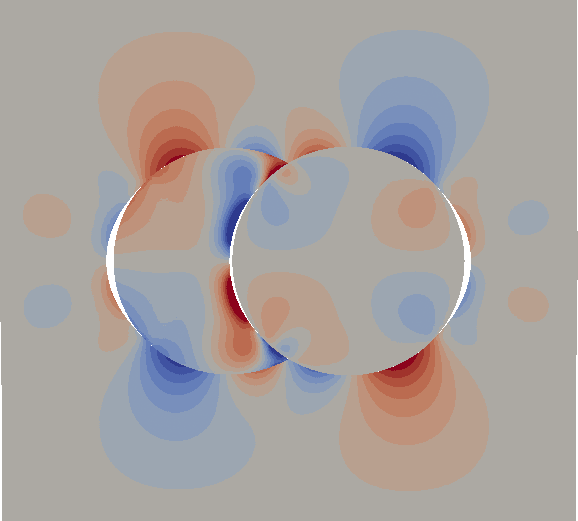}\,\,
\includegraphics[width=.1\textwidth]{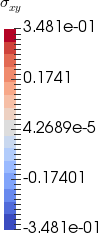}} \\
\subfloat[$\sigma_{yy}$]{\includegraphics[width=.4\textwidth]{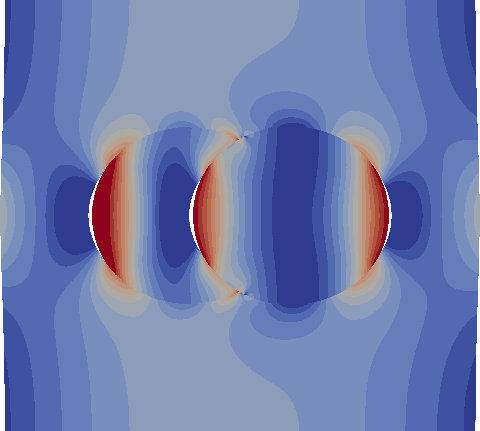}\,\,
\includegraphics[width=.1\textwidth]{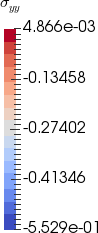}} 
\subfloat[$\sigma_{yy}$]{\includegraphics[width=.4\textwidth]{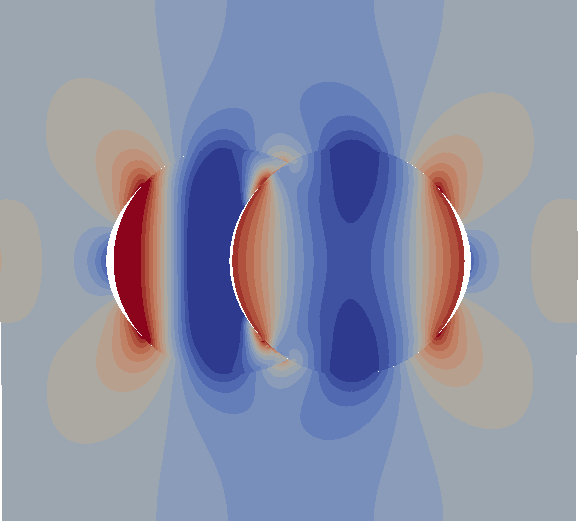}\,\,
\includegraphics[width=.1\textwidth]{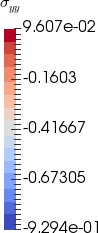}} 
\caption{Displacement and stress components of two inclusion for finest mesh with $E^1=E^2=E^3=1$ (left) and $E^1=1, E^2=E^3=10$ (right).}
 \label{fig: two inclusions unknowns}
\end{figure}

In our next example, we consider two intersecting inclusions in a rectangular domain $\Omega = [-1.2,1.2] \times [-1.2,1.2]$, which is the example referenced throughout this article. The two circular inclusions are described by the level set functions 
\begin{align}
\phi_1(x,y)&=\sqrt{(x-x_M^1)^2+y^2} - r, \\
\phi_2(x,y)&=\sqrt{(x-x_M^2)^2+y^2} - r, 
\end{align}
where $r=0.5$, $x_M^1=-0.25$ and $x_M^2=0.25$. We choose $k^+=k^-=1.0$, $\gamma_g=0.1$,  $\gamma_{\Pi}=0.1$ and $\alpha=10$. The resulting geometry has two triple junctions which will yield high stress concentrations at these points. We apply a displacement of $\V{u} = \icol{0\\-1}$ at the top boundary of the rectangular domain, a zero displacement at the bottom and zero Neumann conditions at the side of the domain (see Figure~\ref{subfig: schematic two inclusions}). Figure~\ref{fig: two inclusions unknowns} shows the vertical displacement, the shear stress $\sigma_{xy}$ and the normal stress $\sigma_{yy}$ for the two inclusion problem with no contrast in $E^1=E^2=E^3=1$ on the left and with a contrast of $E^1=1$ and $E^2=E^3=10$ on the right. In the case of stiffer inclusions than the matrix material, the inclusions deviate only slightly from their circular shape while in the case of no contrast the inclusions undergo a much stronger deformation. 

\begin{figure}
[htbp]
\centering
\subfloat[Energy norm error vs iterations and mesh size.]{\includegraphics[width=.6\textwidth]{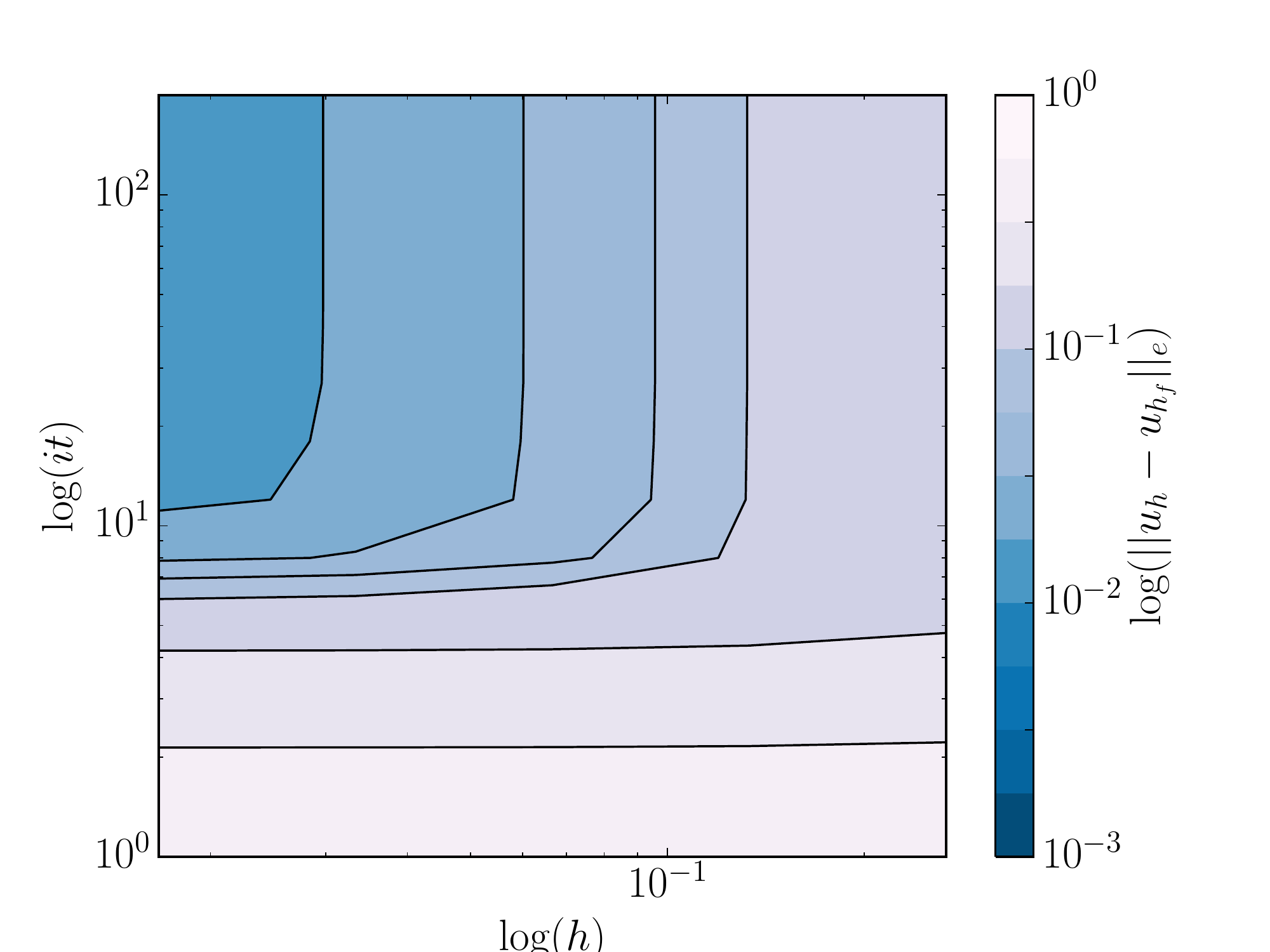}
\label{subfig: two inclusions convergence 3D}} 
\subfloat[Convergence rates.]{\includegraphics[width=.4\textwidth]{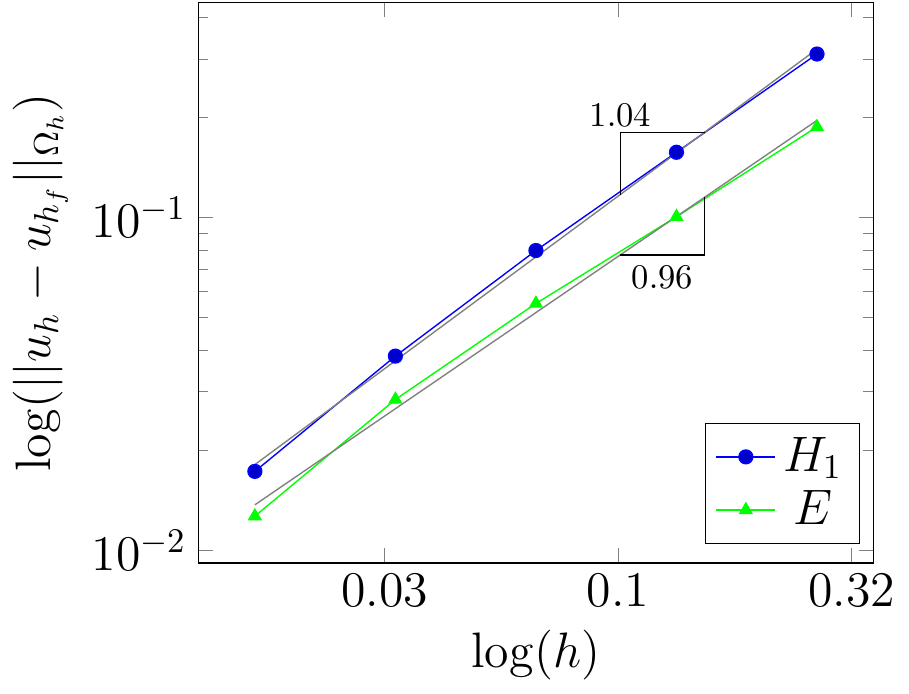}
\label{subfig: two inclusions convergence rate}} \\
\subfloat[Energy norm error with iteration count for coarser mesh sizes.]{\includegraphics[width=.5\textwidth]{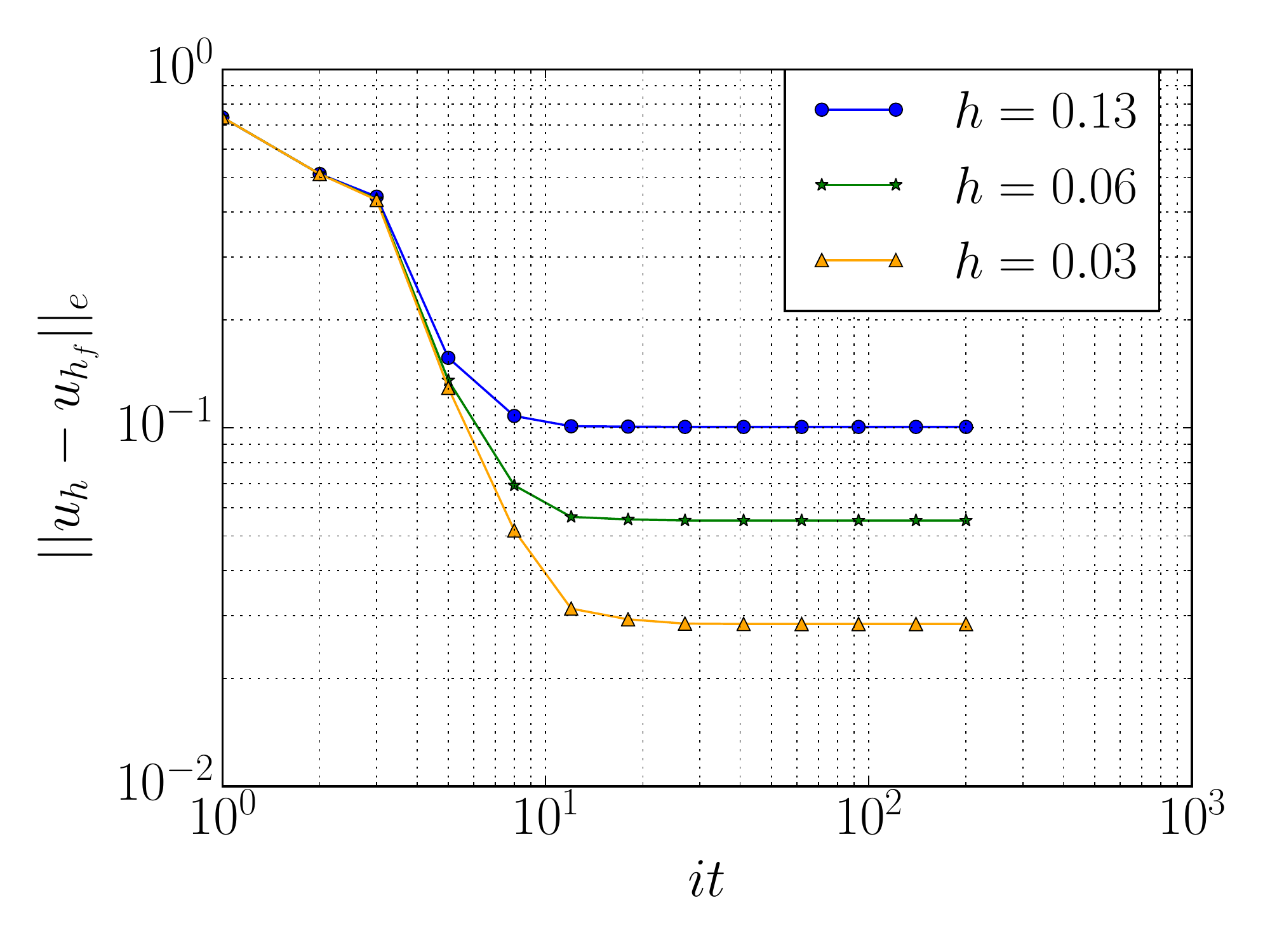}
\label{subfig: two inclusions error it coarse}}
\subfloat[Energy norm error with iteration count for finest solution.]{\includegraphics[width=.5\textwidth]{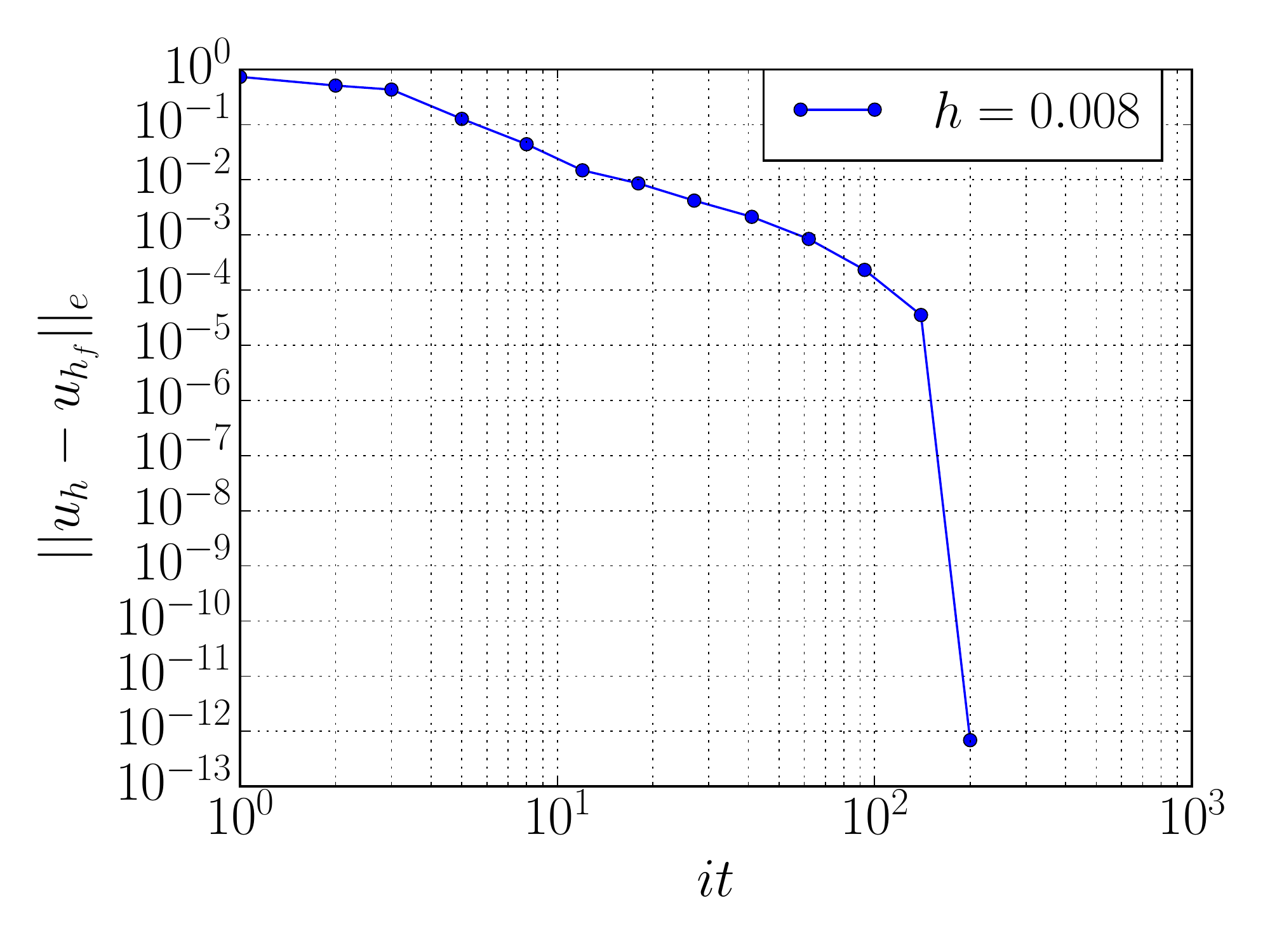}
\label{subfig: two inclusions error it finest}}
\caption{Convergence rates and energy norm error with iteration count for two inclusions with $E^1=E^2=E^3=1$.}
 \label{fig: two inclusions convergence rate}
\end{figure}

\begin{figure}
[htbp]
\centering
\subfloat[Energy norm error vs iterations and mesh size.]{\includegraphics[width=.6\textwidth]{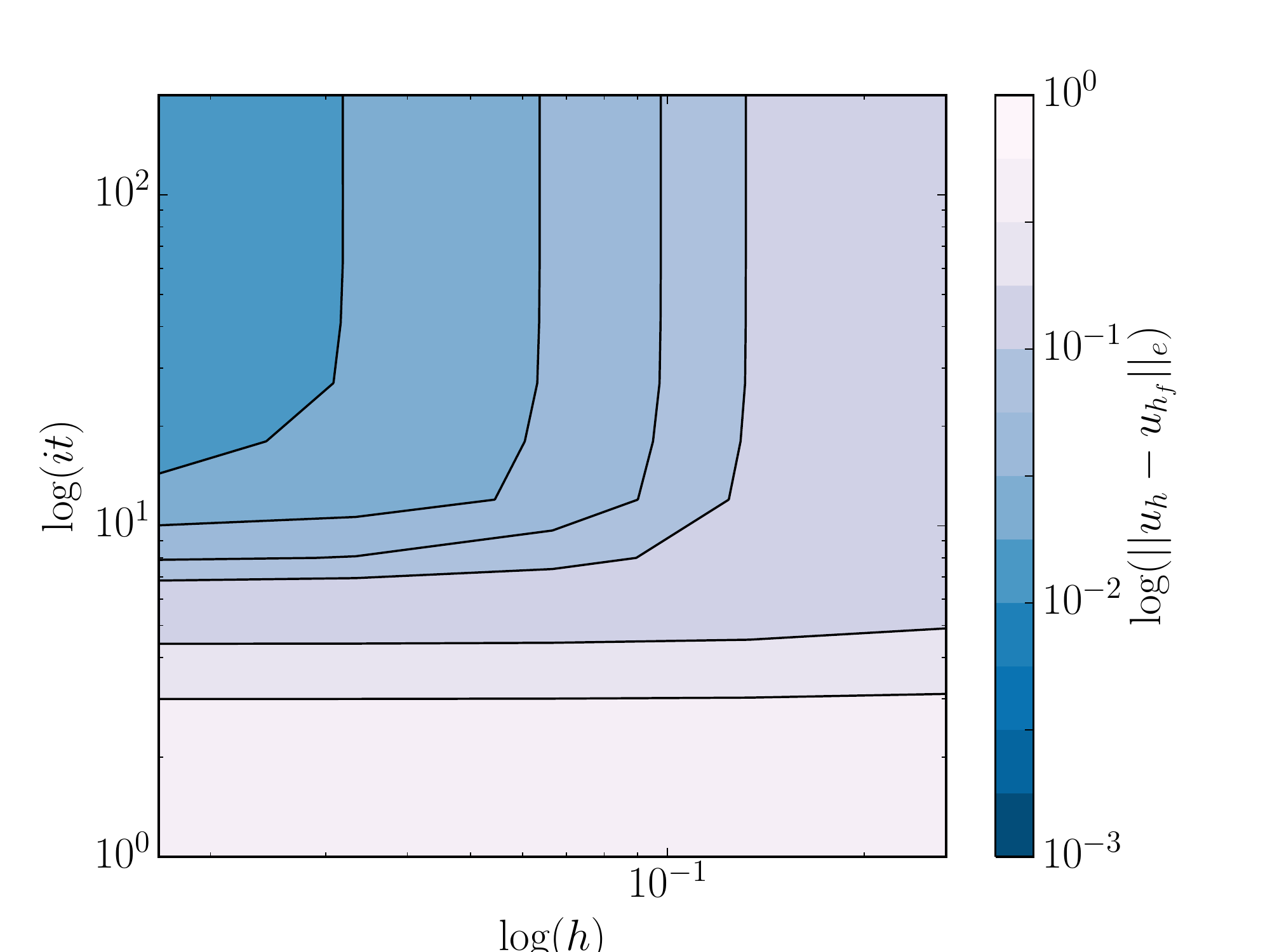}} 
\subfloat[Convergence rates.]{\includegraphics[width=.4\textwidth]{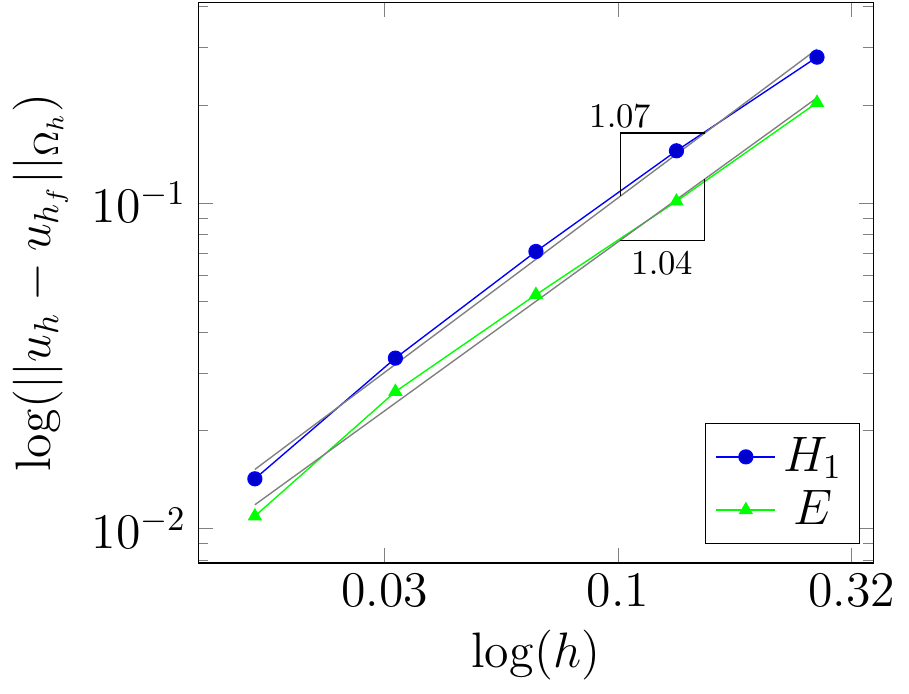}} \\
\caption{Convergence rates and energy norm error with iteration count for two inclusions with $E^2=E^3=10$ in a matrix with $E^1=1$.}
 \label{fig: two inclusions convergence rate 1_10}
\end{figure}

We choose a mesh configuration as displayed in Figure~\ref{subfig: schematic two inclusions mesh} and hierarchically refine the mesh as described in the previous example. We compare our coarse numerical solution with the finest numerical solution ($h=0.008$) at LaTin iteration $200$. The $H_1$ and energy norm errors with mesh refinement are shown in Figure~\ref{fig: two inclusions convergence rate} for a Young's modulus of $E^1=E^2=E^3=1$ in each domain and in Figure~\ref{fig: two inclusions convergence rate 1_10} for a contrast of the Young's modulus of $E^2=E^3=10$ for the inclusions to $E^1=1$ in the background domain. In both cases, with and without contrast in the Young's modulus between the inclusions and the background domain, we obtain convergence rates of first order which is optimal for our strategy, which relies on the P1 finite element method. 

%

\subsection{An advanced 3D geometry: braided composite material}

In this section, we consider a 3D braided composite geometry formed of interwoven threads. Each thread is created through the motion of an ellipsoid, given by 
\begin{equation}
\begin{aligned}
\phi^{thread}_j(x,y,z) &= \sqrt{\left(\frac{x-x_M^j}{a}\right)^2 + \left(\frac{y-y_M^j}{b}\right)^2 + \left(\frac{z-z_M^j}{c}\right)^2}   - r, \\
\end{aligned}
\end{equation}
along a sinusoidal path defined by $(x_M^j,y_M^j,z_M^j)$. For the ellipsoid, we choose $r = 1.0$, $a = 0.2$, $b = 0.2$ (width of the thread) and $c = 0.1$ (height of the thread). The three threads in  $x$-direction in our example geometry are described by $
x_M^j = x$, 
$y_M^j = \left[ - \frac{\lambda}{2}, 0, \frac{\lambda}{2} \right]$, 
$z_M^j = A \sin((\frac{2\pi(x + z_C^j)}{\lambda}))$, 
$z_C^j = \left[  \frac{\lambda}{2}, 0, \frac{\lambda}{2} \right]$. 
Here, $\lambda=1$ is the length of the path oscillation (sinus period), $A = 0.1$ is the amplitude of the path oscillation and $z_C^j$ is the phase shift. Each new thread is shifted in $y$-direction by half the period length and their phase is shifted by half the period length. We form a union of all threads in $x$-direction to obtain one level set function, \textit{i.e.} $\phi_x(x,y,z) = \min_j(\phi^{thread}_j(x,y,z))$.
The threads in the $y$-direction are combined into a second level set function analogously. They are shifted by $\frac{\lambda}{4}$ with respect to the threads in the $x$-direction to build the interwoven geometry.  

\begin{figure}
[htbp]
\centering
\subfloat[]{
\centering
\begin{tikzpicture}
\node[inner sep=0pt] (cont) at (0,0) {\includegraphics[width=.5\textwidth]{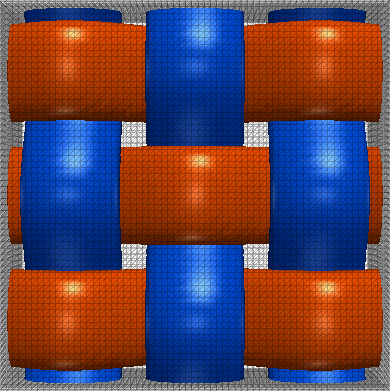}};
\node at (0,4.45) {{\Large $\Disp = \left(\begin{smallmatrix}  0 \\ 1 \\ 0 \end{smallmatrix}\right)$}};


\node at (0,-4.3) {{\Large $\Disp=\V{0}$}};
\node[rotate=90] at (4.3,0) {{\Large $\Stress \cdot \Normal = \V{0}$}};
\node[rotate=90] at (-4.3,0) {{\Large $\Stress \cdot \Normal = \V{0}$}};
\end{tikzpicture}
} \\
\subfloat[]{\includegraphics[width=.5\textwidth]{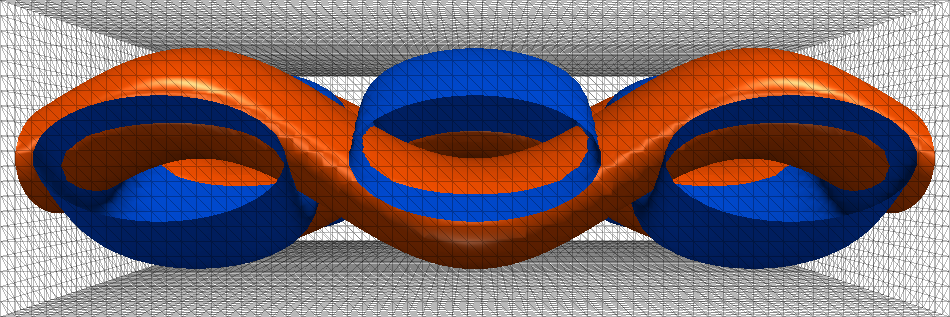}} 
\caption{Zero level set contours of braided composite with mesh and boundary conditions.}
 \label{fig: level set braided composite}
\end{figure}

We apply the boundary conditions displayed in Figure~\ref{fig: level set braided composite} to a domain $\Omega = [-\frac{1}{2},-\frac{3}{4},-\frac{1}{4}] \times [1,\frac{3}{4},\frac{1}{4}]$. We set a Young's modulus of one in all subdomains and pull the threads and the matrix in $y$-direction.  We choose $k^+=k^-=1.0$, $\gamma_g=0.1$,  $\gamma_{\Pi}=0.1$ and $\alpha=10$.

\begin{figure}
[htbp]
\subfloat[$\sigma_{yy}$]{
\centering
\begin{tikzpicture}
\node[inner sep=0pt] (cont) at (0,0) {\includegraphics[width=.4\textwidth]{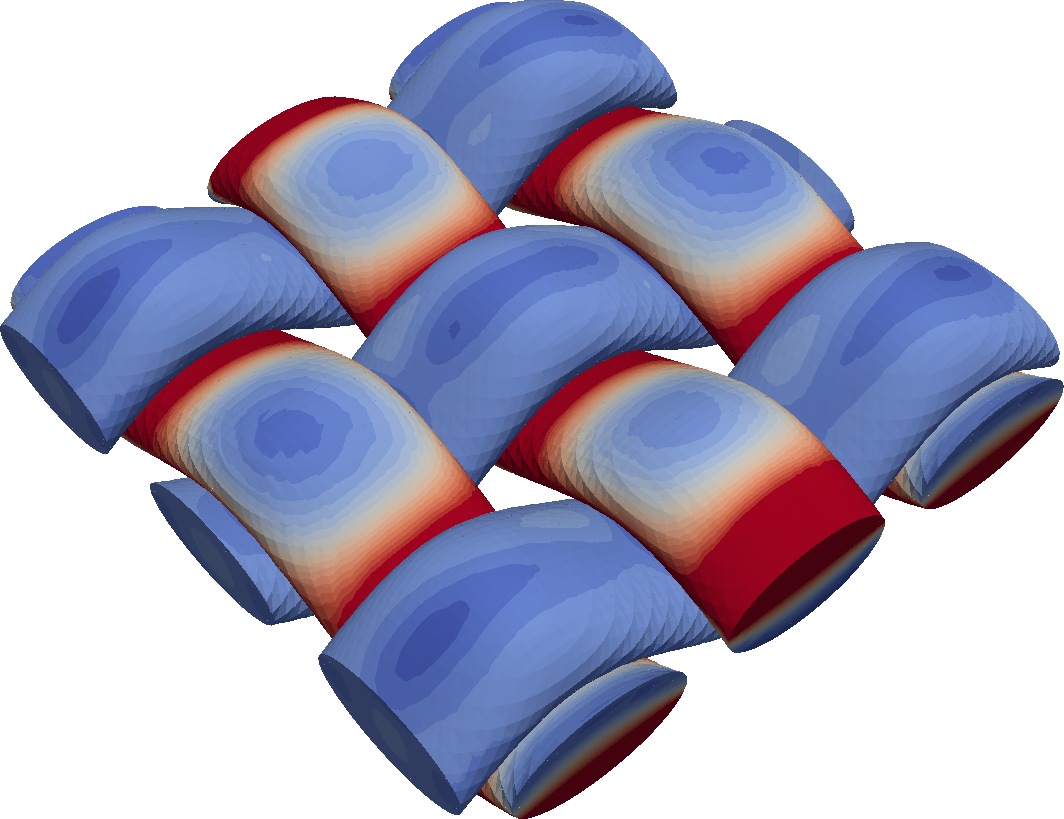}\,\,
\includegraphics[width=.1\textwidth]{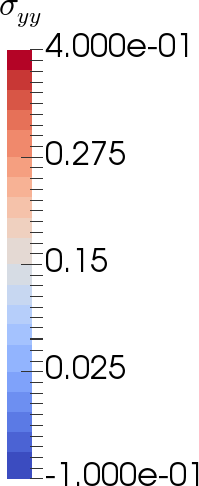}};
\draw[arrows=-latex] (-1,-3) -- (0.6,-1.7) node[below right] {$x$};
\draw[arrows=-latex] (-1,-3) -- (-2.4,-1.7) node[below left] {$y$};
\draw[arrows=-latex] (-1,-3) -- (-1,-2.5) node[right] {$z$};
\end{tikzpicture}
} 
\subfloat[$\sigma_{xy}$]
{\centering
\begin{tikzpicture}
\node[inner sep=0pt] (cont) at (0,0) {\includegraphics[width=.4\textwidth]{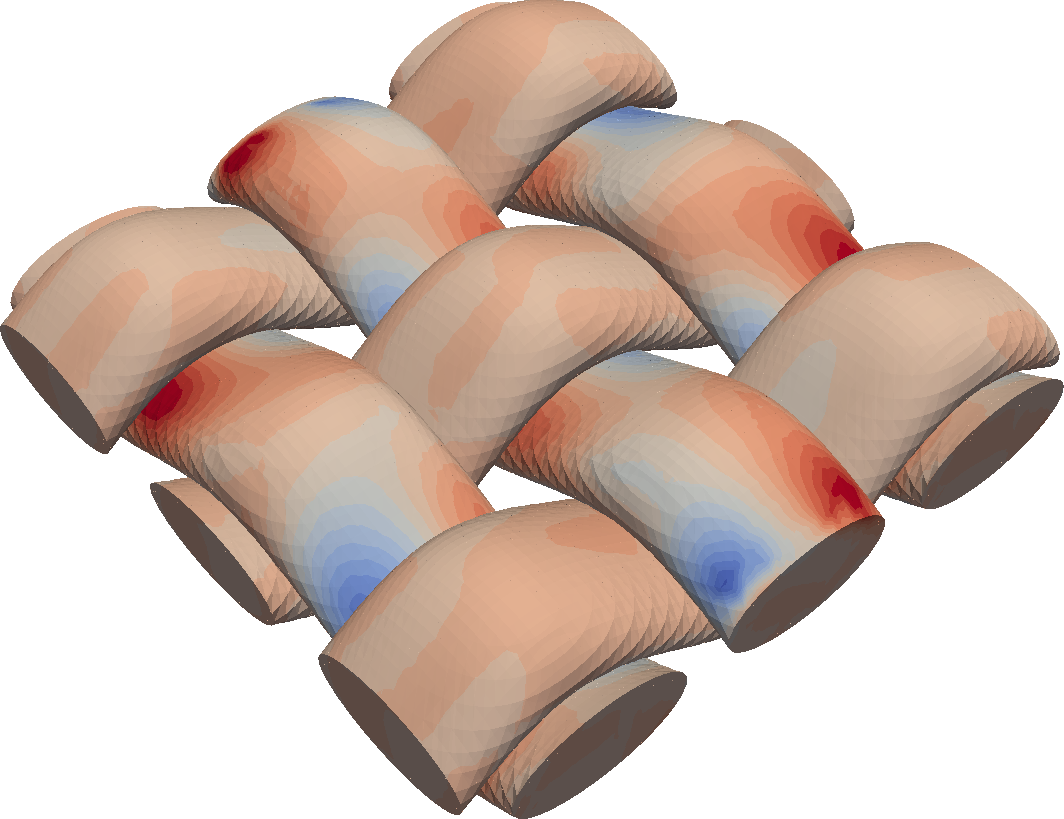}\,\,
\includegraphics[width=.1\textwidth]{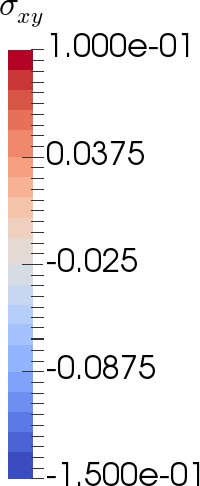}}; 
\draw[arrows=-latex] (-1,-3) -- (0.6,-1.7) node[below right] {$x$};
\draw[arrows=-latex] (-1,-3) -- (-2.4,-1.7) node[below left] {$y$};
\draw[arrows=-latex] (-1,-3) -- (-1,-2.5) node[right] {$z$};
\end{tikzpicture}
}
\caption{The normal stress component $\sigma_{yy}$ and the shear stress component $\sigma_{xy}$ of the braided composite.}
 \label{fig: stress components braided}
\end{figure}
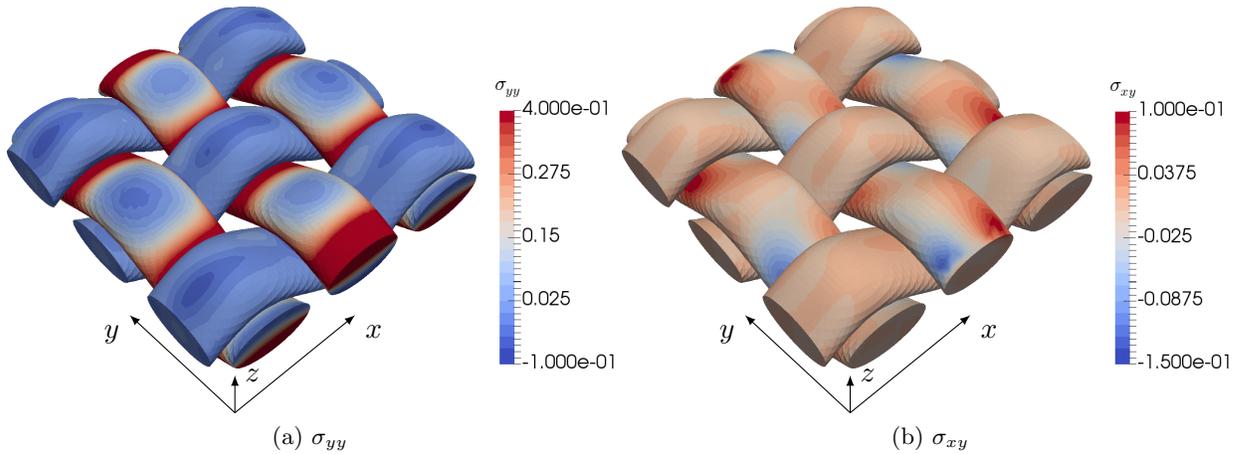

\begin{figure}
[htbp]
\centering
\subfloat[Diplacement in $y$-direction.]
{\centering
\begin{tikzpicture}
\node[inner sep=0pt] (cont) at (0,0) {\includegraphics[width=.4\textwidth]{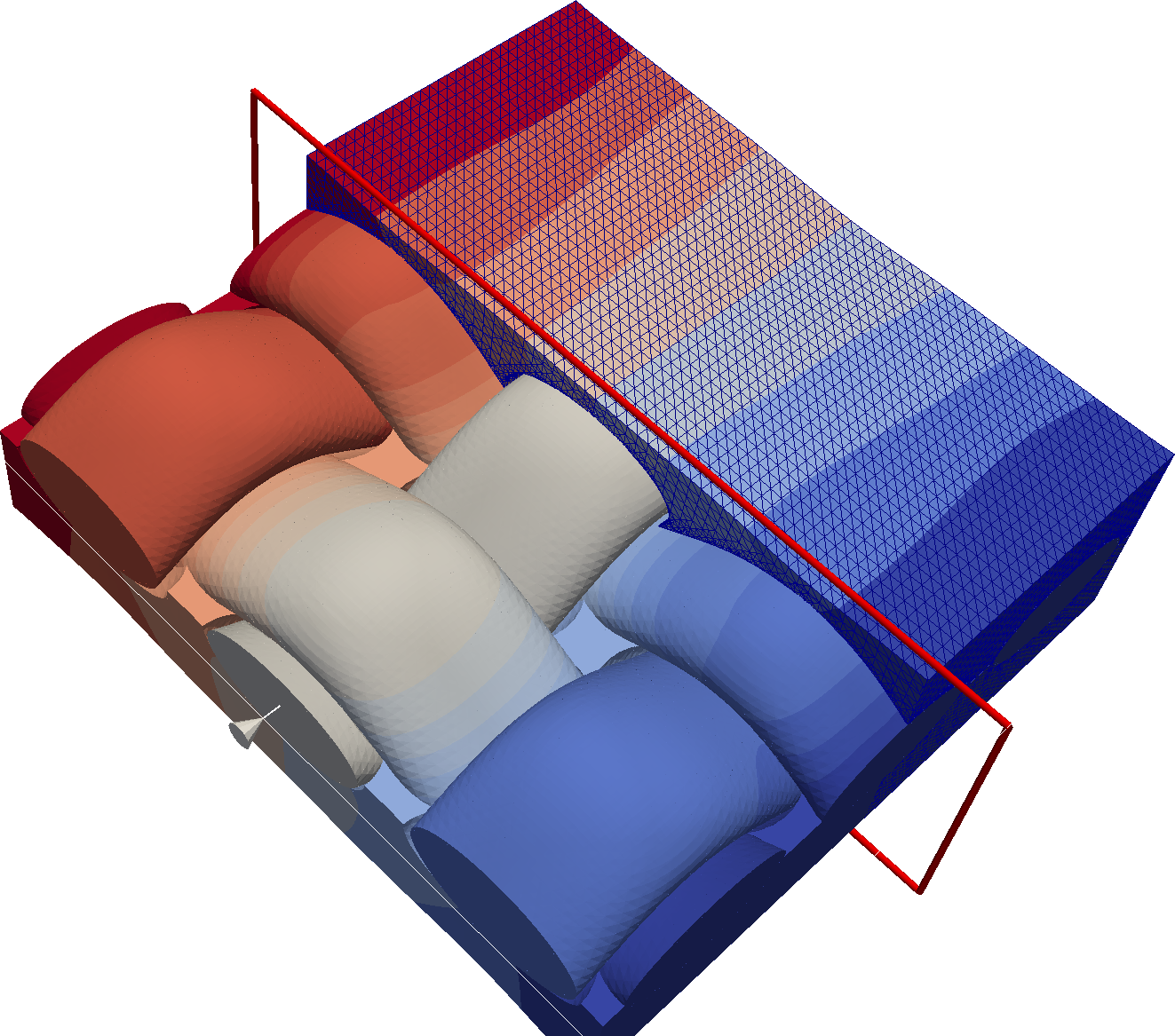}\,\,
\includegraphics[width=.1\textwidth]{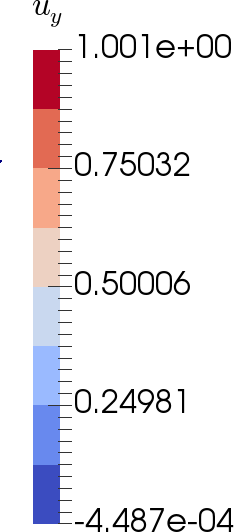}}; 
\begin{scope}[yshift=-.4cm]
\draw[arrows=-latex] (-1,-3) -- (0.6,-1.7) node[below right] {$x$};
\draw[arrows=-latex] (-1,-3) -- (-2.4,-1.7) node[below left] {$y$};
\draw[arrows=-latex] (-1,-3) -- (-1,-2.5) node[right] {$z$};
\end{scope}
\end{tikzpicture}
\label{subfig: displacement braided}
} \\
\subfloat[Stress component $\sigma_{zz}$.]{
\parbox{.8\textwidth}
{\centering
\begin{tikzpicture}
\node[inner sep=0pt] (cont) at (0,0) {\includegraphics[width=.8\textwidth]{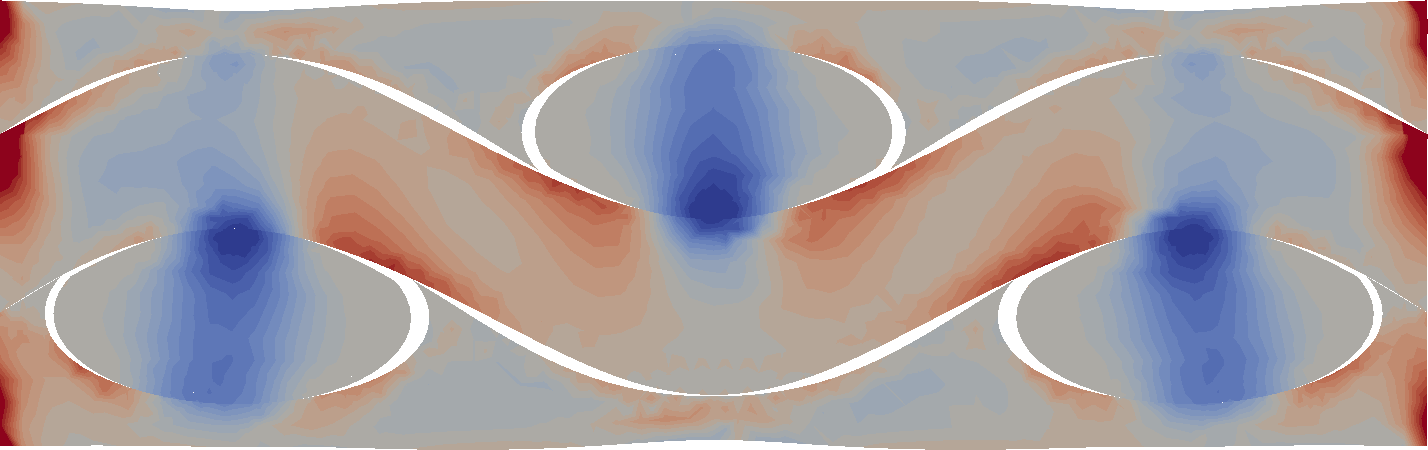}}; 
\begin{scope}[yshift=-.1cm]
\draw[arrows=-latex] (  (6.5,-2) -- (5.5,-2) node[below right] {$y$};
\draw[arrows=-latex] (6.5,-2) -- (6.5,-1) node[below right] {$z$};
\end{scope}
\end{tikzpicture}
\\
\includegraphics[width=.5\textwidth]{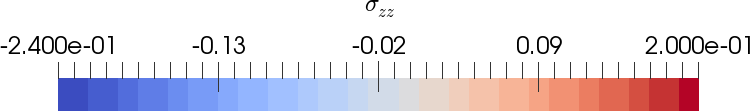}}
\label{subfig: stress szz braided}
}
\caption{Displacement in $y$-direction and normal stress component $\sigma_{zz}$ in a slice of the domain indicated in red in (a).}
 \label{fig: displacement and szz braided} 
\end{figure}

The resulting normal stress component $\sigma_{yy}$ and the shear stress component $\sigma_{xy}$ are displayed in Figure~\ref{fig: stress components braided}. Figure~\ref{fig: displacement and szz braided} shows the displacement in $y$-direction and the stress component $\sigma_{zz}$ in a slice of the domain indicated in red in Figure~\ref{subfig: displacement braided}.
We can qualitatively verify the numerical solution by inspecting the slice displayed in Figure~\ref{subfig: stress szz braided}. A vertical axial stress component (compression) is transferred through interfaces where the tangent plane is horizontal only when the vertical jump of displacement vanishes, which means that composite phases are in contact. Conversely, we see that when phases are not in contact at interface points where the tangent plane is horizontal, then the vertical compressive stress is zero.


Let us also point out that neither in 3D nor in 2D did we observe the development of numerical instabilities around contact zones, which we believe is due to the proposed P1-P1 mixed LaTIn formulation, together with the globalisation of the local stage and of course the ghost-penalty regularisation of the linear stages.



\section{Conclusions}


We have presented a novel LaTIn-CutFEM approach for the simulation of multiple-body contact mechanics. The solver uses the enrichment technique that is the main characteristic of the CutFEM technology: the introduction of overlapping meshes (\textit{i.e.} duplication of interface elements) to allow the representation of discontinuities in the solution fields across embedded interfaces. In order for the proposed scheme to be numerically stable, a P1/P1 reformulation of the LaTIn hybrid-mixed formulation has also been developed. Specifically, a weak enforcement of the contact conditions by a piecewise linear pair of interface fields prevents the growth and onset of spurious solution modes. Moreover, we have employed the ghost-penalty method, which prevents ill-conditioning of system matrices in the presence of ``bad cuts". We have demonstrated that our scheme is optimally convergent with mesh refinement, and that it is as stable as a conforming, primal finite element formulation.

We have also demonstrated the capabilities and notably the high geometrical flexibility of our implementation of the proposed algorithm. The latter property has been obtained through the use of a versatile multiple-level set approach to describe the geometry of complex interacting solids. In particular, we have shown that our algorithm is robust in the case where multiple simultaneous enrichments are performed, which happens when several solids interact at the same point.

This work is expected to be a strong basis for further contributions to the field. For now, our main focus is to develop a robust stopping criterion for LaTIn. Indeed, we have shown that, in contrast to widespread belief, the LaTIn algorithm is relatively fast. In our examples, where, arguably, the default search direction parameter is close to its optimum value, 20 iterations lead to solution states whereby the algorithmic error is blurred into the discretisation error. We hope to be able to make use of this observation to develop rational estimates stopping criteria for LaTIn formulations of complex interface laws.


\section*{Acknowledgments}

The authors gratefully acknowledge the financial support provided by the Welsh Government and Higher Education Funding Council for Wales through the S\^{e}r Cymru National Research Network in Advanced Engineering and Materials and EPSRC under grant EP/J01947X/1: \textit{Towards rationalised computational expense for simulating fracture over multiple scales} (RationalMSFrac). Susanne Claus would like to acknowledge the work of Andr\'{e} Massing in co-developing the CutFEM library.

\bibliographystyle{plain}
\bibliography{bibliography.bib}


\end{document}